\documentclass[final]{siamltex}

\usepackage{latexsym,amsmath,amssymb,amsfonts,amscd,multirow,dsfont}
\usepackage{epsfig,verbatim,epstopdf,graphics}
\usepackage{bm,color}
\usepackage{subfigure}

\newtheorem{thm}{Theorem}[section]

\newtheorem{exa}[thm]{Example}
\newtheorem{rem}[thm]{Remark}

\newcommand{\bk}{\mathbf{k}}
\newcommand{\bx}{\mathbf{x}}

\newcommand{\bE}{\mathbf{E}}

\newfont{\iams}{msbm9}

\newcommand{\commentbis}[1]{}
\newcommand{\be}{\begin{eqnarray}}
\newcommand{\ee}{\end{eqnarray}}
\newcommand{\beno}{\begin{eqnarray*}}
\newcommand{\eeno}{\end{eqnarray*}}

\newcommand{\barr}[1]{\begin{array}{#1}}
\newcommand{\earr}{\end{array}}

\newcommand{\beq}{\begin{equation}}
\newcommand{\eeq}{\end{equation}}
\newcommand{\beqa}{\begin{eqnarray}}
\newcommand{\eeqa}{\end{eqnarray}}

\newcommand{\bv}{{\bf v}}

\newcommand{\bV}{{\bf V}}
\newcommand{\bn}{{\bf n}}

\newcommand{\bzero}{\mathbf{0}}
\newcommand{\bone}{\mathbf{1}}
\newcommand{\bl}{\mathbf{l}}
\newcommand{\bi}{\mathbf{i}}

\newcommand{\bj}{\mathbf{j}}
\newcommand{\bW}{\mathbf{W}}

\newcommand{\ba}{{\bm{a}}}
\newcommand{\bal}{{\bm{\alpha}}}
\newcommand{\bb}{{\bm{\beta}}}
\newcommand{\bq}{\mathbf{q}}

\title
{An Adaptive Multiresoluton Discontinuous Galerkin Method for Time-Dependent Transport
Equations in Multi-dimensions}

\author{
 Wei Guo
\thanks{Department of Mathematics, Michigan State University,
East Lansing, MI 48824 U.S.A.
 {\tt wguo@math.msu.edu}}
\and
 Yingda Cheng
\thanks{Department of Mathematics, Michigan State University,
East Lansing, MI 48824 U.S.A.
 {\tt ycheng@math.msu.edu}. Research is supported by NSF grant  DMS-1453661.}
}

\date{\today}

\begin{document}

\maketitle

\begin{abstract}
In this paper, we develop  an adaptive multiresolution discontinuous Galerkin (DG) scheme for  time-dependent transport equations in multi-dimensions.  The method is constructed using   multiwavlelets  on tensorized nested grids. Adaptivity is realized by   error thresholding based on the hierarchical surplus, and the Runge-Kutta DG (RKDG) scheme is employed as the reference time evolution algorithm.
We show that the scheme performs similarly to a sparse grid DG method when the solution is smooth, reducing computational cost in multi-dimensions. When  the solution is no longer smooth, the adaptive algorithm can automatically capture fine local structures.   The method is therefore very suitable for deterministic kinetic simulations. Numerical results including several benchmark tests, the Vlasov-Poisson  (VP) and oscillatory VP systems are provided.
\end{abstract}


\section{Introduction}

In this paper, we propose  an adaptive multiresolution DG scheme for  time-dependent transport equations in multi-dimensions. This  is a continuation of our previous research on sparse grid DG schemes \cite{sparsedgelliptic,guo_sparsedg}. 
In particular, here we  consider   linear variable-coefficient equations,  aiming at  developing  efficient solvers for  kinetic transport problem   as the eventual goal.   
It is well known that the main bottleneck to solve kinetic equations are their high dimensionality. The equations are posed in the probability space, which is in six dimensions in a realistic setting. A popular framework for high dimensional computations is called sparse grid  \cite{zenger1991sparse,bungartz2004sparse, garcke2013sparse}. The idea is to use a properly truncated subset of the tensor product approximation space to break the curse of dimensionality. In our previous work \cite{guo_sparsedg}, a sparse grid DG method has been formulated and applied to kinetic simulations. The construction is based on Alpert's multiwavelets \cite{alpert1993class,alpert_adaptive_2002} and the method is demonstrated to  save significant computational and storage cost because of the reduced degrees of freedom of the approximation space.  By using the DG framework, many attractive features such as stability and conservation   can be proven.
However,   the scheme's success and the underlying convergence theory still rely heavily on the smoothness of the exact solution. In fact, it was generally understood that any \emph{a priori} type of choice of the sparse grid approximation space will depend on the smoothness assumption of the exact solution, which is often not satisfied in practice. 
For example, for the VP system and many other kinetic models,  small scale structures will often develop over time. Therefore, using the standard sparse grid methods or  any uniform grid based methods may not be optimal. 
The situation is   even worse if the solution contains discontinuities.
 In the literature, adaptive sparse grid methods have been developed \cite{zenger1991sparse,griebel1998adaptive,bungartz2004sparse,bokanowski_adaptive_2012} to address this issue. Such schemes measure the hierarchical coefficients or the so-called hierarchical surplus as a natural indicator for refinement or coarsening. There is a particular connection of this approach with the celebrated adaptive wavelet method \cite{dahmen_wavelet_1997,cohen2000wavelet}. This type of multiresolution schemes have been used to  accelerate the computations for conservation laws under finite difference or finite volume frameworks \cite{harten_multiresolution_1995,bihari_multiresolution_1997,dahmen_multiresolution_2001,alves_adaptive_2002,cohen_fully_2003,chiavassa2005multiresolution}. In recent years, there have been developments of adaptive multiresolution DG schemes  \cite{calle_wavelets_2005,archibald_adaptive_2011,hovhannisyan_adaptive_2014,gerhard_adaptive_2014,gerhard_high-order_2014,cite-key} which use  the multiwavelets of Alpert for computing conservation laws and compressible flows. In the context of adaptive computation for Vlasov equations, closely related work includes semi-Lagrangian type wavelet method \cite{besse2003adaptive,gutnic_vlasov_2004,besse_wavelet-mra-based_2008} and the $h$-adaptive RKDG method \cite{zhuqiu}. 
 
The objective of the present paper is to develop an adaptive multiresolution DG method that also fits under the sparse grid framework in multi-dimensions. When compared with other adaptive multiresolution DG methods in the literature, the main difference is in the multi-dimensional case. Our scheme will naturally go back to a sparse grid DG method, saving computational cost, when the solution possess sufficient smoothness. This is realized by using the fully tensorized basis functions instead of exploring multiwavelet only in local elements.  When the solution is no longer smooth, the adaptive algorithm that uses the  hierarchical surplus as the refinement or coarsening indicator, can automatically capture the local structures, thus removing the smoothness requirement of  \emph{a priori} chosen sparse grid approximation space. We use the hash table as the underlying data structure and can deal with equations in arbitrary dimensions. By using the DG formulation, many nice properties are retained for the transport equations. The numerical scheme is validated by benchmark tests with smooth and nonsmooth solutions,   the standard VP system and oscillatory VP system.

The rest of this paper is organized as follows:
in Section \ref{sec:method},  we  construct the adaptive multiresolution DG scheme.  
The numerical performance is validated in Section \ref{sec:numerical} by   three benchmark tests. Section \ref{sec:kinetic} discusses the application  to Vlasov equations, and  we conclude the paper   in Section  \ref{sec:conclusion}.

 %adaptive DG \cite{hartmann_adaptive_2002,hartmann_adaptive_2003}
%amr \cite{berger_adaptive_1984-1,berger_local_1989,plewa2003adaptive}

\section{Numerical method}
\label{sec:method}

	In this section, we formulate an adaptive multiresolution  DG method for solving time-dependent linear transport equations. First, we review 
the multiresolution analysis and multiwavelets which serve as foundations of the underlying scheme. Then, we discuss an adaptive multiresolution projection method that supplies the numerical initial conditions. The adaptive  time evolution algorithm is introduced at the end of this section after a review of the reference DG method.
		
	\subsection{Multiresolution analysis and multiwavelets}

	In this subsection, we  review multiresolution analysis associated with piecewise polynomials. We focus on box shaped domains in this paper. Without loss of generality, all the discussions in this section are for a unit sized box $\Omega=[0,1]^d$, where $d$ is the dimension of the problem.

	 First, we review the case when $d=1$. We define a set of nested grids, where the $n$-th level grid $\Omega_n$ consists of $2^n$ uniform cells $I_{n}^j=(2^{-n}j, 2^{-n}(j+1)]$, $j=0, \ldots, 2^n-1,$ for any $n \ge 0.$ For notational convenience, we also denote $I_{-1}=[0,1].$

The nested grids result in the nested piecewise polynomial spaces. In particular, let
$$V_n^k:=\{v: v \in P^k(I_{n}^j),\, \forall \,j=0, \ldots, 2^n-1\}$$
  be the usual piecewise polynomials of degree at most $k$ on the $n$-th level grid $\Omega_n$. Then, we have $$V_0^k \subset V_1^k \subset V_2^k \subset V_3^k \subset  \cdots$$
We can now define the multiwavelet subspace $W_n^k$, $n=1, 2, \ldots $ as the orthogonal complement of $V_{n-1}^k$ in $V_{n}^k$ with respect to the $L^2$ inner product on $[0,1]$, i.e.,
\begin{equation*}
V_{n-1}^k \oplus W_n^k=V_{n}^k, \quad W_n^k \perp V_{n-1}^k.
\end{equation*}
For notational convenience, we let
 $W_0^k:=V_0^k$, which is standard piecewise polynomial space of degree $k$ on $[0,1]$. Therefore, we have $V_n^k=\bigoplus_{0 \leq l \leq n} W_l^k$.

 Now we need to supply a set of orthonormal basis associated with the space $W_l^k$. The case of mesh level $l=0$ is trivial. We use the scaled Legendre polynomials and denote the basis by
$v^0_{i,0}(x),\quad i=1,\ldots,k+1.$
When $l>0$, the orthonormal bases in $W_l^k$ are presented in \cite{alpert1993class} and denoted by 
$$v^j_{i,l}(x),\quad i=1,\ldots,k+1,\quad j=0,\ldots,2^{l-1}-1.$$
The construction follows a repeated Gram-Schmidt process and the explicit expression of the multiwavelet basis functions are provided in \cite{alpert1993class}. 
Note that such multiwavelet bases retain the orthonormal property of wavelet bases for different mesh levels, i.e.,
\begin{equation}
\label{ortho1d}
\int_0^1 v^j_{i,l}(x)v^{j'}_{i',l'}(x)\,dx=\delta_{ii'}\delta_{ii'}\delta_{jj'}.
\end{equation}

Now we are ready to review the case when $d>1$. First we recall some basic notations about multi-indices. For a multi-index $\mathbf{\alpha}=(\alpha_1,\cdots,\alpha_d)\in\mathbb{N}_0^d$, where $\mathbb{N}_0$  denotes the set of nonnegative integers, the $l^1$ and $l^\infty$ norms are defined as 
$$
|\bal|_1:=\sum_{m=1}^d \alpha_m, \qquad   |\bal|_\infty:=\max_{1\leq m \leq d} \alpha_m.
$$
The component-wise arithmetic operations and relational operations are defined as
$$
\bal \cdot \bb :=(\alpha_1 \beta_1, \ldots, \alpha_d \beta_d), \qquad c \cdot \bal:=(c \alpha_1, \ldots, c \alpha_d), \qquad 2^\bal:=(2^{\alpha_1}, \ldots, 2^{\alpha_d}),
$$
$$
\bal \leq \bb \Leftrightarrow \alpha_m \leq \beta_m, \, \forall m,\quad
\bal<\bb \Leftrightarrow \bal \leq \bb \textrm{  and  } \bal \neq \bb.
$$

By making use of the multi-index notation, we denote by $\bl=(l_1,\cdots,l_d)\in\mathbb{N}_0^d$ the mesh level in a multivariate sense. We  define the tensor-product mesh grid $\Omega_\bl=\Omega_{l_1}\otimes\cdots\otimes\Omega_{l_d}$ and the corresponding mesh size $h_\bl=(h_{l_1},\cdots,h_{l_d}).$ Based on the grid $\Omega_\bl$, we denote by $I_\bl^\bj=\{\bx:x_m\in(h_mj_m,h_m(j_{m}+1)),m=1,\cdots,d\}$ an elementary cell, and 
$$\bV_\bl^k:=\{\bv: \bv(\bx) \in Q^k(I^{\bj}_{\bl}), \,\,  \bzero \leq \bj  \leq 2^{\bl}-\bone \}= V_{l_1,x_1}^k\times\cdots\times  V_{l_d,x_d}^k$$
the tensor-product piecewise polynomial space, where $Q^k(I^{\bj}_{\bl})$ denotes the collection of polynomials of degree up to $k$ in each dimension on cell $I^{\bj}_{\bl}$. 
If we use equal mesh refinement of size $h_N=2^{-N}$ in each coordinate direction, the  grid and space will be denoted by $\Omega_N$ and $\bV_N^k$, respectively.  

Based on a tensor-product construction, the multi-dimensional increment space can be  defined as
$$\bW_\bl^k=W_{l_1,x_1}^k\times\cdots\times  W_{l_d,x_d}^k.$$ 
Therefore,  the standard tensor-product piecewise polynomial space on $\Omega_N$ can be written as
\begin{equation}\label{eq:hiere_tp}
\bV_N^k=\bigoplus_{\substack{ |\bl|_\infty \leq N\\\bl \in \mathbb{N}_0^d}} \bW_\bl^k,
\end{equation}
while the sparse grid approximation space we used in \cite{sparsedgelliptic,guo_sparsedg} is  
\begin{equation}
\label{eq:hiere_sg}
\hat{\bV}_N^k=\bigoplus_{\substack{ |\bl|_1 \leq N\\\bl \in \mathbb{N}_0^d}}\bW_\bl^k \subset \bV_N^k.
\end{equation}
The dimension of $\hat{\bV}_N^k$  scales as $O((k+1)^d2^NN^{d-1})$ \cite{sparsedgelliptic}, which is significantly less than that of $\bV_N^k$ with exponential dependence on $Nd$. The approximation results for $\hat{\bV}_N^k$ are discussed in \cite{sparsedgelliptic,guo_sparsedg}, which has a stronger smoothness requirement than the traditional $\bV_N^k$ space. In this paper, we will not require the numerical solution to be in $\hat{\bV}_N^k$, but rather in $\bV_N^k$ and to be chosen adaptively.

Finally, we define the   basis functions in multi-dimensions as
\begin{equation}
\label{basis}
v^\bj_{\bi,\bl}(\bx)\doteq\prod_{m=1}^d v^{j_m}_{i_m,l_m}(x_m), %\quad i_m=1,\ldots,k+1,\,j_m=0,\ldots,\max(0,2^{l_m-1}-1).
\end{equation}
for $\bl \in \mathbb{N}_0^d, \bj \in B_\bl \doteq \{\bj\in\mathbb{N}_0^d: \,\mathbf{0}\leq\bj\leq\max(2^{\bl-\mathbf{1}}-\mathbf{1},\mathbf{0}) \}$ and $\mathbf{1}\leq\bi\leq \bk+\mathbf{1}.$ 
The orthonormality of the bases can be established by \eqref{ortho1d}. Furthermore, we note that the support of $v^\bj_{\bi,\bl}$ is      $I_{\bl-\mathbf{1}}^\bj.$
%By defining a shorthand notation 
%$B_\bl = \{\bj\in\mathbb{N}_0^d: \,\mathbf{0}\leq\bj\leq\max(2^{\bl-\mathbf{1}}-\mathbf{1},\mathbf{0}) \}$,
%therefore,  a function in the sparse approximation $\bV_N^k$ can be written as 
%$$u_h(\bx)=\sum_{\substack{ |\bl|_\infty\leq N,\,\bj\in B_\bl\\ 
%%		\mathbf{0}\leq\bj\leq\max(2^{\bl-\mathbf{1}}-\mathbf{1},\mathbf{0})\\ 
%		\mathbf{1}\leq\bi\leq \bk+\mathbf{1}}} u^\bj_{\bi,\bl}v^\bj_{\bi,\bl}(\bx),$$
%where $u^\bj_{\bi,\bl}$ denotes the corresponding degree of freedom, which is also called the surplus in the literature. Similarly, a function in the sparse approximation space $\hat{\bV}_N^k$ can be written as 
%$$u^s_h(\bx)=\sum_{\substack{ |\bl|_1\in N,\,\bj\in B_\bl\\ 
%		%		\mathbf{0}\leq\bj\leq\max(2^{\bl-\mathbf{1}}-\mathbf{1},\mathbf{0})\\ 
%		\mathbf{1}\leq\bi\leq \bk+\mathbf{1}}} u^\bj_{\bi,\bl}v^\bj_{\bi,\bl}(\bx).$$
%

	\subsection{Adaptive multiresolution projection method}
	In this subsection, we formulate an adaptive multiresolution projection algorithm which supplies the numerical initial condition for DG schemes. Given a maximum mesh level $N$ and an accuracy threshold $\varepsilon>0$, we find a projected solution $u_h(\bx) \in \bV^k_N$ of a given function $u(\bx)$ defined on $\Omega$ using an adaptive procedure.  
	
The backbone of the  algorithm is the fact that  each hierarchical basis of space $\bV^k_N$ represents the fine level detail on a specific mesh scale, which naturally provides an error indicator for the design of adaptive algorithms. We first review the mixed derivative norm for a function $u(x).$  For any set $L=\{i_1, \ldots i_r \} \subset \{1, \ldots d\}$, we define $L^c$ to be the complement set of $L$ in $\{1, \ldots d\}.$ For a non-negative integer $\alpha$ and set $L$,  we define the semi-norm on any domain denoted by $\Omega$
$
|u|_{H^{\alpha,L}(\Omega)} :=  \left \| \left ( \frac{\partial^{\alpha}}{\partial x_{i_1}^{\alpha}} \cdots \frac{\partial^{\alpha}}{\partial x_{i_r}^{\alpha}}  \right ) u \right \|_{L^2(\Omega)}
$
and 
$
|u|_{\mathcal{H}^{q+1}(\Omega)} :=\max_{1 \leq r \leq d} \left ( \max_{\substack{L\subset\{1,2,\cdots,d\} \\|L|=r}} |u|_{H^{t+1, L}(\Omega)} \right ),$
which is the norm for the mixed derivative of $u$ of at most degree $q+1$ in each direction. For a function $u(\bx) \in \mathcal{H}^{p+1}(\Omega),$ we showed that \cite{guo_sparsedg} $u(\bx)=\sum_{\bl \in \mathbb{N}_0^d} \sum_{\bj \in B_\bl, \mathbf{1}\leq\bi\leq \bk+\mathbf{1}} u^\bj_{\bi,\bl} v^\bj_{\bi,\bl}(\bx),$ and
$$
\left (\sum_{\bj \in B_\bl, \mathbf{1}\leq\bi\leq \bk+\mathbf{1}} |u^\bj_{\bi,\bl}|^2 \right )^{1/2} \leq C2^{-(q+1)|\bl|_1}|u|_{\mathcal{H}^{q+1}(\Omega)},
$$
where  $u^\bj_{\bi,\bl}= \int_{\Omega}u(\bx)v^\bj_{\bi,\bl}(\bx)d\bx,$  $q=\min\{p,k\},$ and $C$ is a constant independent of mesh level $\bl.$ Henceforth, the hierarchical coefficient $u^\bj_{\bi,\bl}$ (also called hierarchical surplus) serves as a natural indicator for the local smoothness of $u(\bx)$. The main idea of the adaptive algorithm is to choose only coefficients above a prescribed threshold value $\varepsilon$. In this paper, we experiment on    error indicators  
$\left\|\sum_{\mathbf{1}\leq\bi\leq\bk+\mathbf{1}}u^\bj_{\bi,\bl}v^\bj_{\bi,\bl}(\bx)\right\|_{L^s(\Omega)}$ using different norms,
where $|| \cdot ||_{L^s(\Omega)}$ denotes the broken Sobolev $L^s(\Omega)$ norm for a function in $\bV_N^k$, with $s=1, 2,  \infty.$ When $s=2,$ due to orthonormality of the basis, $\left\|\sum_{\mathbf{1}\leq\bi\leq\bk+\mathbf{1}}u^\bj_{\bi,\bl}v^\bj_{\bi,\bl}(\bx)\right\|_{L^2(\Omega)}$ is equivalent to $\left(\sum_{\mathbf{1}\leq\bi\leq\bk+\mathbf{1}}|u^\bj_{\bi,\bl}|^2\right)^{\frac12}.$   In other cases, for simplicity, we use instead $\sum_{\mathbf{1}\leq\bi\leq\bk+\mathbf{1}}|u^\bj_{\bi,\bl}|\|v^\bj_{\bi,\bl}(\bx)\|_{L^1(\Omega)}$ for $s=1$ and $\sum_{\mathbf{1}\leq\bi\leq\bk+\mathbf{1}}|u^\bj_{\bi,\bl}|\|v^\bj_{\bi,\bl}(\bx)\|_{L^\infty(\Omega)}$ for $s=\infty$. The values of $|\|v^\bj_{\bi,\bl}(\bx)\|_{L^1(\Omega)}$ and $|\|v^\bj_{\bi,\bl}(\bx)\|_{L^\infty(\Omega)}$ can be precomputed and stored. Overall,  $\|v^\bj_{\bi,\bl}(\bx)\|_{L^1(\Omega)}$ scales as $2^{-|\bl|_1/2}$, and $\|v^\bj_{\bi,\bl}(\bx)\|_{L^\infty(\Omega)}$ scales as $2^{|\bl|_1/2}.$ 

In summary, we flag an element $V^\bj_{\bl}:=\{v^\bj_{\bi,\bl}, \mathbf{1}\leq\bi\leq\bk+\mathbf{1}\}$  if
	\begin{align}
	&\sum_{\mathbf{1}\leq\bi\leq\bk+\mathbf{1}}|u^\bj_{\bi,\bl}|\|v^\bj_{\bi,\bl}(\bx)\|_{L^1(\Omega)}>\varepsilon,\quad\text{if}\quad s=1\label{eq:l1}\\
	&\left(\sum_{\mathbf{1}\leq\bi\leq\bk+\mathbf{1}}|u^\bj_{\bi,\bl}|^2\right)^{\frac12}>\varepsilon,\quad\text{if}\quad s=2\label{eq:l2}\\
	&\sum_{\mathbf{1}\leq\bi\leq\bk+\mathbf{1}}|u^\bj_{\bi,\bl}|\|v^\bj_{\bi,\bl}(\bx)\|_{L^\infty(\Omega)}>\varepsilon,\quad\text{if}\quad s=\infty\label{eq:l8},
	\end{align}
	where $\varepsilon$ is a prescribed error threshold.  Similar to \cite{griebel1998adaptive}, we use a top down approach,  starting recursively from the coarsest level. 	Once an element is flagged, then we consider adding its children elements for improvement of accuracy.
	In particular, if a element $V^{\bj'}_{\bl'}$ satisfies the following conditions:
	\begin{itemize}
		\item There exists an integer m such that $1\le m\le d$ and $\bl'=\bl + \mathbf{e}_m$, where $\mathbf{e}_m$ denotes the unit vector in $x_m$ direction, and the support of $V^{\bj'}_{\bl'}$ is within the support of $V^{\bj}_{\bl}.$%$\bj'=\bj+j_m \mathbf{e}_m + r,$ where $j_m$ is $m$-th component of $\bj$ and $r=0$ or 1.
		\item $|\bl'|_\infty\leq N$,	
	\end{itemize} 
	then it is called a child element of $V_\bl^\bj$. Accordingly, element $V_\bl^\bj$ is called a parent element of $V_{\bl'}^{\bj'}$. In this notation, we can see an element can have multiple children and multiple parents.

	The last component of the algorithm is  an efficient data structure. As suggested in  \cite{griebel1998adaptive}, we use
 the hash table approach which is easy to implement, requires little storage overhead, and allows one to conveniently deal with hierarchical index $(\bl,\bj)$ in the implementation. Specifically, by a prescribed hash-function, a hierarchical index $(\bl,\bj)$ is mapped to  a hash-key (an integer), which serves as an address in the hash table. Then, given a hierarchical index, the associated data can be easily stored and retrieved by computing the hash-key. For more details about the hash table including how to choose proper hash-function and other implementation details, readers are referred to \cite{griebel1998adaptive}.
	
	Finally, we summarize the adaptive projection algorithm as follows. 

\medskip	
	
\noindent\rule{16.5cm}{1pt}\\
\noindent{\bf Algorithm 1: Adaptive projection}\\
\noindent\rule{16.5cm}{0.4pt}
 
 {\bf Input:} Function $u(\bx)$.
 
  {\bf Parameters:} Maximum level $N,$ polynomial degree $k,$  error threshold   $\varepsilon.$
 
  {\bf Output:} Hash table H, leaf table L and projected  solution $u_h(\bx) \in \bV_{N,H}^k.$
 
 \medskip
	\begin{enumerate}
		\item Project $u(\bx)$ onto the coarsest level of mesh, e.g., level 0. Add all elements to the hash table $H$ (active list). Define an element without children as a leaf element,
		and add all the leaf elements to the leaf table $L$ (a smaller hash table).
		\item For each leaf element $V_\bl^\bj$ in the leaf table, if \eqref{eq:l1}, \eqref{eq:l2} or \eqref{eq:l8} holds,
		 then we consider its child  elements: for a child element $V_{\bl'}^{\bj'}$, if it has not been added to the table $H$, then  compute the detail coefficients $ \{u^{\bj'}_{\bi,\bl'}, \mathbf{1}\leq\bi\leq\bk+\mathbf{1}\}$  and add $V_{\bl'}^{\bj'}$ to both table $H$ and table $L$.  For its parent elements in $H$, we increase the number of children by one. 
		\item Remove the parent elements from table $L$ for all the newly added elements.
		\item Repeat step 2 - step 3, until no element can be further added.
	\end{enumerate}
\noindent\rule{16.5cm}{1pt}

\bigskip
 
Once the adaptive projection algorithm completes, it will generate a final hash table H and a numerical approximation $u_h(\bx)=\sum_{v^\bj_{\bi,\bl} \in H}u^\bj_{\bi,\bl} v^\bj_{\bi,\bl}(\bx)$.  We denote the approximation space  $\bV^k_{N,H}=\textrm{span}\{v^\bj_{\bi,\bl} \in H\}$ and it is a subspace of $\bV^k_N.$
As noticed in \cite{griebel1998adaptive}, this top down approach may terminates too early and does not resolve the large coefficients on the very fine mesh levels. An alternative way is to find the $L^2$ projection of $u(\bx)$ in the finest level $\bV^k_N,$ and then truncate the elements with small coefficients as done in \cite{hovhannisyan_adaptive_2014}. However, this will effectively increase the computational cost and we do not pursue it in this work.

\subsection{The reference DG scheme}

In this subsection, we   review the standard RKDG method defined on space ${\bV}_N^k$ for
the following $d$-dimensional linear transport equation with variable coefficients  
\begin{equation}
\label{eq:model}
\left\{\begin{array}{l}
u_t + \nabla\cdot(\ba(t,\bx) \,u) =0,\quad \bx\in\Omega,\\[2mm]
u(0,\bx) = u_0(\bx),
\end{array}\right.
\end{equation}
subject to periodic boundary conditions. Other types of boundary conditions can be accommodated in a similar way.

 First, we review some basic notations about jumps and averages for   piecewise functions defined on the grid $\Omega_N$. Let $T_h$ be the collection of all elementary cell $I^{\bj}_{N}, \quad 0 \leq j_m  \leq 2^{N}-1, \forall \,m=1, \ldots, d$.
 $\Gamma:=\bigcup_{T \in \Omega_N} \partial_T$ be the union of the interfaces for all the elements in $\Omega_N$ (here we have taken into account the periodic boundary condition when defining $\Gamma$) and $S(\Gamma):=\Pi_{T\in \Omega_N} L^2(\partial T)$ be the set of $L^2$ functions defined on $\Gamma$. For any $q \in S(\Gamma)$ and $\bq \in [S(\Gamma)]^d$,  we define their   averages $\{q\}, \{\bq\}$ and jumps $[q], [\bq]$ on the interior edges as follows. Suppose
$e$ is an interior edge shared by elements $T_+$ and $T_-$, we define the unit normal vectors $\bm{n}^+$ and  $\bm{n}^-$ on $e$ pointing exterior of $T_+$ and $T_-$, respectively, then
\begin{flalign*}
[ q] \  =\  \, q^- \bm{n}^- \, +  q^+ \bm{n}^+, & \quad \{q\} = \frac{1}{2}( q^- + 	q^+), \\
[ \bq] \  =\  \, \bq^- \cdot \bm{n}^- \, +  \bq^+ \cdot \bm{n}^+, & \quad \{\bq\} = \frac{1}{2}( \bq^- + \bq^+).
\end{flalign*}

The semi-discrete DG formulation for   \eqref{eq:model} is defined as follows: find $u_h\in {\bV}_N^k$, such that
\begin{align}  
\label{eq:DGformulation}
\int_{\Omega}(u_h)_t\,v_h\,d\bx =& \int_{\Omega} u_h\ba\cdot\nabla v_h\,d\bx - \sum_{\substack{e \in \Gamma}}\int_{e} \widehat{\ba u_h} \cdot [v_h]\,ds,\quad  	\\
:= & A(u_h,v_h) \notag
\end{align}
for $\forall \,v_h \in {\bV}_N^k,$ where $\widehat{\ba u_h}$ is defined on the element interface denotes a monotone numerical flux to ensure the $L^2$ stability of the scheme. In this paper, we use the upwind flux
\begin{equation}
\widehat{\ba u_h} = \ba\{u_h\} + \frac{|\ba\cdot \bn|}{2}[u_h],
\end{equation}
with $\bn = \bn^+$ or $\bn^-$ for the constant coefficient case.
More generally, for variable coefficients problems, we adopt the global Lax-Friedrichs flux
\begin{equation}
	\widehat{\ba u_h} = \{\ba u_h\} + \frac{\alpha}{2}[u_h],
\end{equation}
where $\alpha=\max_{\bx}{|\ba(\bx,t)\cdot \bn|}$, the maximum is taken for all possible $\bx$ at time $t$ in the computational domain. 

We use the  total variation diminishing (TVD) Runge-Kutta methods \cite{Shu_1988_JCP_NonOscill} to solve the ordinary differential  equations resulting from the semidiscrete formulation \eqref{eq:DGformulation}, $(u_h)_t = R(u_h).$ A commonly used third-order TVD Runge-Kutta method is given by
\begin{align}
u_h^{(1)} &= u_h^{n} + \Delta t R(u^n_h), \notag \\
u_h^{(2)} &= \frac{3}{4}u_h^{n} + \frac14 u_h^{(1)} +\frac14 \Delta t R(u_h^{(1)}),\label{eq:tvd}\\
u_h^{n+1} &= \frac{1}{3}u_h^{n} + \frac23 u_h^{(1)} +\frac23 \Delta t R(u_h^{(2)}), \notag
\end{align}
where $u_h^{n}$ denotes the numerical solution at time level $t=t^n$. 

%Finally, we would like to make some remarks on the implementation issues. Unlike the traditional piecewise polynomial space, for which one element can only interact with itself and its immediate neighbors, the basis functions in the sparse space $\hat{\bV}^k_N$ are no longer locally defined due to the hierarchical structure, leading to additional challenges in implementation. In fact, it is crucial to  take full advantage of such a hierarchical (tree-like) structure when implementing the scheme to save computational cost. 
%As for the numerical flux, the global Lax-Friedrichs flux is adopted since we are able efficiently compute the interface integral in \eqref{eq:DGformulation} by using the \emph {unidirectional principle}. Such an idea  has been used in the sparse IPDG method for solving variable coefficient elliptic problems. In particular,  we first project $\ba$
%into space $\hat{\bV}^k_N$ and denote the resulting projection by $\ba_h$. Since $\ba_h$ is a separable function, the multi-dimensional interface integral in \eqref{eq:DGformulation} can be computed by evaluating multiplication of one-dimensional integrals. An advantage of this procedure is that we do not rely on numerical quadratures to compute the interface integrals, which can become quite complicated in the sparse grid setting. 

\subsection{Adaptive multiresolution DG evolution algorithm}
Based on the previous subsections, we are now ready to formulate the adaptive multiresolution DG evolution algorithm which consists of several key steps.

The first step is the prediction step, which means
given the hash table $H$ that stores the numerical solution $u_h$ at time step $t^n$ and the associated leaf table $L$, we need to predict the location where the details becomes significant at the next time step $t^{n+1}$, then add more elements in order  to capture the fine structures.  The time step size $\Delta t$ is chosen as follows. We denote by
	$l_m^n$ the largest mesh level in the $x_m$ direction in the current hash table $H$, and $l_m^{n,p} = \min(l_m^n+1, N)$ for the sake of possible refinement after prediction. Accordingly, we denote  $h_m^{n,p}=2^{-l_m^{n,p}}$. The time step $\Delta t$ for at time $t^n$ is given by
	\begin{align}
	\displaystyle\Delta t &= \frac{\text{CFL}}{\displaystyle\sum_{m=1}^d \frac{c_m}{h_m^{n,p}}},
	\end{align}
	where $c_m$ is the maximum wave propagation speed in $x_m$-direction and we use $\text{CFL}=0.1$ in our simulation.
We then solve for $u_h\in \bV_{N,H}^k$ from $t^n$ to $t^{n+1}$, such that
$
\int_{\Omega}(u_h)_t\,v_h\,d\bx =A(u_h,v_h)  
$
for $\forall \,v_h \in \bV_{N,H}^k,$ where $A(u_h,v_h)$ has been defined in \eqref{eq:DGformulation}.
The forward Euler discretization is used as the time integrator in this step and we denote the predicted solution at $t^{n+1}$  by $u_h^{(p)}.$  We remark that the standard global time stepping method is employed in the current adaptive framework for simplicity. It is nontrivial to develop a local time stepping method for the proposed adaptive multiresolution DG method due to the distinct hierarchical basis functions, and this subject is left for future study.
%Combined with the Forward Euler time stepping, we are able to obtain the prediction solution at time step $t^{n+1}$ as follows.
%\begin{equation}
%\label{eq:forwardeuler}
%u_h^{(p)} = u^{n} + \Delta t R(u^n_h),
%\end{equation}
%where $R$ denotes the DG spatial discretization operator resulting from \eqref{eq:DGformulation1}, and $u_h^{(p)}$ denotes the prediction solution at $t^{n+1}$. 

The second step is the refinement step according to $u_h^{(p)}$. We traverse the hash table $H$ and if an element $V_\bl^\bj$ satisfies the refinement criteria \eqref{eq:l1}, \eqref{eq:l2} or \eqref{eq:l8}, indicating that such an element becomes significant at the next time step, then we need to refine the mesh by adding its children elements to $H$. The detailed procedure is described as follows. For a child element $V_{\bl'}^{\bj'}$ of $V_\bl^\bj$, if it has been already added to $H$, i.e. $V_{\bl'}^{\bj'}\in H$, we do nothing; if not, we add the element $V_{\bl'}^{\bj'}$ to $H$ and let the associated detail coefficients $u^{\bj'}_{\bi,\bl'}=0,\,\mathbf{1}\leq\bi\leq\bk+\mathbf{1}$. Moreover, we need to increase the number of children by one for all elements that has $V_{\bl'}^{\bj'}$ as its child element and remove the parent elements of $V_{\bl'}^{\bj'}$ from the leaf table   if they have been added. Finally, we obtain a larger hash table $H^{(p)}$ and the associated approximation space $\bV_{N,H^{(p)}}^k$ and the leaf table $L^{(p)}$.

Then,  based on the updated hash table $H^{(p)}$, we   evolve the numerical solution by the DG formulation with space $\bV_{N,H^{(p)}}^k$. Namely, we solve for  $\bV_{N,H^{(p)}}^k$ from $t^n$ to $t^{n+1}$, such that
$
\int_{\Omega}(u_h)_t\,v_h\,d\bx =A(u_h,v_h)  
$
for $\forall \,v_h \in \bV_{N,H^{(p)}}^k,$ where $A(u_h,v_h)$ has been defined in \eqref{eq:DGformulation}. The semidiscrete equation is solved by  the TVD-RK scheme \eqref{eq:tvd} to generate the pre-coarsened numerical solution $\tilde{u}_h^{n+1}$. We notice that the first inner stage of the  Runge-Kutta method is actually the forward Euler prediction step. Moreover, recall that the detail coefficients for the newly added elements are set to zero. Therefore, after the time evolution of the first inner stage, the coefficients for original elements for $u_h^{(1)}$ should be the same as the prediction solution $u_h^{(p)}$, which can be reused to save computational cost. We only need to calculate the coefficients of newly added elements for $u_h^{(1)}$.

The last step   is to coarsen the mesh by removing elements that become insignificant at time level $t^{n+1}.$  The hash table $H^{(p)}$ that stores the numerical solution $\tilde{u}_h^{n+1}$  is recursively coarsened by the following procedure. 
The leaf table $L^{(p)}$ is traversed, and if an element $V_\bl^\bj\in L^{(p)}$ satisfies the coarsening criterion
\begin{align}
	&\sum_{\mathbf{1}\leq\bi\leq\bk+\mathbf{1}}|u^\bj_{\bi,\bl}|\|v^\bj_{\bi,\bl}(\bx)\|_{L^1(\Omega)}<\eta,\quad\text{if}\quad s=1\label{eq:l1_c}\\
	&\left(\sum_{\mathbf{1}\leq\bi\leq\bk+\mathbf{1}}|u^\bj_{\bi,\bl}|^2\right)^{\frac12}<\eta,\quad\text{if}\quad s=2\label{eq:l2_c}\\
	&\sum_{\mathbf{1}\leq\bi\leq\bk+\mathbf{1}}|u^\bj_{\bi,\bl}|\|v^\bj_{\bi,\bl}(\bx)\|_{L^\infty(\Omega)}<\eta,\quad\text{if}\quad s=\infty\label{eq:l8_c},
\end{align}
where $\eta$ is a prescribed error constant, then we remove the element from both table $L^{(p)}$ and $H^{(p)}$, and set the associated coefficients $u^{\bj'}_{\bi,\bl'}=0,\,\mathbf{1}\leq\bi\leq\bk+\mathbf{1}$. For each of its parent elements in table $H^{(p)}$, we decrease the number of children by one. If the number becomes zero, i.e, the element has no child any more, then it is added to the leaf table $L^{(p)}$ accordingly. Repeat the coarsening procedure until no element can be removed from the table $L^{(p)}$. By removing only the leaf element at each time, we avoid generating ``holes" in the hash table. The output of this coarsening procedure are the updated hash table and leaf table, denoted by $H$ and $L$ respectively, and the compressed numerical solution $u_h^{n+1} \in \bV_{N,H}^k$. In practice, $\eta$ is chosen to be smaller than $\varepsilon$ for safety. In the simulations presented in this paper, we use $\eta = \varepsilon/10$.

In summary, the following algorithm advances the numerical solution for one time step.

\medskip	

\noindent\rule{16.5cm}{1pt}\\
\noindent{\bf Algorithm 2: Adaptive evolution from $t^n$ to $t^{n+1}$}\\
\noindent\rule{16.5cm}{0.4pt}\\

 {\bf Input:} Hash table H and leaf table L at $t^n$, numerical solution $u_h^{n} \in \bV_{N,H}^k.$
 
   {\bf Parameters:} Maximum level $N,$ polynomial degree $k,$  error constants   $\varepsilon, \eta,$ CFL constant.

  {\bf Output:} Hash table H and leaf table L at $t^{n+1}$, numerical solution $u_h^{n+1} \in \bV_{N,H}^k.$
 
 \medskip
\begin{enumerate}
	\item {\bf Prediction.} Given a hash table $H$ that stores the numerical solution $u_h$ at time step $t^n$, calculate $\Delta t$. Predict the solution by the DG scheme using space $\bV_{N,H}^k$ and the forward Euler time stepping method. Generate the predicted solution  $u_h^{(p)}$.
	\item {\bf Refinement.} Based on the predicted solution $u_h^{(p)}$, screen all elements in the hash table $H$. If for element $V_\bl^\bj$, the refining criteria \eqref{eq:l1_c}, \eqref{eq:l2_c}, or \eqref{eq:l8_c}  hold,
		%
		%$$\sum_{\mathbf{1}\leq\bi\leq\bk+\mathbf{1}}|u^\bj_{\bi,\bl}|^2<\epsilon$$
		%
	then add its children elements
	to $H$ and $L$  provided they are not added yet, and set the associated detail coefficients to zero. We also need to make sure that all the parent elements of the newly added element are in $H$ (i.e., no ``hole" is allowed in the hash table) and increase the number of children for all its parent elements by one. This step  generates the updated hash table $H^{(p)}$ and leaf table $L^{(p)}$.
	\item{\bf Evolution.} Given the predicted table $H^{(p)}$ and   the leaf table $L^{(p)}$, we evolve the solution from  $t^n$ to $t^{n+1}$ by the DG scheme using space $\bV_{N,H^{(p)}}^k$ and  the third order Runge-Kutta time stepping method \eqref{eq:tvd}. This step  generates the pre-coarsened numerical solution $\tilde{u}_h^{n+1}.$
	\item{\bf Coarsening.} For each element in the leaf table, if the coarsening criteria \eqref{eq:l1_c}, \eqref{eq:l2_c} or \eqref{eq:l8_c}
	 hold, then remove the element from table $H^{(p)}$ and $L^{(p)}$. For each of its parent elements in $H^{(p)}$, we decrease the number of children by one. If the number becomes zero, i.e, the element has no child, then it will be added to leaf table $L^{(p)}$. Repeat the coarsening procedure until no element can  be removed from the leaf list. Denote the resulting hash table and leaf table   by $H$ and $L$ respectively, and the compressed numerical solution $u_h^{n+1} \in \bV_{N,H}^k$.
\end{enumerate}
\noindent\rule{16.5cm}{1pt}\\

\bigskip

\begin{rem}
	The optimal choice of the maximum mesh level $N$ and error parameters $\varepsilon$ and $\eta$ is problem dependent and  the performance of the adaptive scheme is closely related to their choice. For example, an excessively small $\varepsilon$ may result in unnecessary refinement and hence larger computational cost, but little gain in accuracy.
	 On the other hand, if an excessively large $N$ is chosen, then we may need a very small time step for the stability consideration, which degrades the efficiency of the proposed scheme.  
\end{rem}

%It is very important to keep tracking the number of children for each element in $H$, since in what follows we will introduce an  adaptive evolution algorithm, for which we need to remove unnecessary elements from the hash table to reduce the degrees of freedom. We prescribe that only leaf elements are allowed to be removed for accuracy consideration. In other words, 'holes' are not permitted in the hash table. We denote the approximation space that includes all elements in the hash table by $\bV^k_{N,H}$.

\section{Numerical results}
\label{sec:numerical}

In this section, we present benchmark numerical results to demonstrate the performance of the proposed  scheme for solving linear transport equations. For all test examples, we consider both smooth and non smooth initial profiles.

%\begin{exa}[Linear advection with constant coefficient]
%	\label{ex:linear}
%	We consider  
%	\begin{equation}
%	\label{eq:linear_adv}
%	\left\{\begin{array}{l} \displaystyle u_t + \sum_{m=1}^d u_{x_m} = 0,\quad \bx\in[0,1]^d,\\[2mm]
%	\displaystyle u(0,\bx) = \prod_{m=1}^{d} \sin^4\left(\pi x_m\right),
%	\end{array}\right.
%	\end{equation}
%	with periodic boundary conditions.
%\end{exa}
%

\begin{exa}[Linear advection with constant coefficient]
	\label{ex:linear}
	We consider  
	\begin{equation}
	\label{eq:linear_adv}
	 u_t + \sum_{m=1}^d u_{x_m} = 0,\quad \bx\in[0,1]^d
	\displaystyle 
	\end{equation}
	with periodic boundary conditions.
\end{exa}

We first consider a smooth initial condition
\begin{equation}
\label{eq:linear_init_smooth}
u(0,\bx) = \prod_{m=1}^{d} \sin^4\left(\pi x_m\right),
\end{equation}
with $d=2, 3, 4$
and investigate the   accuracy of the scheme using $L^2$ norm based refinement and coarsening criteria  \eqref{eq:l2} and \eqref{eq:l2_c}. We run the simulations with a fixed maximum mesh level  $N=7,$ different  $\varepsilon$ values, and report the $L^2$ errors and the number of active degrees of freedom at final time $T=1$ in Table \ref{table:linear}.   The following rates of convergence are calculated,  
\begin{align*}
	\mbox{convergence rate with respect to the error threshold}  \quad &R_{\varepsilon_l}=\frac{\log(e_{l-1}/e_l)}{\log(\varepsilon_{l-1}/{\varepsilon_l})}\\
	\mbox{convergence rate with respect to DOF } \quad & R_{\text{DOF}_l}=\frac{\log(e_{l-1}/e_l)}{\log(\text{DOF}_l/\text{DOF}_{l-1})}, 
\end{align*} 
where $e_l$ is the standard $L^2$ error with refinement parameter $\varepsilon_l$, and $\text{DOF}_l$ is the associated number of active degrees of freedom at final time. 
For comparison purpose, recall the standard DG schemes with the tensor product grid yields $R_{\epsilon}\approx 1$ and $R_{\text{DOF}}\approx \frac{k+1}{d}$. From Table \ref{table:linear}, we observe that for the proposed scheme,   $R_{\epsilon}$ is slightly smaller than 1, and  $R_{\text{DOF}}$ is much larger than $\frac{k+1}{d}$ but still smaller than $k+1$. This demonstrates the effectiveness of the adaptive algorithm, as well as the computational saving of the multiresolution scheme in this case. We also experiment on varying both $N$ and $\varepsilon$ values at the same time. To save space, the results are not reported but we remark that  
if a excessively small $\varepsilon$ is taken with a small mesh level $N$, the performance of the  scheme will be very similar to the tensor product DG method and the efficiency of the scheme will be adversely affected. We also test the code with $L^1$ and $L^\infty$ based criteria  \eqref{eq:l1}, \eqref{eq:l1_c} and  \eqref{eq:l8}, \eqref{eq:l8_c}, little difference is observed in the convergence order. To save space, the results are omitted in the paper. For the rest of the paper, unless otherwise noted, the refinement and coarsening criteria based on $L^2$ norms \eqref{eq:l2}, \eqref{eq:l2_c} will be used.

%To achieve a smaller desired accuracy, besides an appropriately chosen $\varepsilon$, to be consistent, it may be better to take a larger maximum mesh level $N$.  

\begin{table}[htp]
		\caption{ Example \ref{ex:linear} with initial condition \eqref{eq:linear_init_smooth}. Numerical error and convergence rate. $N=7$. $T=1$.
		}
		%\vspace{2 mm}
		\centering
		\begin{tabular}{|c| c c c c| c c c c| c c c c|}
			\hline
			$\varepsilon$ &  DOF& $L^2$ error &  $R_{\text{DOF}}$& $R_\varepsilon$ &  DOF& $L^2$ error & $R_{\text{DOF}}$& $R_\varepsilon$ &    DOF& $L^2$ error &  $R_{\text{DOF}}$& $R_\varepsilon$ \\
			\hline
			
			&\multicolumn{4}{|c|}{$ k=1$, $ d=2$}& \multicolumn{4}{|c|}{$ k=1$, $ d=3$}&\multicolumn{4}{|c|}{$ k=1$, $ d=4$}  \\
			\hline
			1E-03	&	312	&	1.47E-02	&		&		&	1168	&	2.62E-02	&		&		&	2592	&	2.87E-02	&		&		\\
			5E-04	&	404	&	8.90E-03	&	1.93	&	0.72	&	1840	&	1.87E-02	&	0.75	&	0.49	&	4512	&	2.32E-02	&	0.39	&	0.31	\\
			1E-04	&	1148	&	1.70E-03	&	1.59	&	1.03	&	3920	&	7.26E-03	&	1.25	&	0.59	&	14976	&	9.49E-03	&	0.75	&	0.56	\\
			5E-05	&	1688	&	1.04E-03	&	1.28	&	0.71	&	6440	&	4.16E-03	&	1.12	&	0.80	&	23776	&	6.60E-03	&	0.79	&	0.53	\\
			1E-05	&	3588	&	2.42E-04	&	1.93	&	0.90	&	18624	&	8.83E-04	&	1.46	&	0.96	&	62368	&	2.13E-03	&	1.17	&	0.70	\\
			5E-06	&	4636	&	1.37E-04	&	2.23	&	0.82	&	25496	&	5.10E-04	&	1.75	&	0.79	&	111424	&	1.18E-03	&	1.02	&	0.86	\\
			%1E-06	&	8132	&	6.06E-05	&	1.45	&	0.51	&		&		&		&		&		&		&		&		\\
			
			\hline
			&\multicolumn{4}{|c|}{$ k=2$, $ d=2$}& \multicolumn{4}{|c|}{$ k=2$, $ d=3$}&\multicolumn{4}{|c|}{$k=2,$   $d=4$}  \\
			\hline
%			1E-03	&	252	&	7.15E-03	&		&		&	1269	&	1.03E-02	&		&		&	3807	&	1.16E-02	&		&		\\
%			5E-04	&	288	&	2.16E-03	&	8.95	&	1.72	&	1593	&	9.40E-03	&	0.40	&	1.32	&	7533	&	9.09E-03	&	0.36	&	0.36	\\
%			1E-04	&	621	&	7.45E-04	&	1.39	&	0.66	&	4023	&	2.36E-03	&	1.49	&	0.86	&	15876	&	2.63E-03	&	1.66	&	0.77	\\
			5E-05	&	774	&	3.61E-04	&		&		&	4428	&	1.30E-03	&		&		&	26244	&	1.48E-03	&		&		\\
			1E-05	&	1584	&	8.78E-05	&	1.97	&	0.88	&	9585	&	2.58E-04	&	2.10	&	1.01	&	51840	&	5.30E-04	&	1.51	&	0.64	\\
			5E-06	&	1998	&	4.58E-05	&	2.80	&	0.94	&	13716	&	1.74E-04	&	1.09	&	0.57	&	69012	&	2.60E-04	&	2.49	&	1.03	\\
			1E-06	&	4023	&	1.43E-05	&	1.67	&	0.73	&	27081	&	4.15E-05	&	2.11	&	0.89 &	168723	&	9.46E-05	&	1.13	&	0.63	\\
			5E-07	&	5157	&	7.20E-06	&	2.76	&	0.99	&	40446	&	2.45E-05	&	1.32	&	0.76	&	226719	&	4.89E-05	&	2.23	&	0.95	\\
			1E-07	&	9072	&	1.80E-06	&	2.46	&	0.86	&	77463	&	7.06E-06	&	1.91	&	0.77	&	531684	&	1.24E-05	&	1.61	&	0.85	\\
			%5E-08	&	11034	&	9.29E-07	&	3.37	&	0.95	&	107757	&	4.21E-06	&	1.56	&	0.74	&		&		&		&		\\
			\hline
			&\multicolumn{4}{|c|}{$ k=3$, $ d=2$}& \multicolumn{4}{|c|}{$ k=3$, $ d=3$}&\multicolumn{4}{|c|}{$ k=3$, $ d=4$}  \\
			\hline
			%1E-03	&		&		&		&		&	1664	&	2.18E-03	&		&		&	12032	&	2.20E-03	&		&		\\
			%5E-04	&		&		&		&		&	1664	&	2.06E-03	&		&	0.08	&	12032	&	2.18E-03	&		&	1.18E-02	\\
			%1E-04	&		&		&		&		&	3584	&	4.42E-04	&	2.01	&	0.96	&	12032	&	2.18E-03	&		&	0.00	\\
			%5E-05	&	512	&	1.59E-04	&		&		&	5504	&	2.41E-04	&	1.41	&	0.87	&	31744	&	4.73E-04	&	1.57	&	2.20E+00	\\
			1E-05	&	1120	&	3.71E-05	&		&		&	10496	&	5.72E-05	&		&		&	58368	&	1.26E-04	&		&	\\
			5E-06	&	1184	&	2.92E-05	&	4.32	&	0.35	&	12032	&	4.91E-05	&	1.12	&	0.22	&	97280	&	7.53E-05	&	1.01	&	0.74	\\
			1E-06	&	2208	&	9.87E-06	&	1.74	&	0.67	&	18688	&	1.31E-05	&	3.00	&	0.82	&	129024	&	3.73E-05	&	2.49	&	0.44	\\
			5E-07	&	2864	&	4.85E-06	&	2.73	&	1.03	&	25984	&	1.09E-05	&	0.56	&	0.27	&	204800	&	1.34E-05	&	2.21	&	1.47	\\
			1E-07	&	3968	&	1.31E-06	&	4.02	&	0.82	&	43840	&	2.71E-06	&	2.66	&	0.86	&	409600	&	6.14E-06	&	1.13	&	0.49	\\
			5E-08	&	5760	&	7.88E-07	&	1.36	&	0.73	&	57472	&	1.50E-06	&	2.20	&	0.86	&	521216	&	2.79E-06	&	3.27	&	1.14	\\
%			1E-08	&	8064	&	9.48E-08	&	6.30	&	1.32	&	101632	&	3.47E-07	&	2.57	&	0.91	&		&		&		&		\\		\hline
		\hline	
		\end{tabular}
		\label{table:linear}
	\end{table}

\begin{figure}[htp]
	\begin{center}
		\subfigure[]{\includegraphics[width=.42\textwidth]{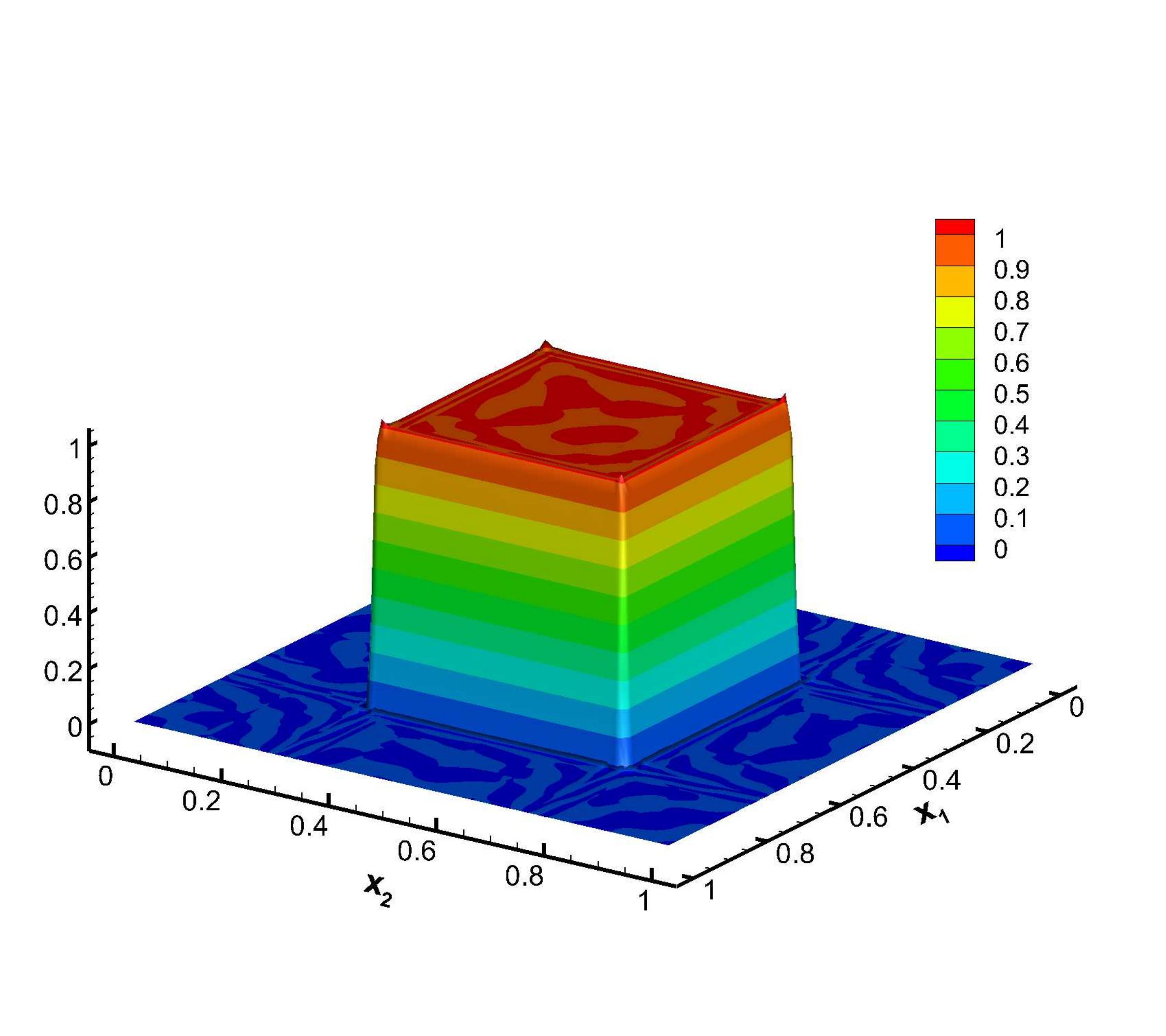}}
		\subfigure[]{\includegraphics[width=.42\textwidth]{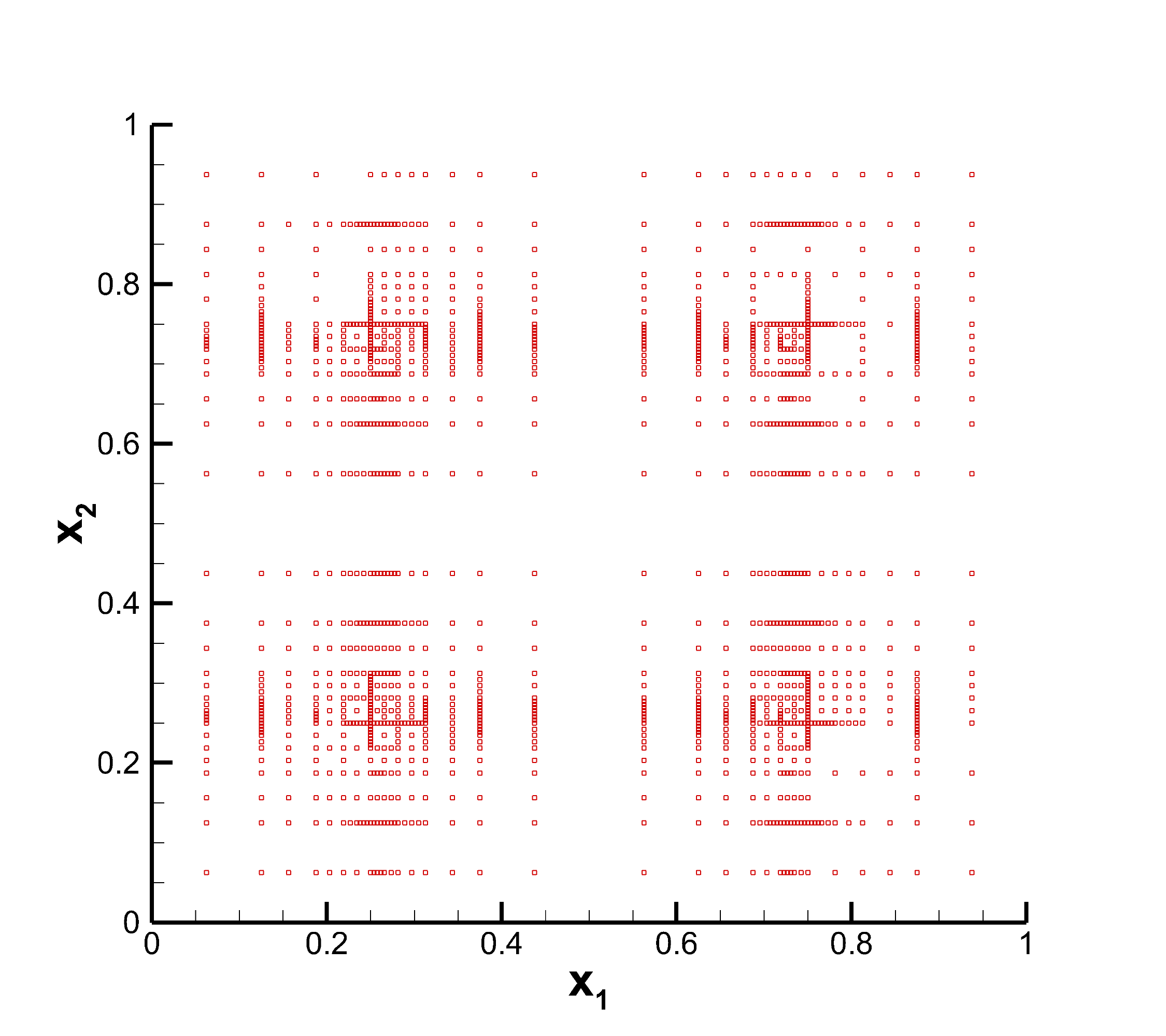}}\\
		\subfigure[]{\includegraphics[width=.42\textwidth]{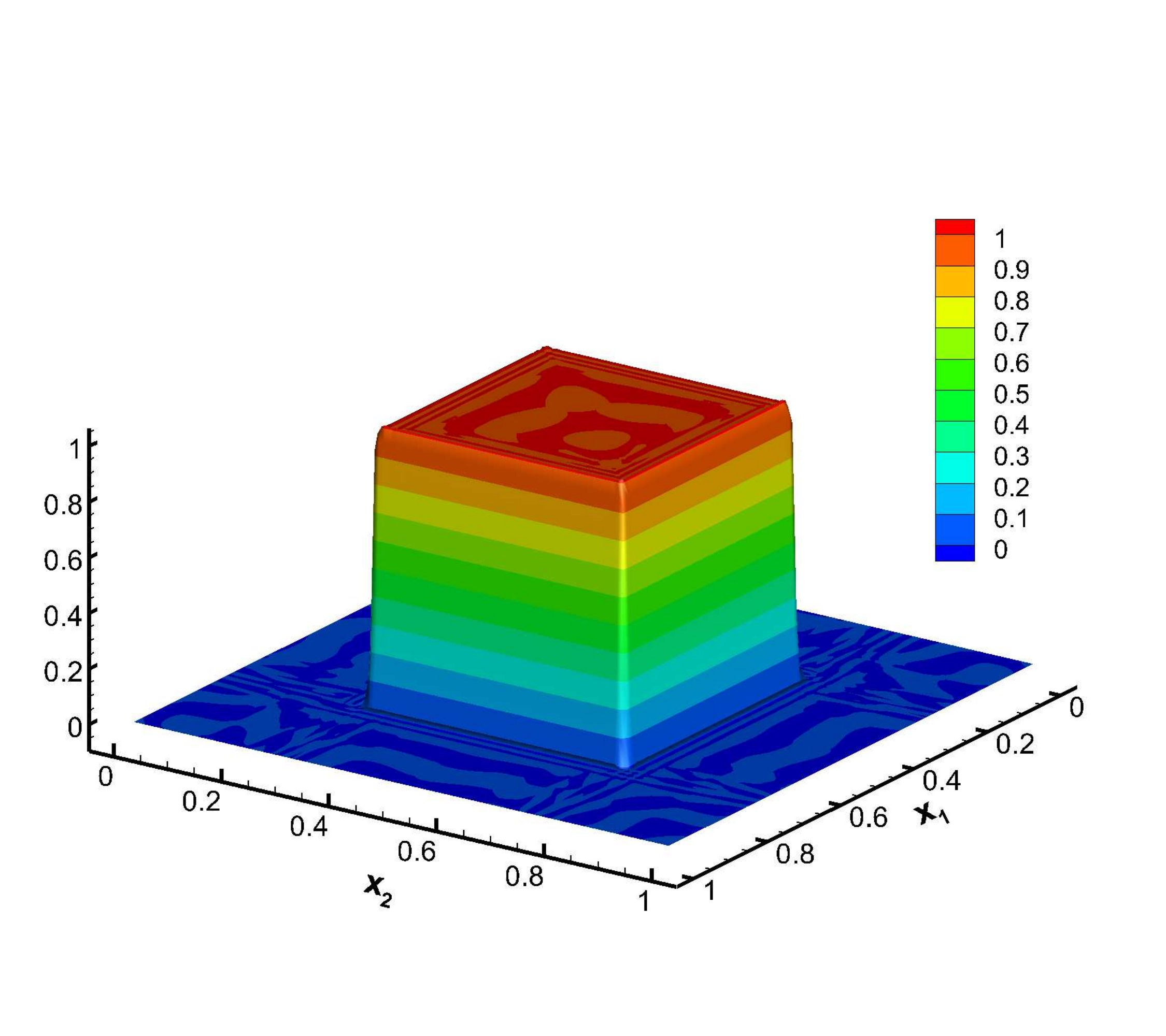}}
		\subfigure[]{\includegraphics[width=.42\textwidth]{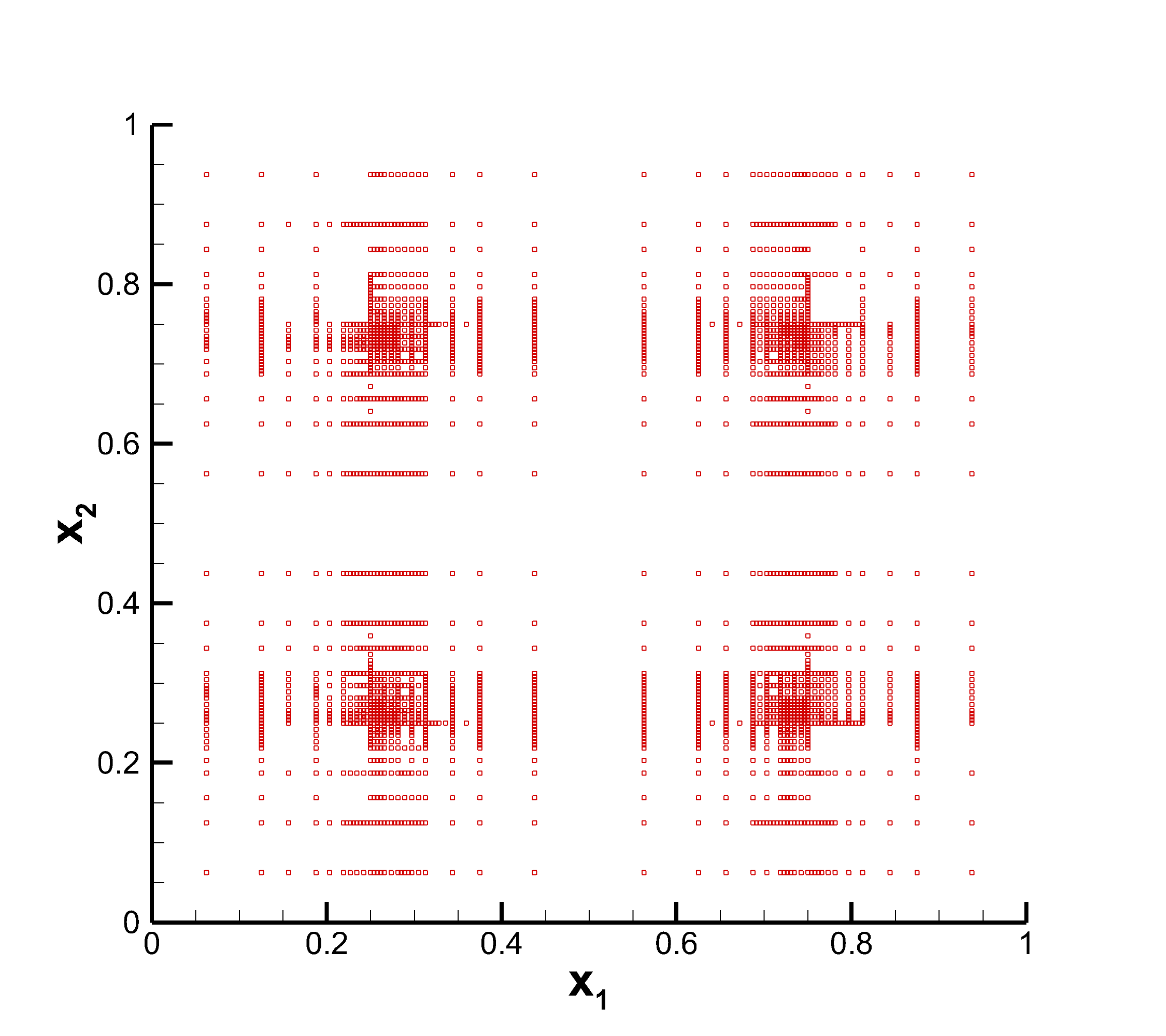}}\\
		\subfigure[]{\includegraphics[width=.42\textwidth]{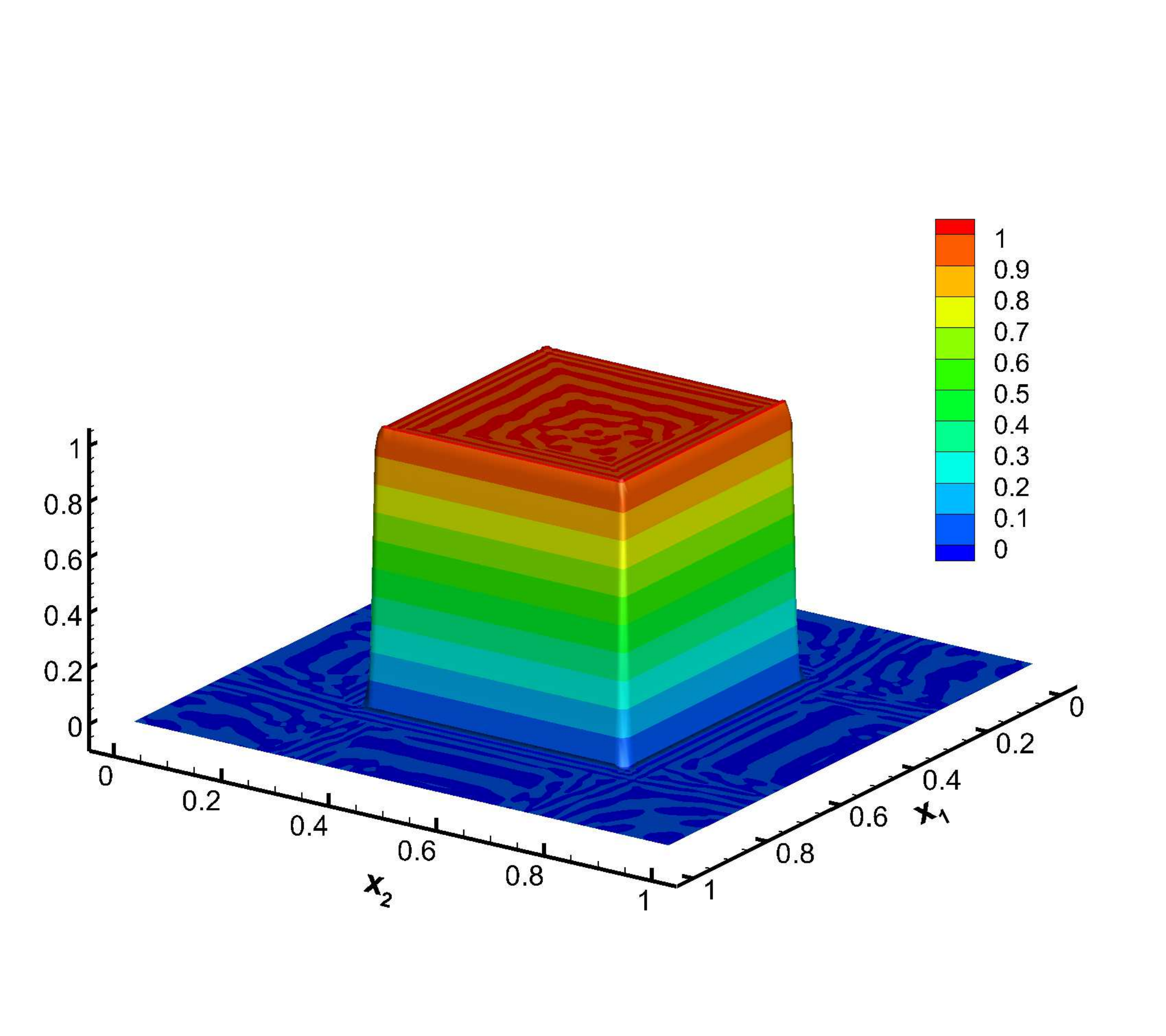}}
		\subfigure[]{\includegraphics[width=.42\textwidth]{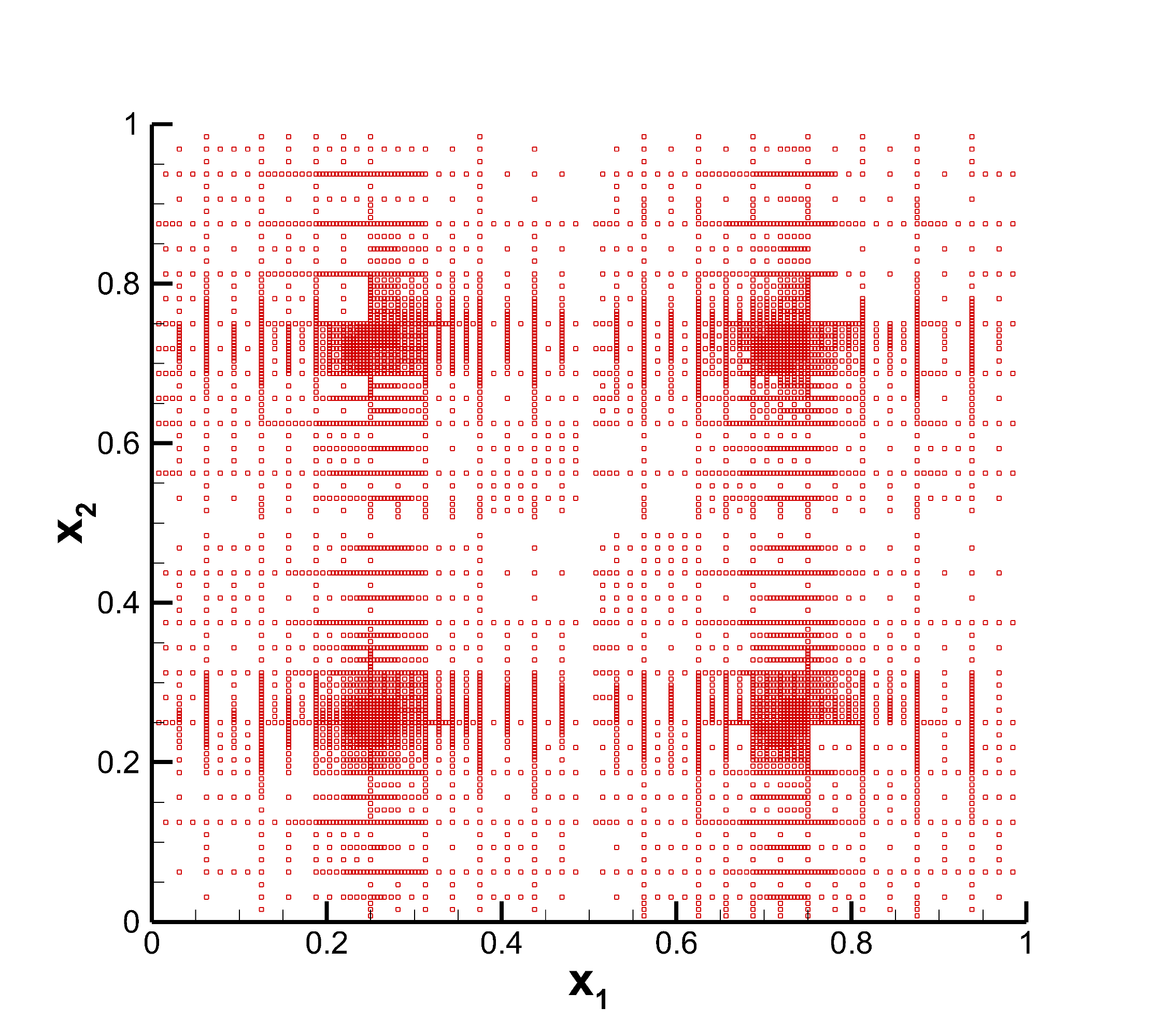}}
	\end{center}
	\caption{Example \ref{ex:linear} with initial condition \eqref{eq:discontinuous}.  (a) Numerical solution with $L^1$ norm based criteria. (b) Active elements with $L^1$ norm based criteria. (c) Numerical solution with $L^2$ norm based criteria. (d) Active elements with $L^2$ norm based criteria. (e) Numerical solution with $L^\infty$ norm based criteria. (f) Active elements with $L^\infty$ norm based criteria. $T=1$. $d=2$. $N=7$. $k=3$. $\varepsilon=10^{-5}$.}
	\label{fig:linear_dis}
\end{figure}

Next, we consider a discontinuous initial condition 
\begin{equation}
\label{eq:discontinuous}
u(0,\bx)=\left\{\begin{array}{ll}1& (x_1,x_2)\in[\frac12-\frac{\sqrt{6}}{2},\frac12+\frac{\sqrt{6}}{2}]^2.\\[2mm]
0& \text{otherwise},\end{array}\right.
\end{equation}
when $d=2.$
It is well known that the standard sparse grid method without adaptivity cannot resolve  such discontinuous solution profiles. 
In our simulations, we fix $N=7, \varepsilon=10^{-5}$ and compare the performance of the scheme with   $L^1$, $L^2$ and $L^\infty$ based refinement/coarsening criteria up to final time $T=1$. The  numerical solutions and the associated active elements are reported in Figure \ref{fig:linear_dis}. We only plot the center of support for each active basis, while noting that the basis functions contains different size of support in the scheme. The method with all three types of refinement/coarsening criteria provides well resolved solution profiles. Active elements all cluster towards the discontinuities.  However, the $L^\infty$ norm based criteria has the most degrees of freedom, increasing computational cost while not improving   numerical performance. Similar comments are also made in \cite{griebel1998adaptive}. The $L^1$ norm based criteria is the most sparse, but the solution is slightly  more oscillatory. This is natural since no limiting procedure has been employed in this paper.   %Mild oscillations are observed in the vicinity of discontinuities since no limiting procedure is used.
 %generates qualitatively similar numerical results. It is observed that the active elements cluster towards the discontinuities and the associated $L^1$ errors are comparable. However, the degrees of freedom with the $L^1$ norm as the error indicator is less than those with the other two norms indicating that there may be unnecessary refinement that cause more computational cost while not improving the numerical performance. Hence, we claim that if the solution contains discontinuities, the $L^1$ norm is preferred as a refinement indicator. Below, we only report the results with the $L^2$ norm as the refinement indicator for brevity. Also note that, some mild oscillations are observed in the vicinity of discontinuities since no limiting procedure is added to remove such undesired oscillations. In our future work, we will design a effective limiter that is well-suited to our proposed adaptive MRA DG scheme.

\begin{figure}[htp]
	\begin{center}
		\subfigure[]{\includegraphics[width=.42\textwidth]{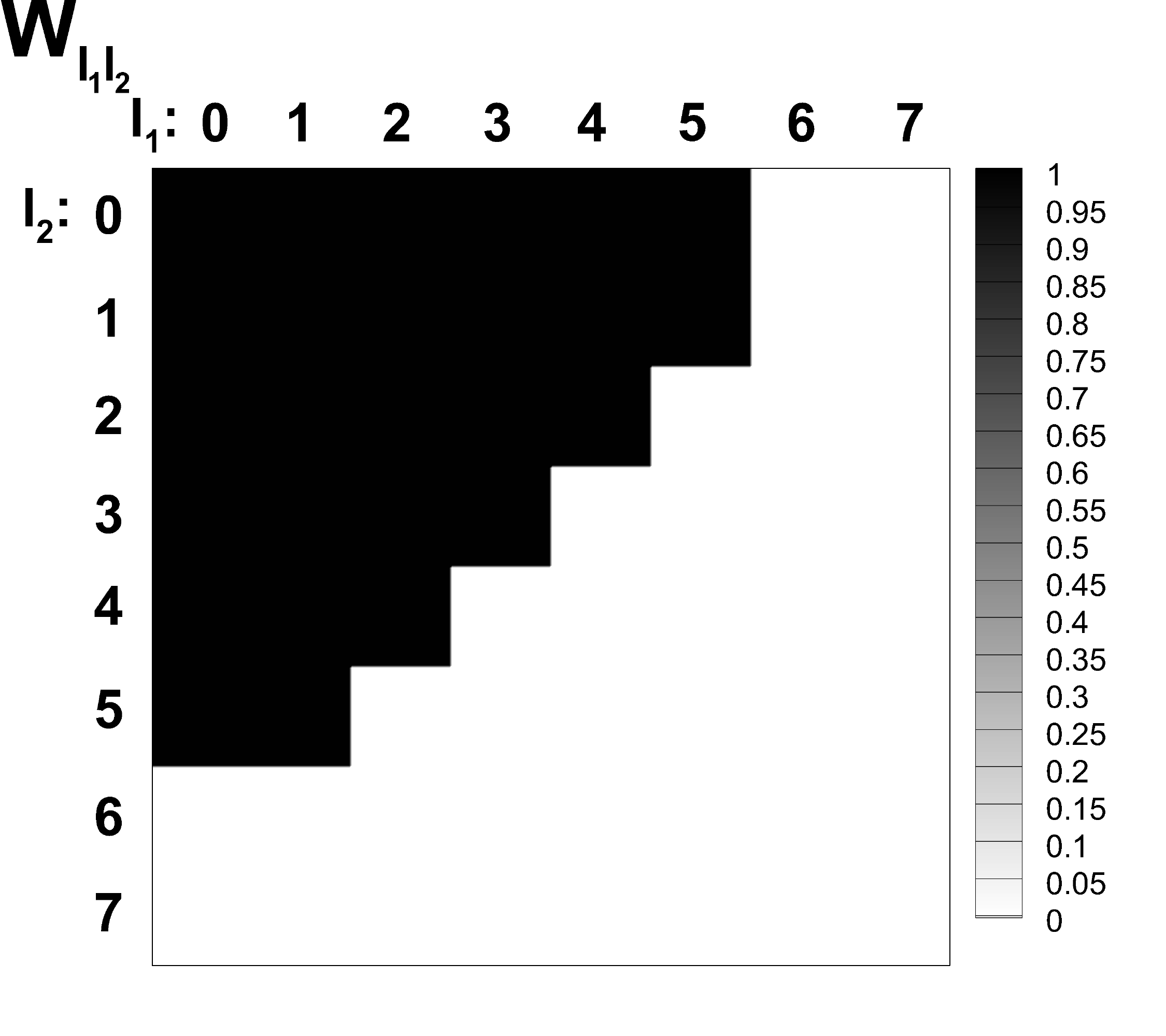}}
		\subfigure[]{\includegraphics[width=.42\textwidth]{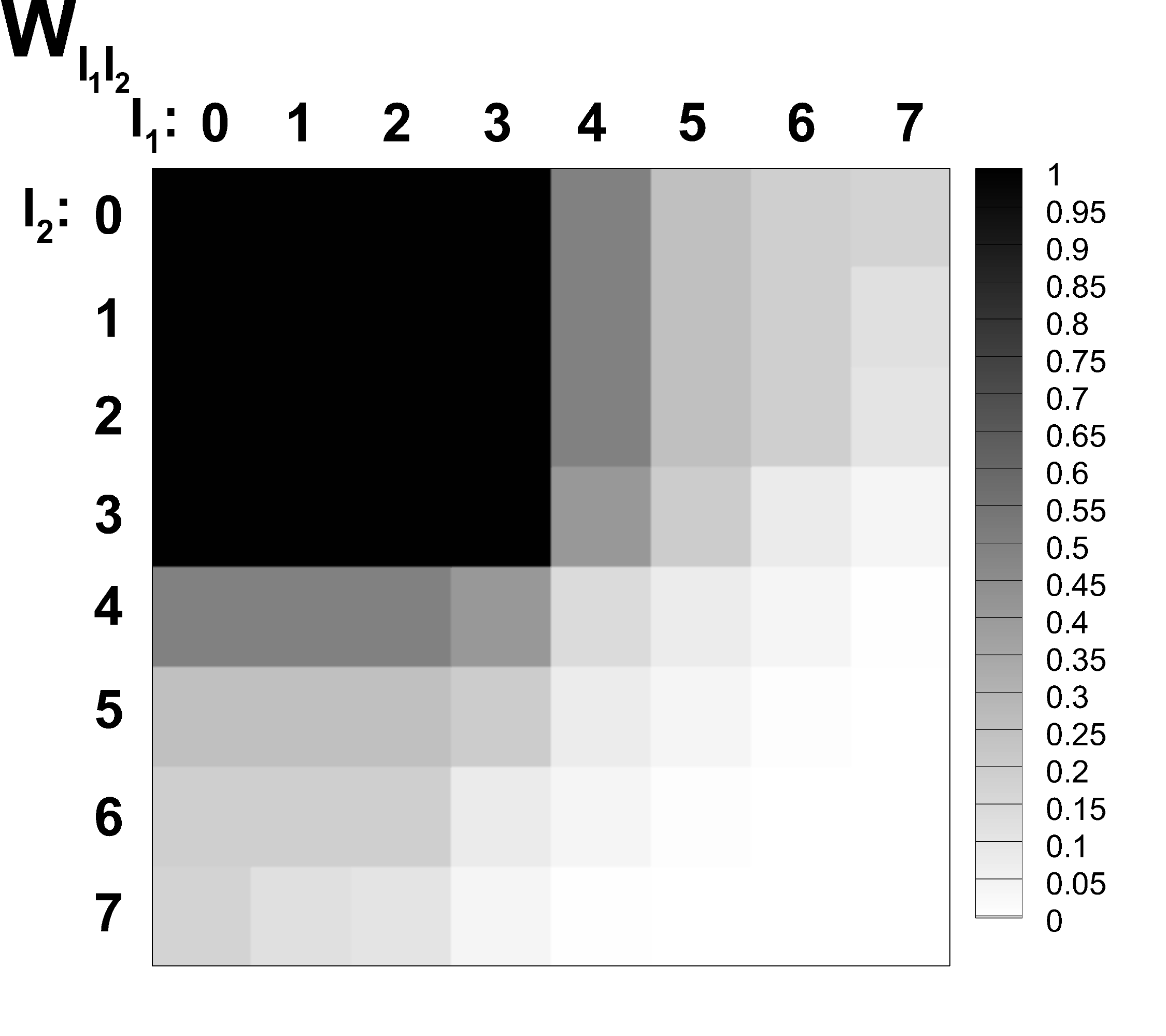}}\\
		\subfigure[]{\includegraphics[width=.42\textwidth]{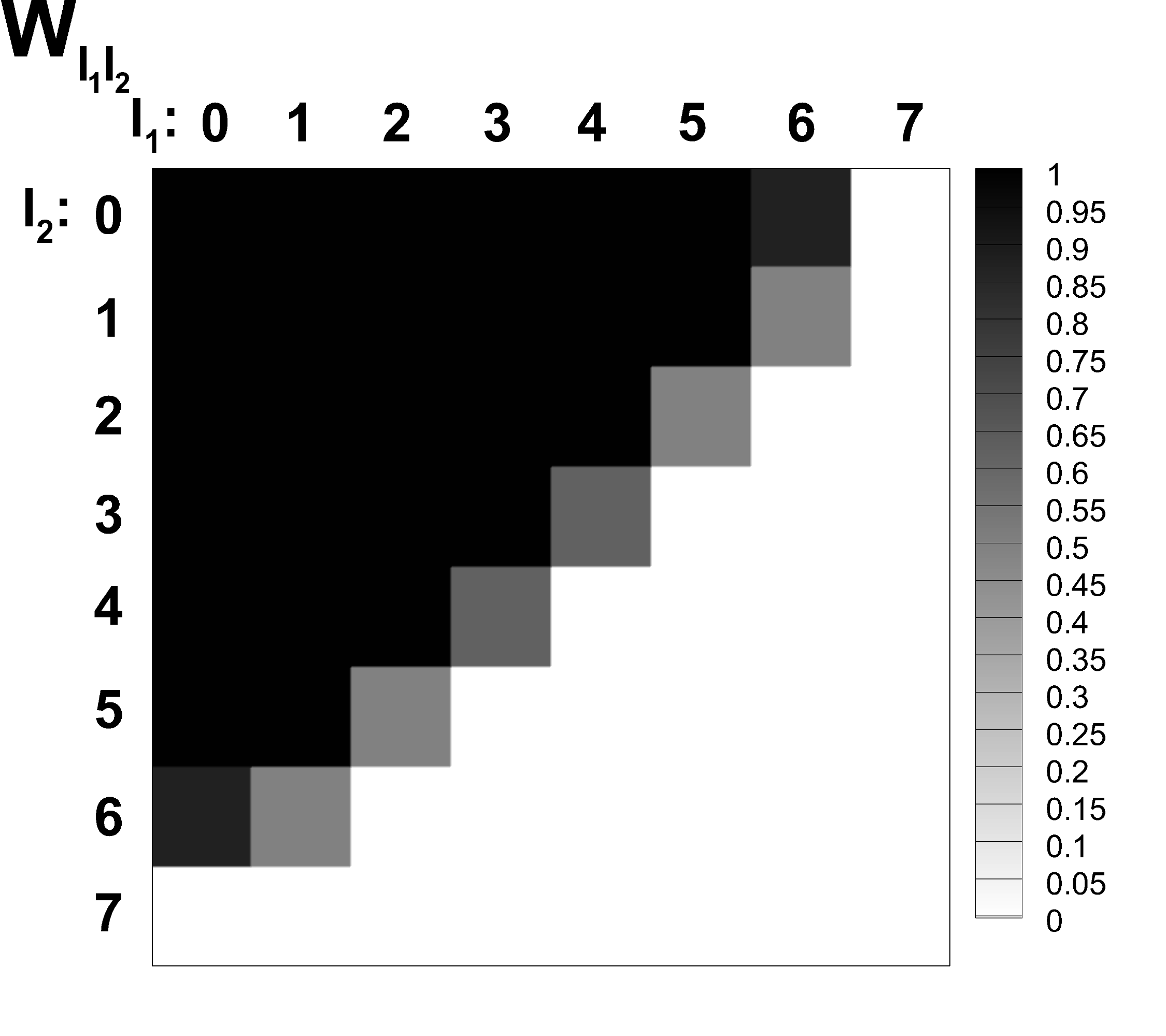}}
		\subfigure[]{\includegraphics[width=.42\textwidth]{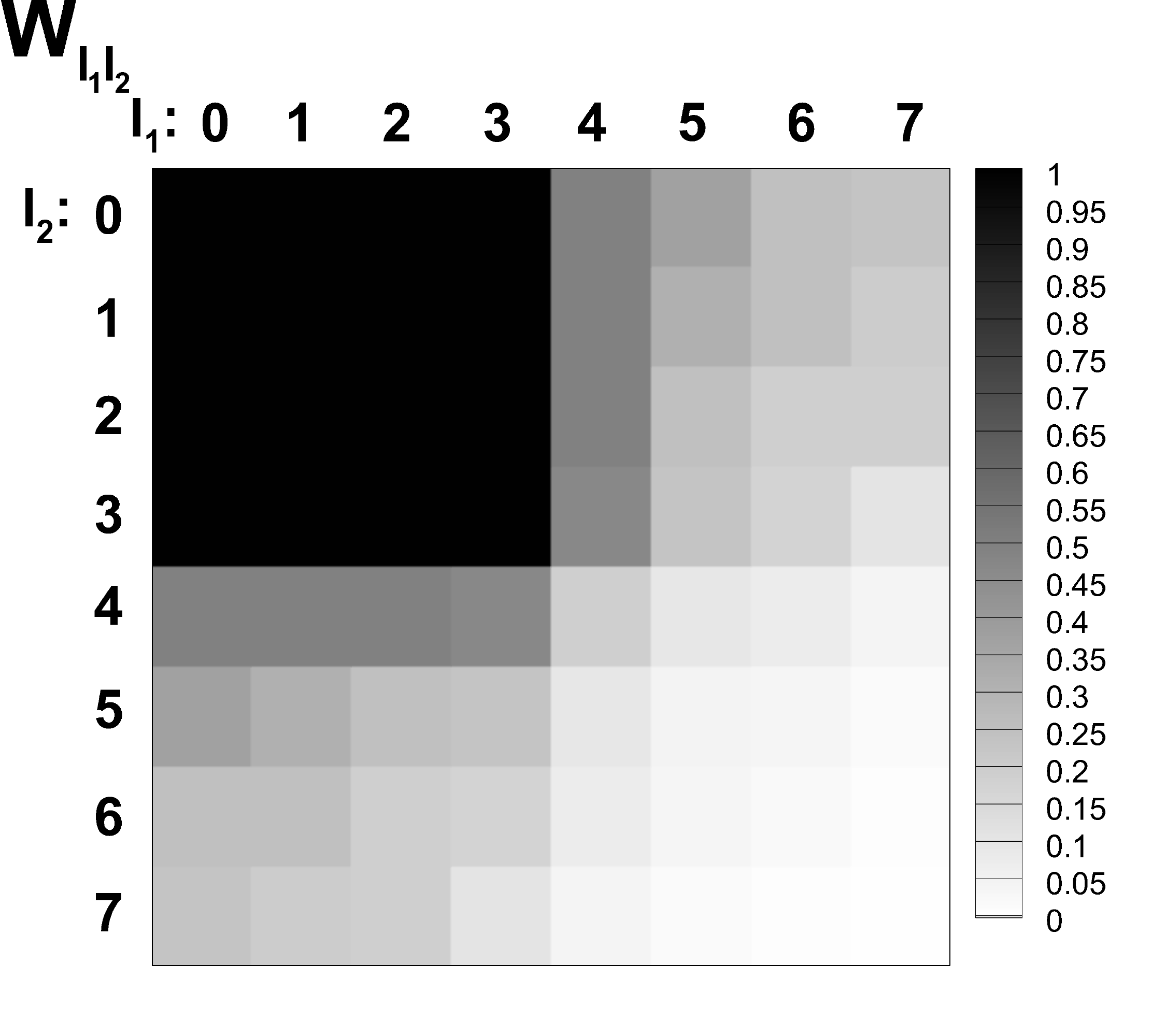}}\\
		\subfigure[]{\includegraphics[width=.42\textwidth]{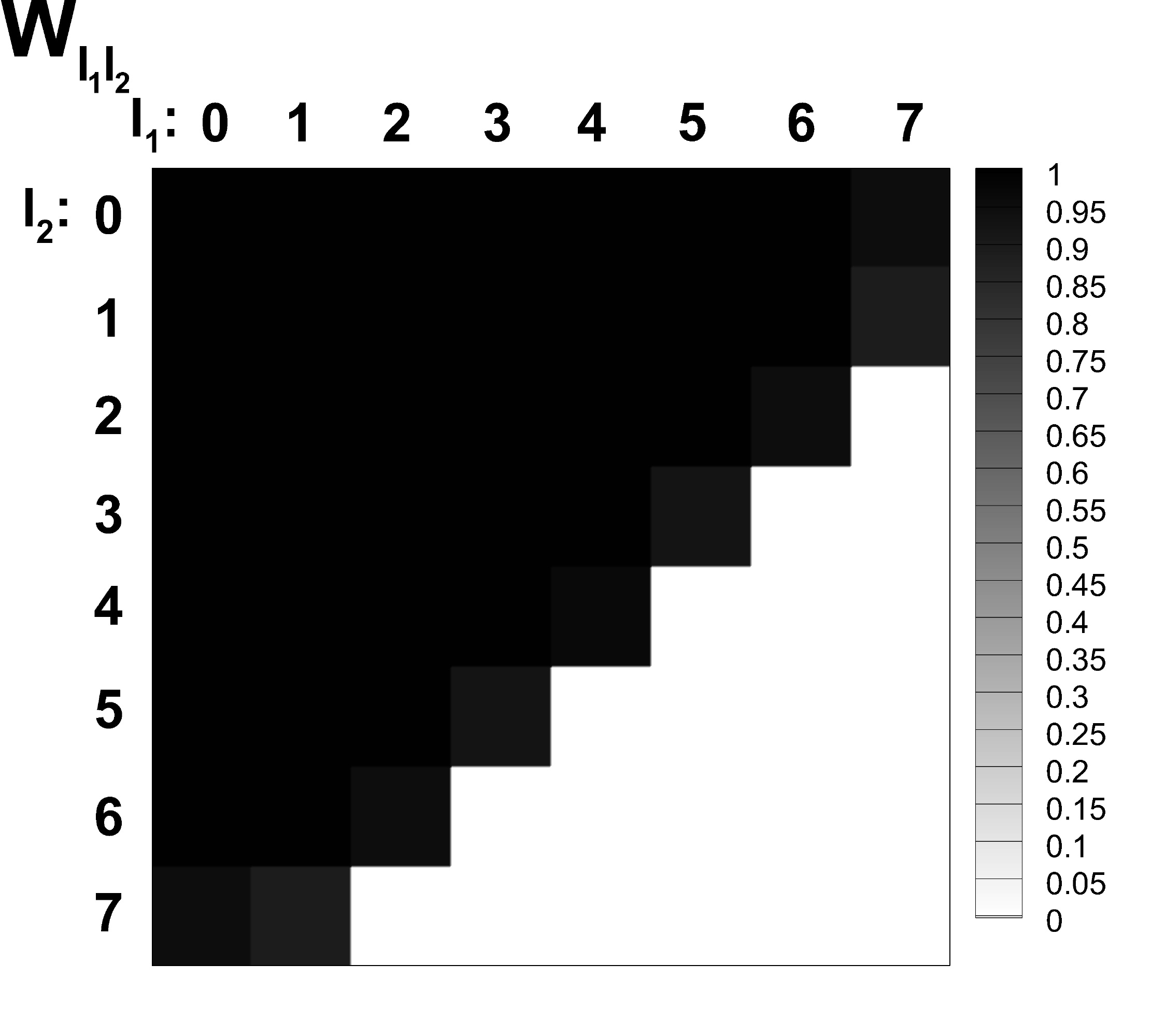}}
		\subfigure[]{\includegraphics[width=.42\textwidth]{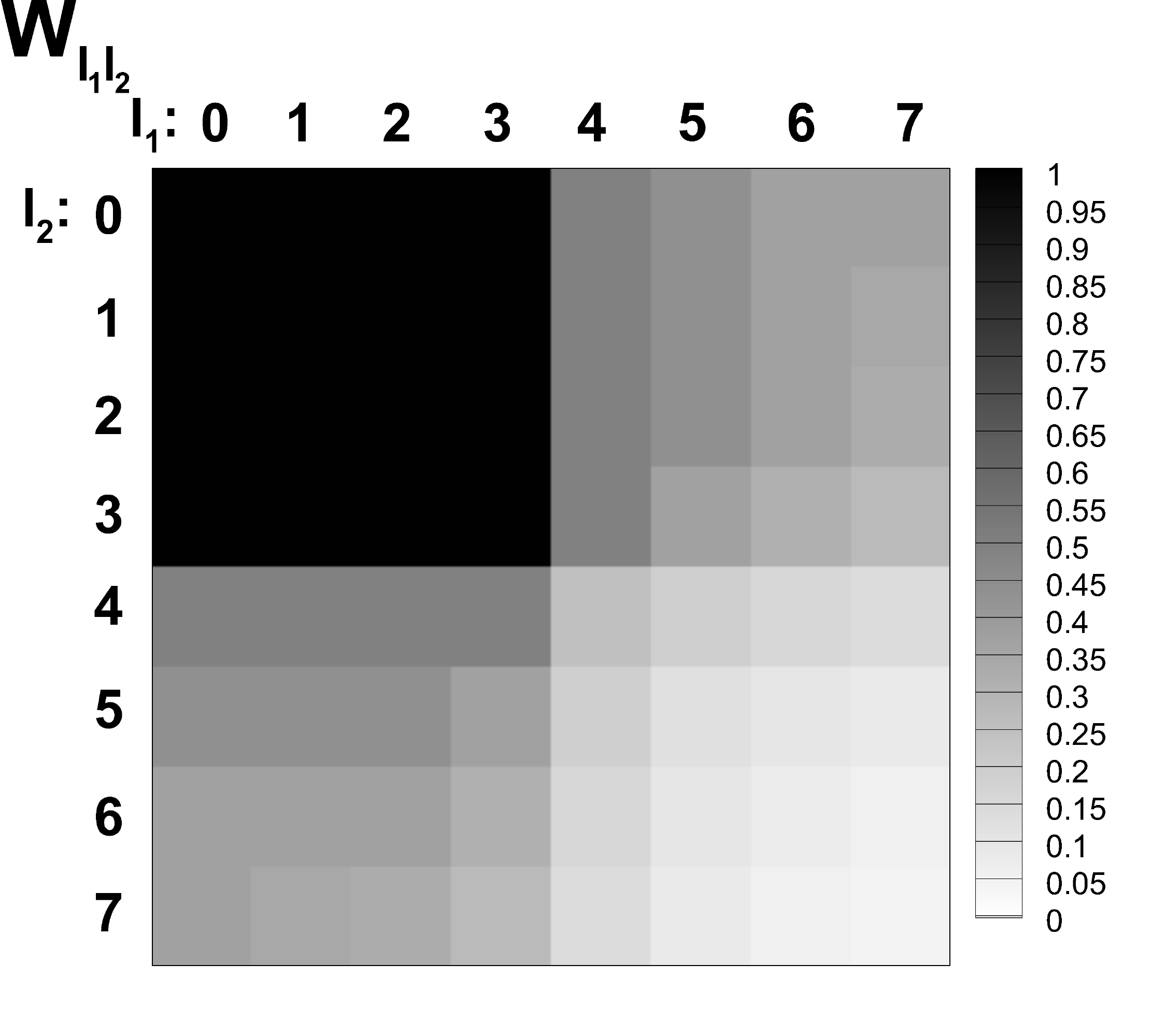}}
	\end{center}
	\caption{Example \ref{ex:linear}. The percentage of active  elements for each incremental space $\bW_\bl$, $\bl=\{l_1,l_2\}$ and $|\bl|_\infty\leq N$ at final time. $T=1$.  $N=7$. $k=3$. $\varepsilon=10^{-7}$ and $\varepsilon=10^{-5}$ for the smooth  (a, c, e) and discontinuous (b, d, f) initial conditions, respectively. We use $L^1$ norm based criteria (a, b), $L^2$ norm  based criteria (c, d), and $L^\infty$ norm  based criteria (e, f).}
	\label{fig:linear_percentage}
\end{figure}

We then perform a detailed comparison for smooth and non smooth solution
 to  demonstrate an important property of the proposed scheme. We fix $d=2$,  $N=7$, $k=3$ and consider initial conditions \eqref{eq:linear_init_smooth} and  \eqref{eq:discontinuous}. We take $\varepsilon=10^{-7}$ and $\varepsilon=10^{-5}$ for the smooth and discontinuous problems, respectively.
In Figure \ref{fig:linear_percentage}, we   plot the percentage of active elements for each incremental space $\bW_\bl$, $\bl={l_1,l_2}$ at final time $T=1$ with all three norms as adaptive indicators. If the percentage is $1,$ it means all the elements on that level is enacted. If the percentage is $0,$ it means no element on that level is enacted. A full grid approximation corresponds to percentage being $1$ on all levels, while a sparse grid approximation \cite{guo_sparsedg} corresponds to percentage being $1$ when $|\bl|_1 \le N,$ and $0$ otherwise. For the adaptive scheme, there is no longer a clean cutoff and we visualize the variation of percentages among all levels when $L^1$, $L^2$ and $L^\infty$ norm based criteria are used.
 We observe from Figure \ref{fig:linear_percentage}  that  only the upper left corner of incremental spaces are active, similar to the sparse grid DG method with space approximation $\hat{\bV}_N^k$ when the solution is smooth. This is true for all refinement/coarsening criteria.
 If the solution is discontinuous,   more elements are incorporated to fully resolve the discontinuities.  The $L^1$ norm based criteria is the most sparse among the three as expected. From this plot, we can conclude that if the solution is globally smooth, then the scheme will go back to a sparse grid DG method proposed in \cite{guo_sparsedg}, leading to great savings in computational cost; otherwise, the adaptive algorithm will automatically use more elements in the refined levels to capture local fine structures.

\begin{figure}[htp]
	\begin{center}
		\subfigure[]{\includegraphics[width=.42\textwidth]{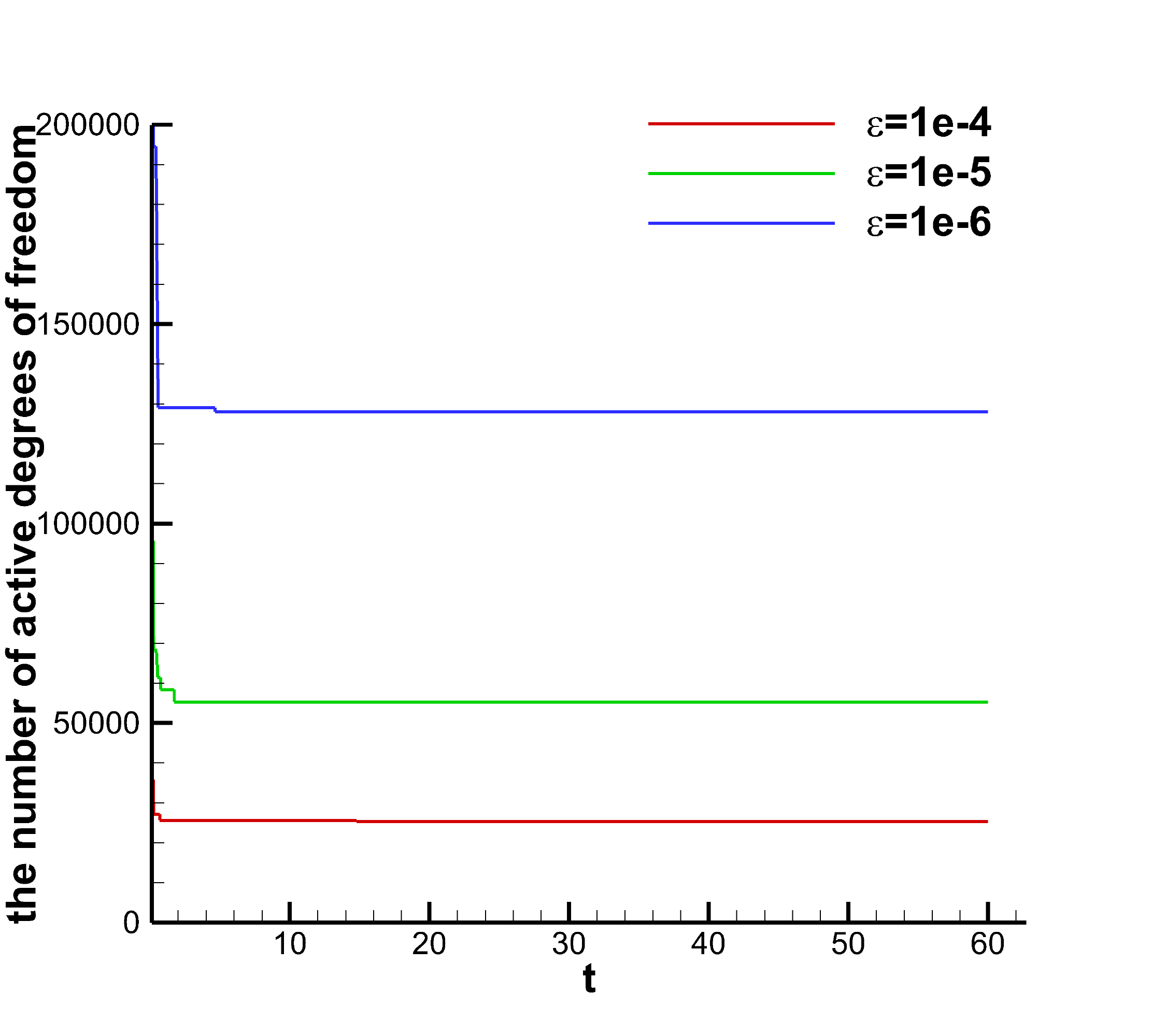}}
		\subfigure[]{\includegraphics[width=.42\textwidth]{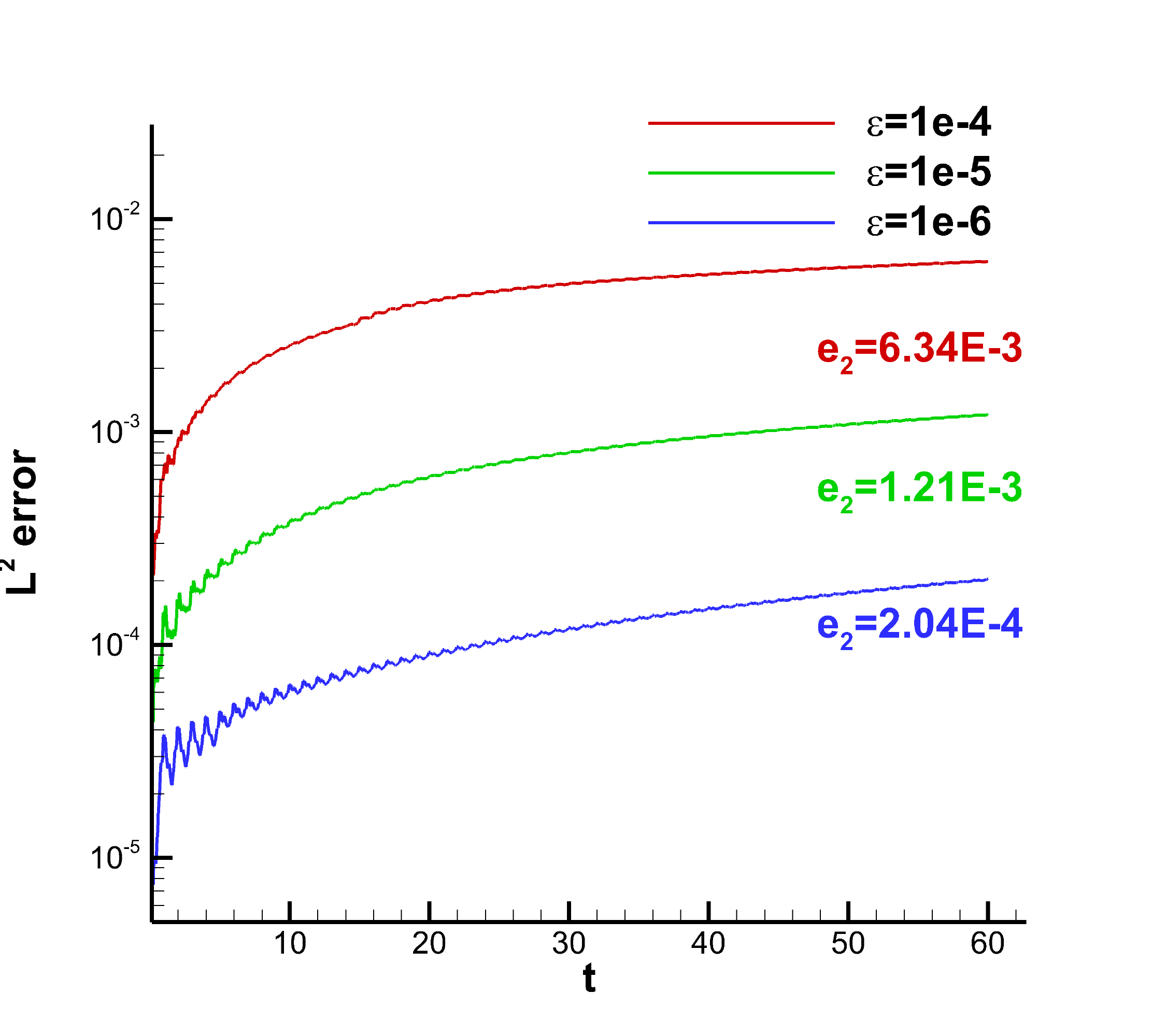}}
	\end{center}
	\caption{Example \ref{ex:linear} with initial condition \eqref{eq:linear_init_smooth} and $L^2$ norm based criteria.   (a) Time histories of the number of active degrees of freedom. (b) Time histories of $L^2$ errors. $T=60$. $d=4$. $N=7$. $k=3$.}
	\label{fig:linear_d4_n7}
\end{figure}

\begin{figure}[htp]
	\begin{center}
		\subfigure[]{\includegraphics[width=.42\textwidth]{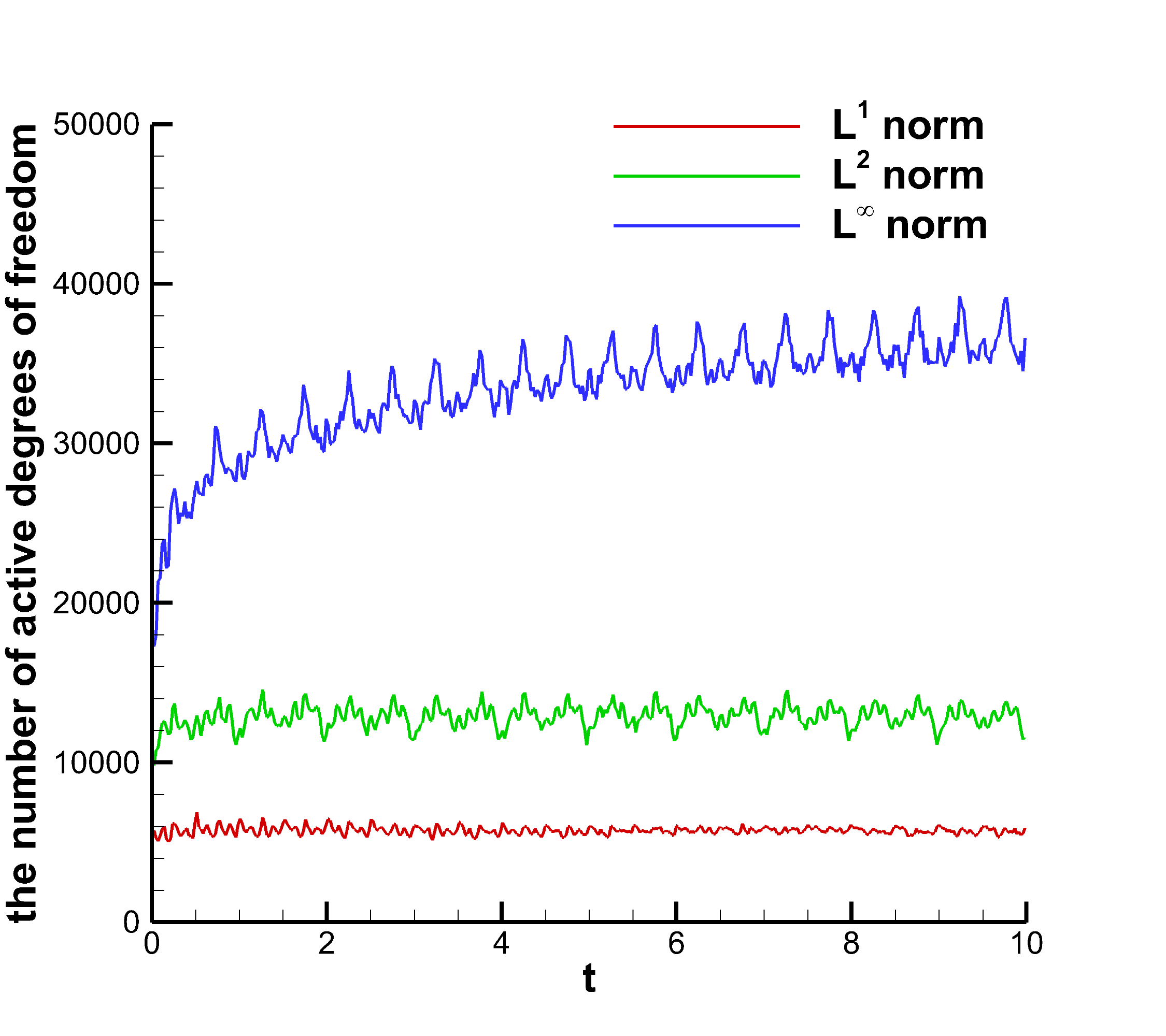}}
		\subfigure[]{\includegraphics[width=.42\textwidth]{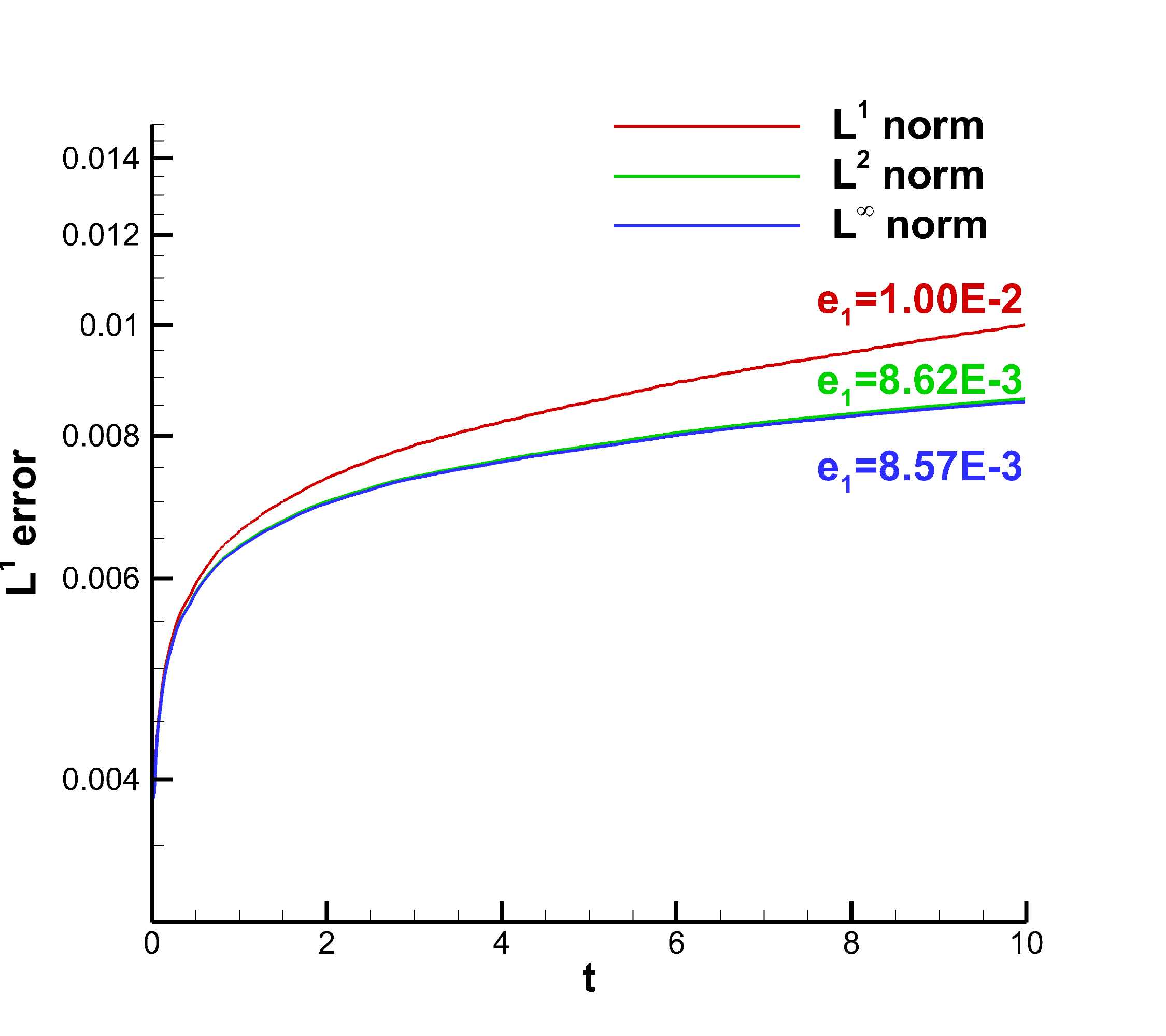}}
	\end{center}
	\caption{Example \ref{ex:linear} with initial condition \eqref{eq:discontinuous}. (a) Time histories of the number of active degrees of freedom. (b) Time histories of $L^1$ errors. $T=10$. $d=2$. $N=7$. $k=3$. $\varepsilon=10^{-5}$.}
	\label{fig:linear_dis_time}
\end{figure}

An additional point we are concerned with is the long time performance of the scheme. For the smooth initial condition \eqref{eq:linear_init_smooth}, 
we set $T=60, d=4, N=7, k=3$    and keep track of the time evolution of  $L^2$ errors and the numbers of active degrees of freedom with $\varepsilon=10^{-4},\,10^{-5},\,10^{-6}$ along time evolution as shown in Figure \ref{fig:linear_d4_n7}. It is observed that, for this linear transport problem, the active degrees of freedom decrease at the very beginning of the simulations, then they nearly remain constant as time evolves for all $\varepsilon$. This is because the profile of solution does not change over time and it is only advected along the characteristic direction. The $L^2$ error demonstrates sub-linear growth beyond the initial stage. The maximum $L^2$ errors over time are reported in the figure. For the discontinuous initial condition \eqref{eq:discontinuous}, we set $T=10, d=2, N=7, k=3, \varepsilon=10^{-5}$ and  report  the time evolution the numbers of active degrees of freedom and $L^1$ errors in Figure \ref{fig:linear_dis_time} with all three norm based criteria. The scheme with the $L^\infty$ norm based criteria yields the smallest error but also involves the largest number of degrees of freedom. The error performance of  the $L^2$ norm based criteria is qualitatively the same as  the $L^\infty$ norm based criteria, while much less degrees of freedom are used. The $L^1$ norm based criteria leads to the largest error while it is uses the least amount of elements among the three. The maximum $L^1$ errors over time are reported in the figure. 

%In Table \ref{table:linear_d2}, we report the degrees of freedom of the associated space, $L^2$ errors and  orders of accuracy for $k=1, 2, 3$ and up to dimension four. The degrees of freedom of the computational method are significantly reduced when compared with the traditional DG space. As for accuracy, we observe half order reduction from the optimal $(k+1)$-th order for high-dimensional computations (d=4). The order is slightly better for lower dimensions. The conclusions from this example agree well with the error estimate in Theorem \ref{thm:conv}.

\begin{exa}[Solid body rotation]
	\label{ex:rotation}
	We consider solid body rotation, which is in the form of \eqref{eq:model} with
	$$\ba = \left(-x_2+\frac12, x_1-\frac12\right),\quad\text{when}\quad d=2,$$
	$$\ba = \left(-\frac{\sqrt{2}}{2}\left(x_2-\frac12\right), \frac{\sqrt{2}}{2}\left(x_1-\frac12\right)+ \frac{\sqrt{2}}{2}\left(x_3-\frac12\right),-\frac{\sqrt{2}}{2}\left(x_2-\frac12\right)\right),\quad\text{when}\quad d=3,$$
	subject to periodic boundary conditions.
\end{exa}

This benchmark test is used to assess the performance of the sparse grid DG transport schemes \cite{guo_sparsedg}. The initial condition is set to be the following smooth cosine bells (with $C^5$ smoothness),
\begin{equation}\label{eq:cosine} u(0,\bx)=\left\{\begin{array}{ll}b^{d-1}\cos^6\left(\frac{\pi r}{2b}\right),& \text{if}\quad r\leq b,\\
0,&\text{otherwise},
\end{array}\right.\end{equation}
where $b=0.23$ when $d=2$ and $b=0.45$ when $d=3$, and $r=|\bx-\bx_c|$ denotes the distance between $\bx$ and the center of the cosine bell with $\bx_c=(0.75,0.5)$ for $d=2$ and  $\bx_c=(0.5,0.55,0.5)$ for $d=3.$
As time evolves, the cosine bell traverses along circular trajectories centered at $(1/2,1/2)$ for $d=2$ and about the axis $\{x_1=x_3\}\cap \{x_2=1/2\}$ for $d=3$ without deformation. We start with the investigation of the convergence rate of the adaptive  scheme. Similar to the previous example, we run the simulation up to $T=1$ with different $\varepsilon$ and summarize the $L^2$ errors, the number of active degrees of freedom and corresponding convergence rates $R_\varepsilon$ and $R_{\text{DOF}}$   in Table \ref{table:solid_d2}. The maximum mesh level is set as $N=7$. For both $d=2,\,3$, it is observed that the rate $R_\varepsilon$ is slightly less than 1 and $R_{\text{DOF}}$ is larger than $\frac{k+1}{d}$ but smaller than $k+1$, which is similar to the previous example.

\begin{table}[htp]
	\caption{Example \ref{ex:rotation} with initial condition \eqref{eq:cosine}. Numerical error and convergence rate. $N=7$. $T=1$.  
	}
	%\vspace{2 mm}
	\centering
	\begin{tabular}{|c| c c c c| c c c c|}
		\hline
$\varepsilon$ &  DOF& $L^2$ error &  $R_{\text{DOF}}$& $R_\varepsilon$ &  DOF& $L^2$ error & $R_{\text{DOF}}$& $R_\varepsilon$ \\
\hline

&\multicolumn{4}{|c|}{$ k=1$, $ d=2$}& \multicolumn{4}{|c|}{$ k=1$, $d=3$}\\
		\hline
5E-04	&	260	&	6.71E-03	&		&		&	928	&	2.17E-03	&		&   	\\
1E-04	&	604	&	1.53E-03	&	1.76	&	0.92	&	3280	&	6.32E-04	&	0.98	&	1.78	\\
5E-05	&	764	&	8.07E-04	&	2.72	&	0.92	&	4912	&	4.34E-04	&	0.93	&	0.24	\\
1E-05	&	1832	&	2.37E-04	&	1.40	&	0.76	&	12744	&	1.14E-04	&	1.40	&	1.93	\\
5E-06	&	2332	&	1.24E-04	&	2.69	&	0.94	&	20416	&	6.17E-05	&	1.31	&	0.38	\\
1E-06	&	3440	&	3.71E-05	&	3.10	&	0.75	&	47496	&	1.99E-05	&	1.34	&	1.63	\\

		\hline
&\multicolumn{4}{|c|}{$ k=2$, $ d=2$}& \multicolumn{4}{|c|}{$ k=2$, $d=3$}\\
		\hline
1E-04	&	747	&	7.15E-04	&		&		&	4779	&	2.28E-04	&		&		\\
5E-05	&	855	&	4.43E-04	&	3.54	&	0.69	&	6345	&	1.54E-04	&	1.38	&	0.57	\\
1E-05	&	1908	&	1.78E-04	&	1.14	&	0.57	&	14418	&	4.67E-05	&	1.46	&	0.74	\\
5E-06	&	2376	&	8.55E-05	&	3.34	&	1.06	&	19845	&	2.35E-05	&	2.14	&	0.99	\\
1E-06	&	4095	&	1.51E-05	&	3.18	&	1.08	&	37395	&	9.07E-06	&	1.50	&	0.60	\\
5E-07	&	4914	&	9.12E-06	&	2.77	&	0.73	&	50355	&	4.94E-06	&	2.04	&	0.88	\\

		\hline
&\multicolumn{4}{|c|}{$ k=3$, $ d=2$}& \multicolumn{4}{|c|}{$ k=3$, $d=3$}\\
		\hline
		5E-06	&	1952	&	6.88E-05	&		&			&	16384	&	1.02E-05	&		&		\\
		
		1E-06	&	3136	&	1.19E-05	&	3.70	&	1.09	&	29440	&	3.36E-06	&	1.90	&	0.69	\\
		
		5E-07	&	3696	&	5.79E-06	&	4.40	&	1.04	&	39616	&	1.62E-06	&	2.45	&	1.05	\\
		
		1E-07	&	4992	&	1.53E-06	&	4.43	&	0.83	&	59456	&	6.23E-07	&	2.36	&	0.60	\\
		
		5E-08	&	6288	&	6.19E-07	&	3.92	&	1.30	&	80832	&	3.53E-07	&	1.84	&	0.82	\\
		
		1E-08	&	9184	&	1.44E-07	&	3.85	&	0.91   &	129088	&	2.62E-08	&	5.56	&	1.62	\\

		\hline
		
	\end{tabular}
	\label{table:solid_d2}
\end{table}

\begin{figure}[htp]
	\begin{center}
%		\subfigure[]{\includegraphics[width=.42\textwidth]{figures/element_sbr_p1_n5-eps-converted-to}}
%		\subfigure[]{\includegraphics[width=.42\textwidth]{figures/error_sbr_p1_n5-eps-converted-to}}\\
		\subfigure[]{\includegraphics[width=.42\textwidth]{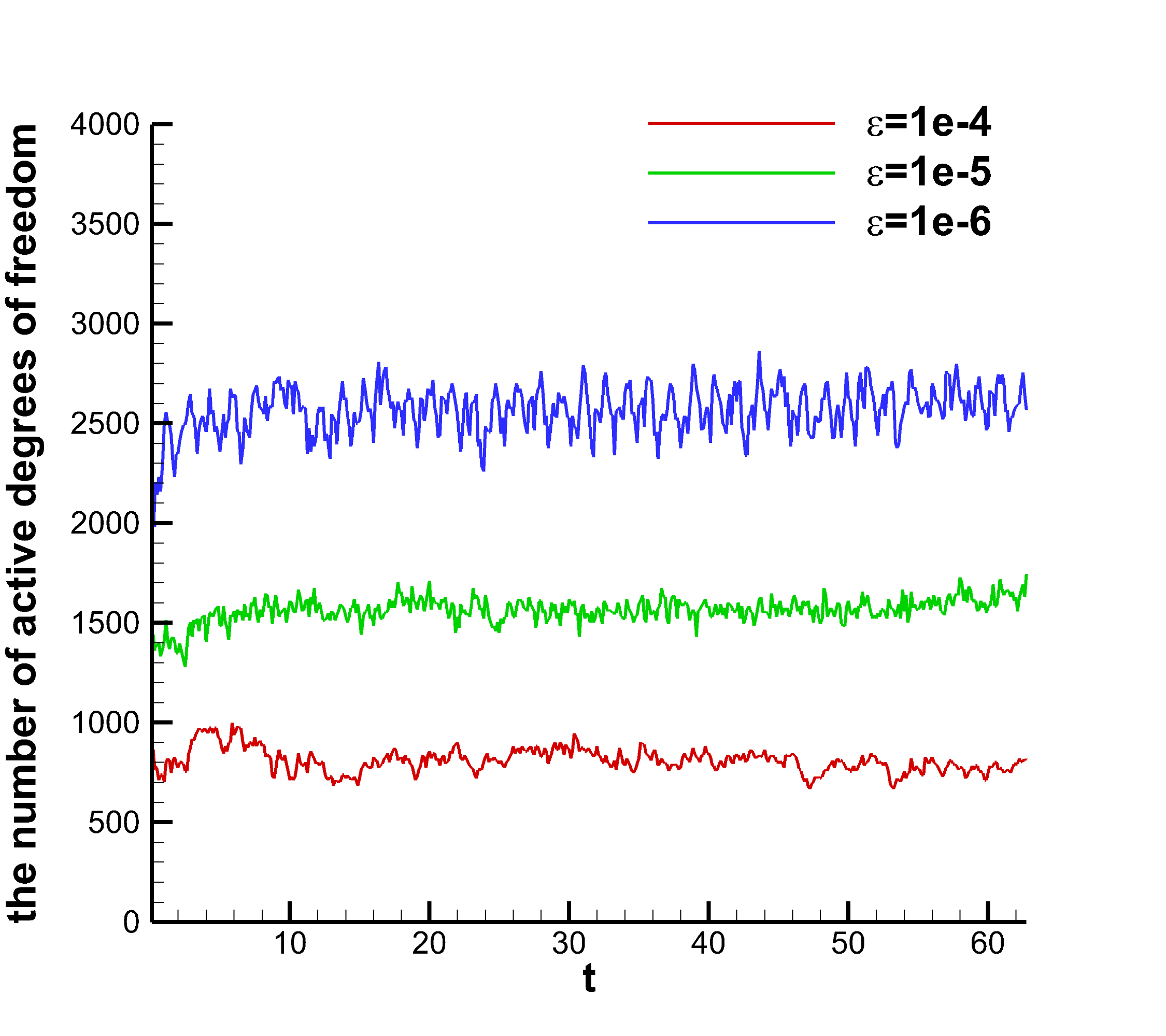}}
		\subfigure[]{\includegraphics[width=.42\textwidth]{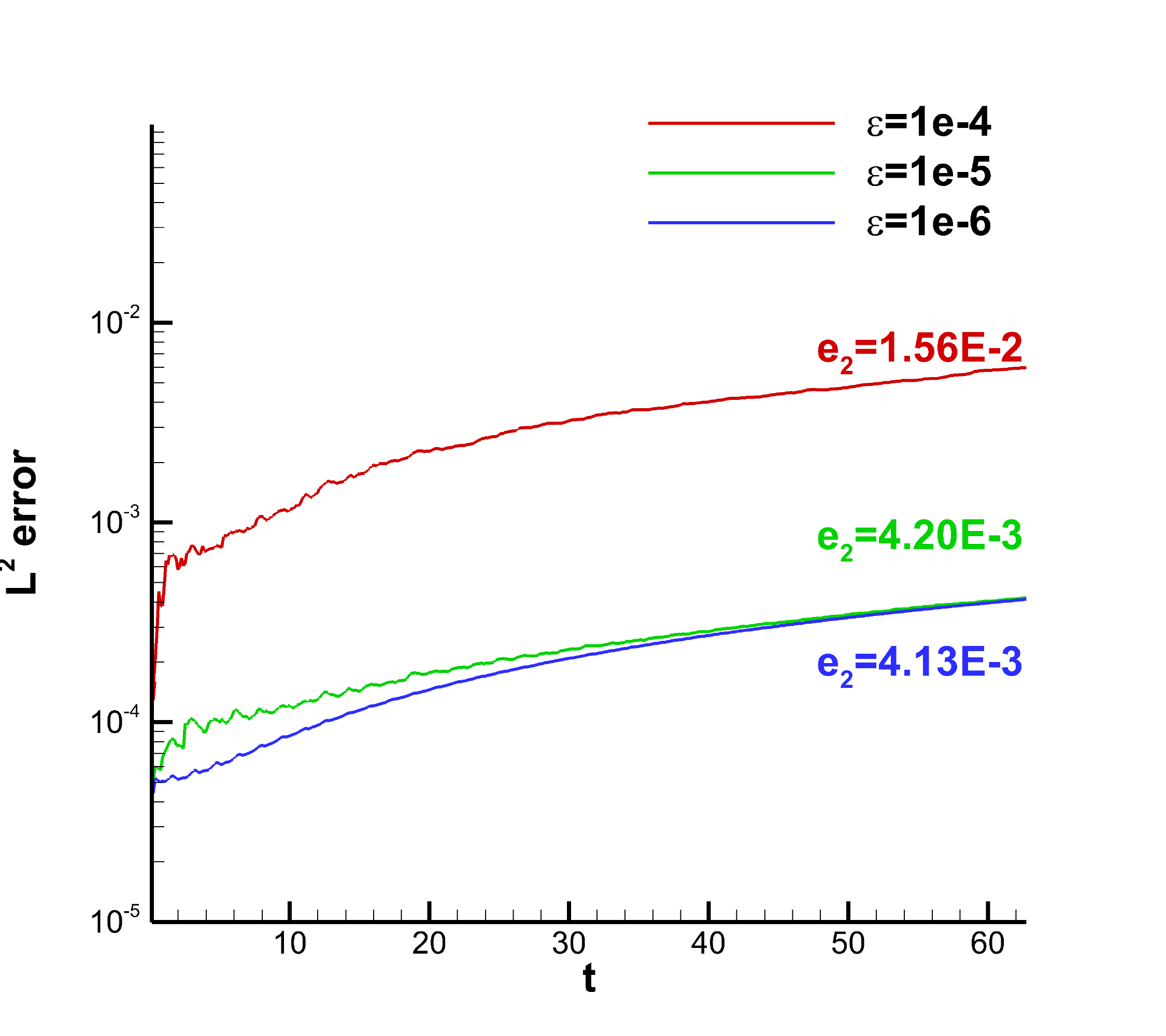}}\\
		\subfigure[]{\includegraphics[width=.42\textwidth]{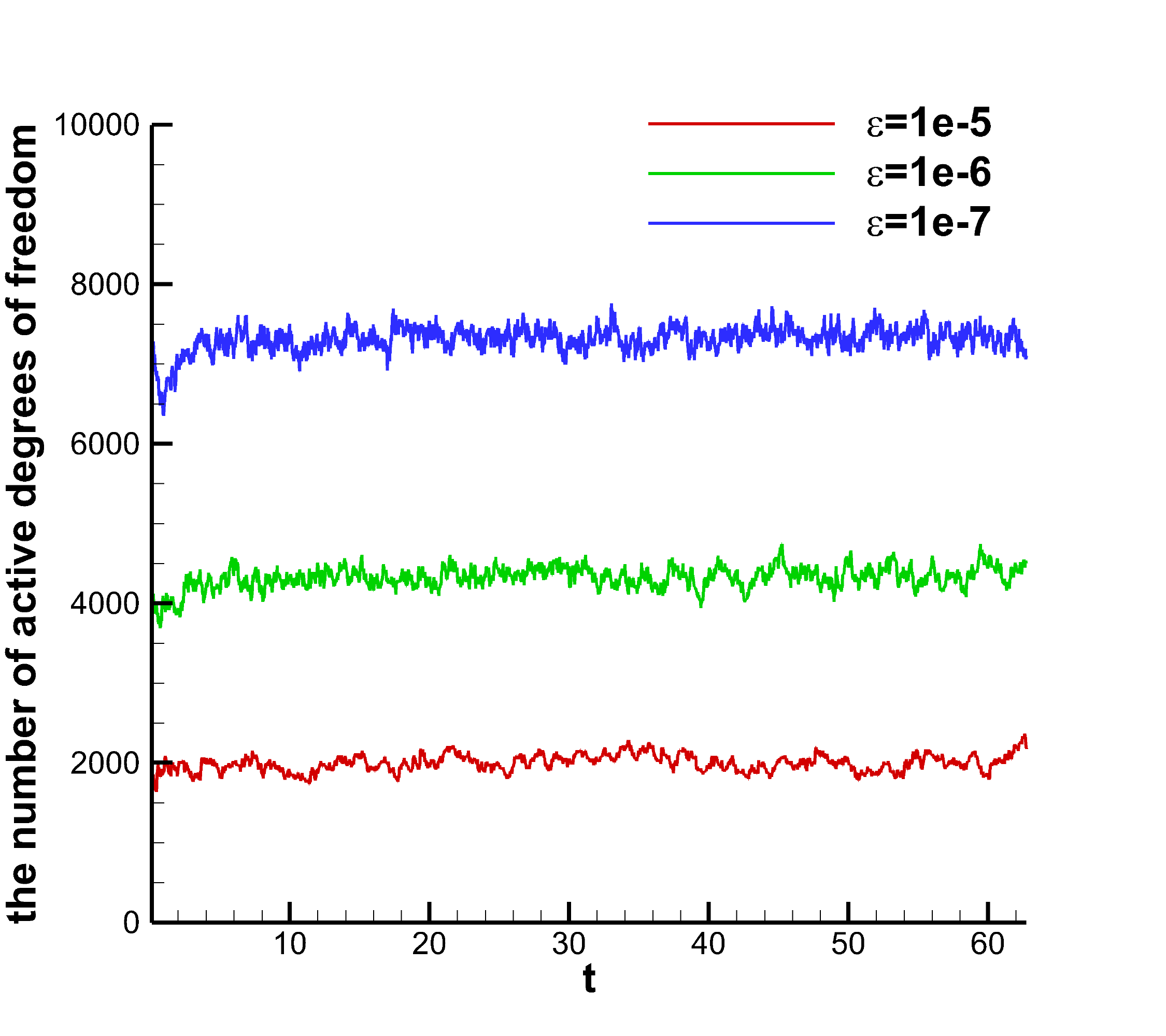}}
		\subfigure[]{\includegraphics[width=.42\textwidth]{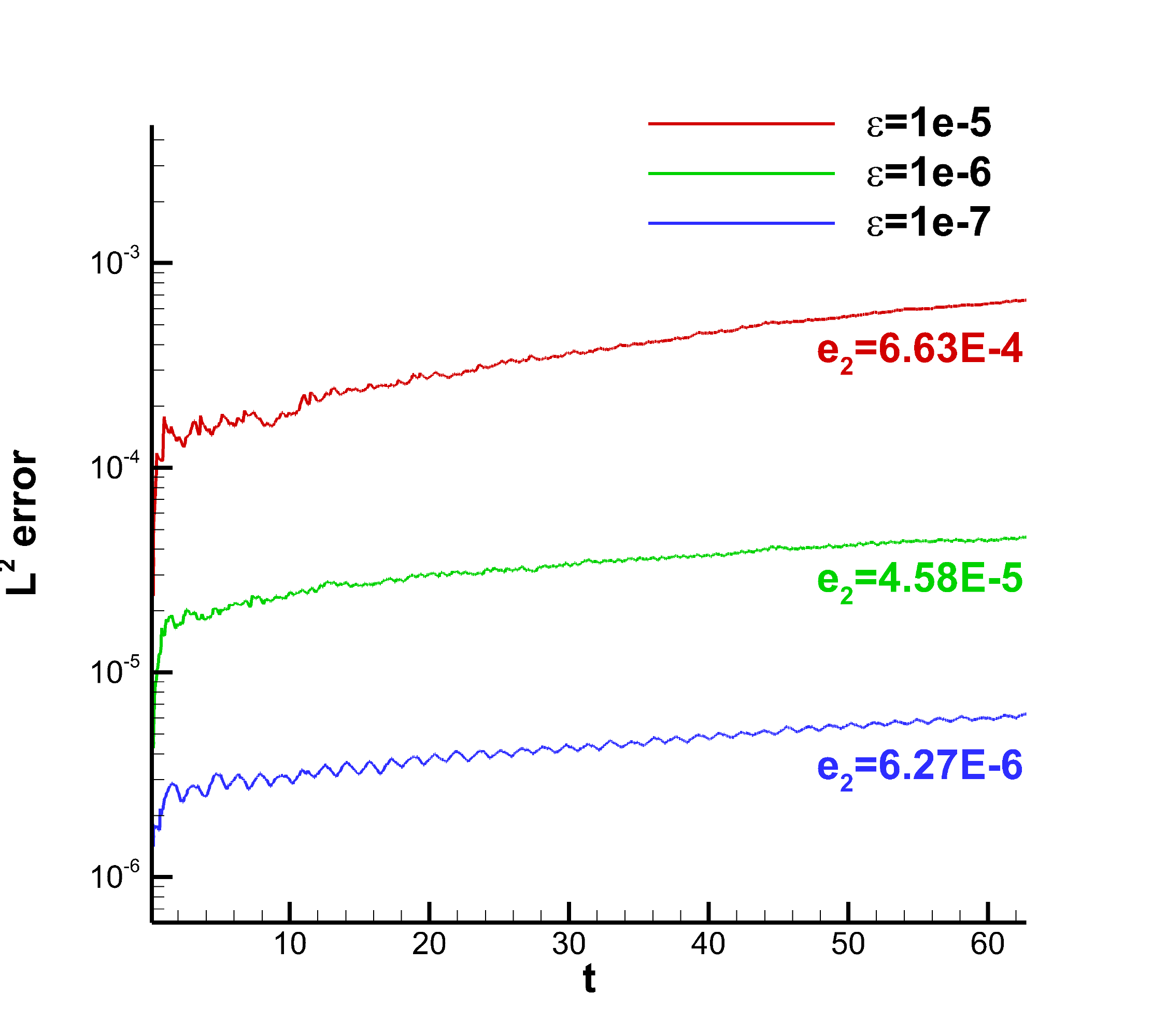}}\\
%		\subfigure[]{\includegraphics[width=.42\textwidth]{figures/element_sbr_p3_n5-eps-converted-to}}
%		\subfigure[]{\includegraphics[width=.42\textwidth]{figures/error_sbr_p3_n5-eps-converted-to}}
	\end{center}
	\caption{Example \ref{ex:rotation} with initial condition \eqref{eq:cosine}. Time history of the numbers of active degrees of freedom and $L^2$ errors with various $\varepsilon$. $k=2$.  (a-b) $N=5$.  (c-d) $N=7$. }
	\label{fig:sbr}
\end{figure}

We also  use this example to compare the performance of the scheme with different configurations of $N$ and $\varepsilon$. We let $d=2$, $k=2$ and compute the solutions up to ten periods and plot the time evolution of $L^2$ errors and the number of active degrees of freedom in Figure \ref{fig:sbr}. In particular, we compare maximum mesh level $N=5$ and $N=7,$ and run the simulations with three different values of $\varepsilon$. Since the cosine bell keeps its initial profile as time evolves,   the degrees of freedom to resolve the solution for a fixed accuracy threshold should remain the same. We observe that if an excessively small $\varepsilon$ is taken, the used degrees of freedom will increase, but the error may not decrease much, see Figure \ref{fig:sbr} (a-b). %It is worth pointing out that the performance of the adaptive MRA DG scheme depends closely on the choice of  the refinement parameter $\varepsilon$, the maximum mesh level $N$ and degree of polynomials $k$. 
This shows the importance of the choice of   $\varepsilon$ and $N$ for the computational efficiency of the scheme.

\begin{figure}[htp]
	\begin{center}
		\subfigure[]{\includegraphics[width=.42\textwidth]{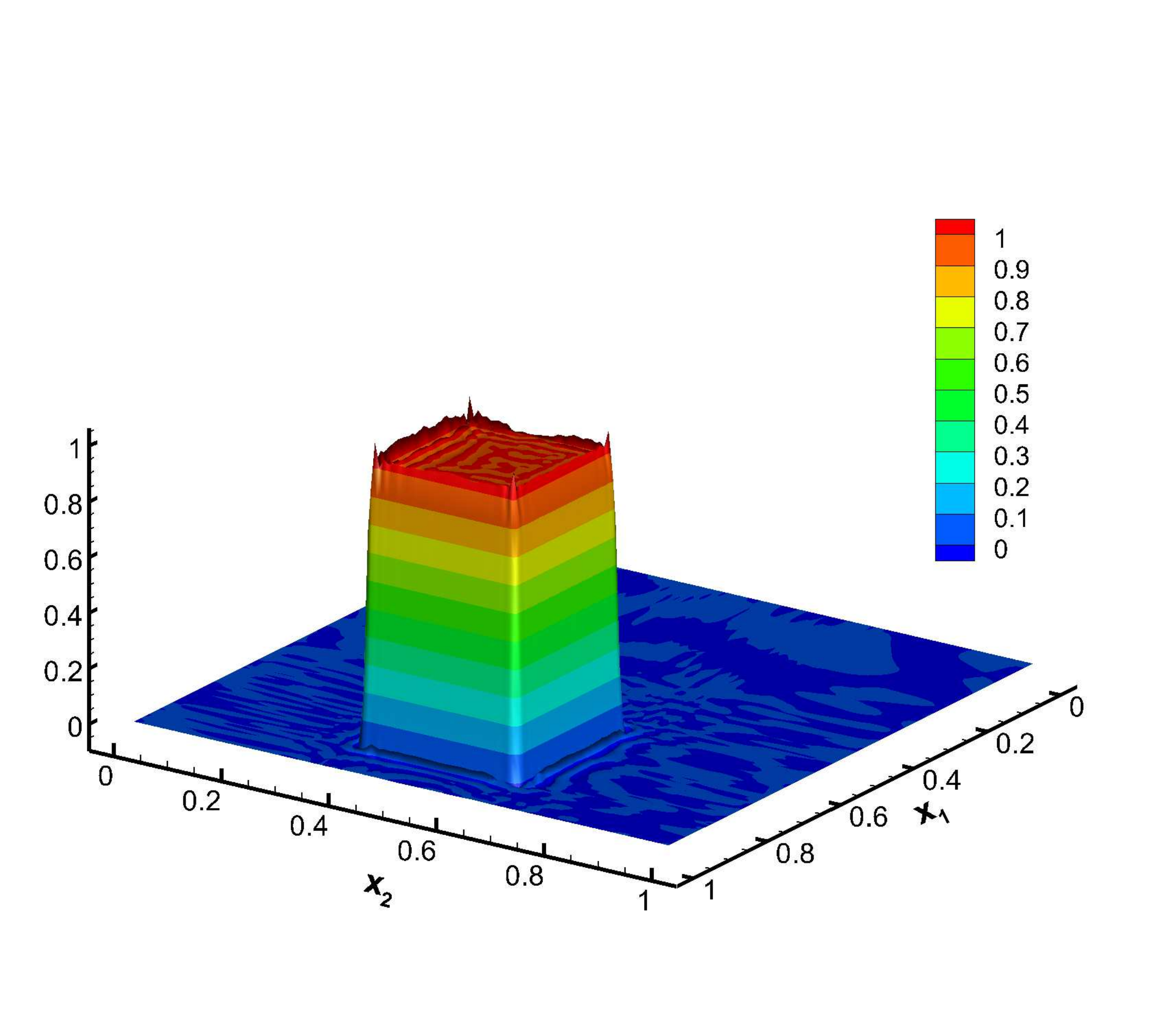}}
		\subfigure[]{\includegraphics[width=.42\textwidth]{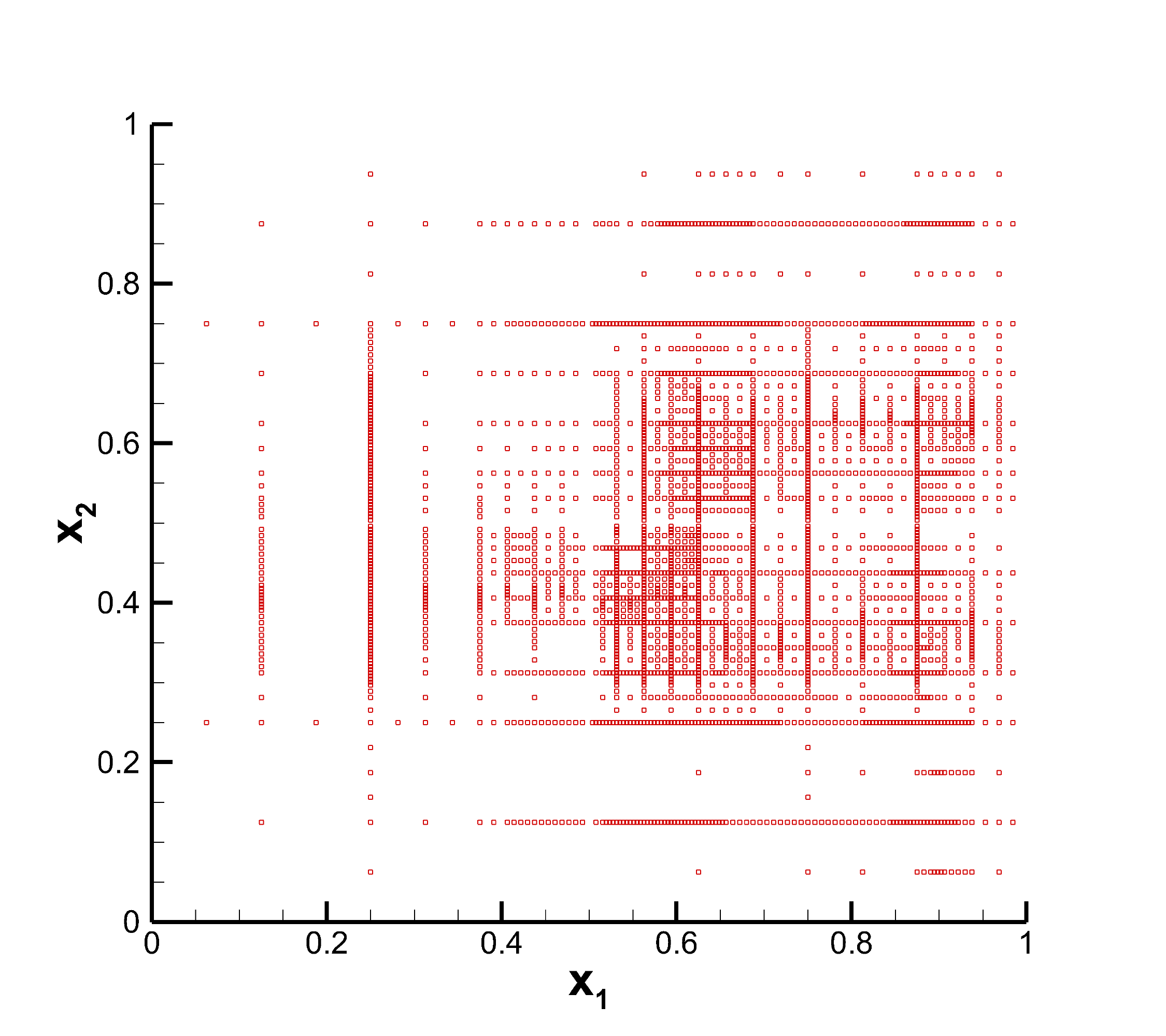}}\\
\subfigure[]{\includegraphics[width=.42\textwidth]{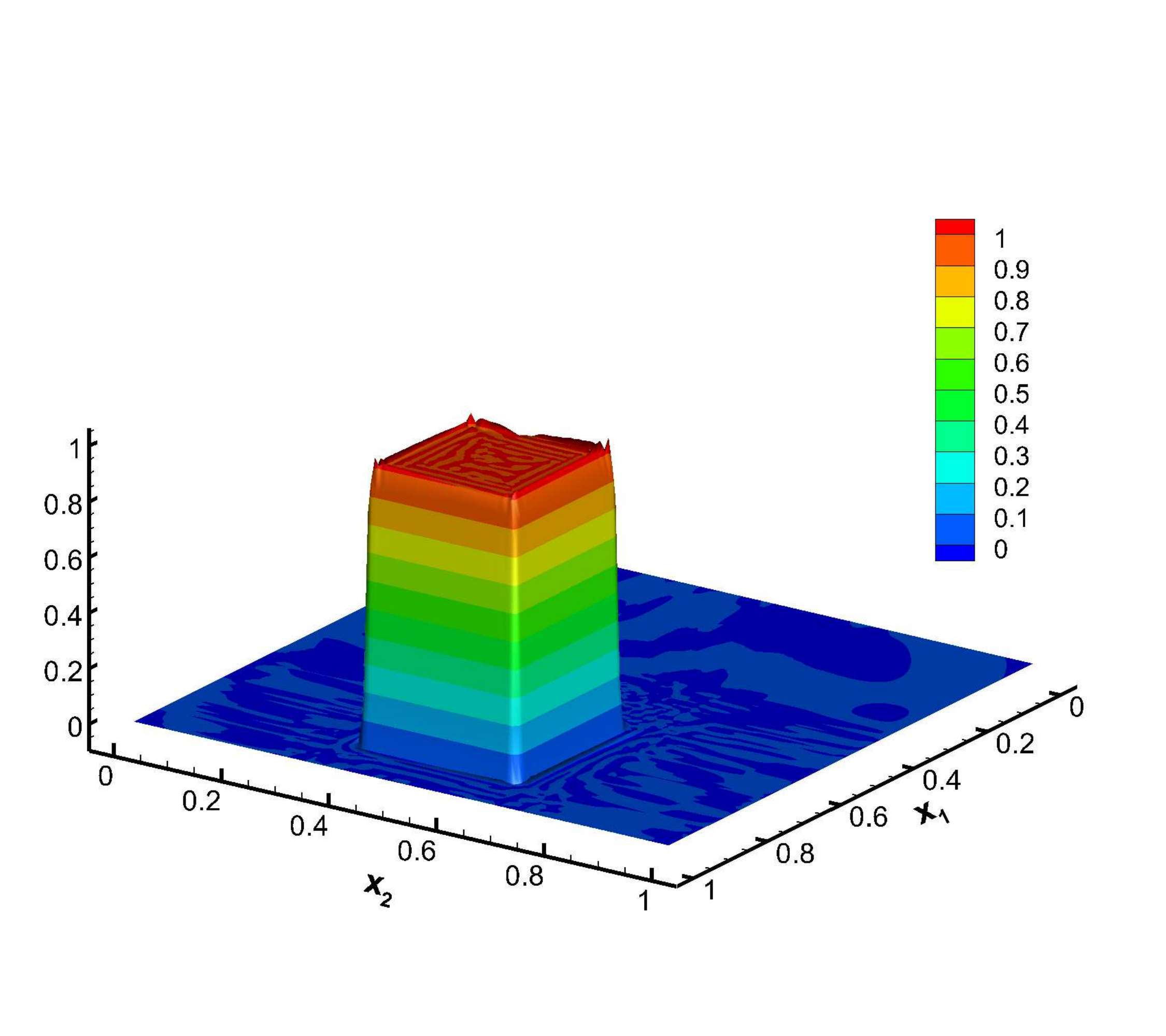}}
\subfigure[]{\includegraphics[width=.42\textwidth]{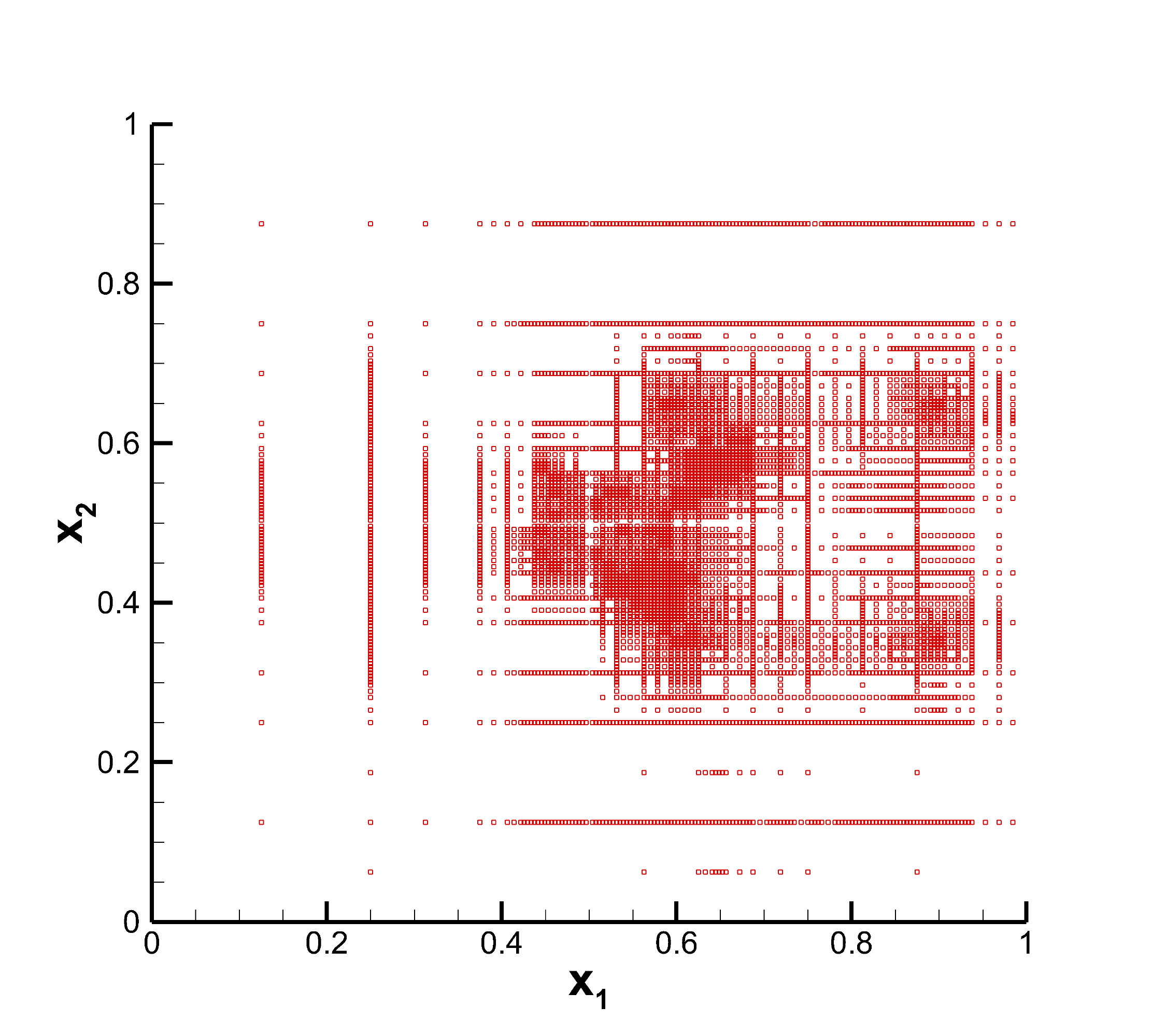}}
	\end{center}
	\caption{Example \ref{ex:rotation} with initial condition \eqref{eq:sbr_discontinuous}.  (a) Numerical solution with $L^1$ norm based criteria. (b) Active elements with $L^1$ norm based criteria. (c) Numerical solution with $L^2$ norm based criteria. (d) Active elements with $L^2$ norm based criteria. $T=2\pi$. $d=2$. $N=7$. $k=3$. $\varepsilon=10^{-5}$.}
	\label{fig:sbr_n7_dis}
\end{figure}

We then consider following discontinuous initial condition:

\begin{equation}
\label{eq:sbr_discontinuous}
 u(0,\bx)=\left\{\begin{array}{ll}\displaystyle 1& (x_1,x_2)\in[\frac34-\frac{\sqrt{2}}{10},\frac34+\frac{\sqrt{2}}{10}]\times[\frac12-\frac{\sqrt{2}}{10},\frac12+\frac{\sqrt{2}}{10}],\\[2mm]
0& \text{otherwise},\end{array}\right.
\end{equation}
when $d=2.$ In the simulation, we set $N=7$, $k=3$, $\varepsilon=10^{-5}$ and consider both $L^1$ and $L^2$ norm based criteria. In Figure \ref{fig:sbr_n7_dis}, we report the numerical solutions and the associated active elements at $T=2\pi$.
Similar to the previous example,  elements clusters towards the discontinuities and the scheme with both criteria is able to well resolve the discontinuities. %Furthermore, the  $L^1$ norm based criteria result in less degrees of freedom than the  $L^2$ norm based criteria. 
However, more severe localized numerical oscillations are observed when compared with the previous example. 

%\begin{align}
%R_{\epsilon_l}&=\frac{\log(e_{l-1}/e_l)}{\log(\epsilon_{l-1}/{\epsilon_l})}\\
%R_{\text{DOF}_l}&=\frac{\log(e_{l-1}/e_l)}{\log(\text{DOF}_l/\text{DOF}_{l-1})}, 
%\end{align} 
%where $e_l$ is the standard $L^2$ error with refinement parameter $\epsilon_l$, and $\text{DOF}_l$ is the associated number of active degrees of freedom at final time.
%

%\begin{figure}[htp]
%	\begin{center}
%		\subfigure[]{\includegraphics[width=.42\textwidth]{figures/element_sbr_p1_n7}}
%		\subfigure[]{\includegraphics[width=.42\textwidth]{figures/error_sbr_p1_n7}}\\
%		\subfigure[]{\includegraphics[width=.42\textwidth]{figures/element_sbr_p2_n7}}
%		\subfigure[]{\includegraphics[width=.42\textwidth]{figures/error_sbr_p2_n7}}\\
%		\subfigure[]{\includegraphics[width=.42\textwidth]{figures/element_sbr_p3_n7}}
%		\subfigure[]{\includegraphics[width=.42\textwidth]{figures/error_sbr_p3_n7}}
%	\end{center}
%	\caption{Example \ref{ex:rotation}. Solid body rotation test. Time history of the numbers of active degrees of freedom and $L^2$ errors. $k=1$ (a-b); $k=2$ (c-d); $k=3$ (e-f). $N=7$. }
%	\label{fig:sbr_n7}
%\end{figure}

\begin{exa}[Deformational flow]
	\label{ex:deformational}
	We consider  two-dimensional deformational flow with velocity field
	$$\ba=(\sin^2(\pi x_1)\sin(2\pi x_2)g(t),-\sin^2(\pi x_2)\sin(2\pi x_1)g(t)),$$
	where $g(t)=\cos(\pi t/T)$ with $T=1.5$. 
\end{exa}
First, we choose the cosine bell \eqref{eq:cosine}  as the initial condition, but with $\bx_c=(0.65,0.5)$ and $b=0.35$. The cosine bell deforms into a crescent shape at $t=T/2$, then goes back to its initial state at $t=T$ as the flow reverses. 
We perform a similar convergence study as in the previous two examples, which is summarized in Table \ref{table:defor_d2}.  We observe similar convergence pattern that  the rate $R_{\text{DOF}}$ is slightly smaller than $k+1$ and $R_\varepsilon$ is close to 1.

\begin{table}[htp]
	\caption{Example \ref{ex:deformational} with initial condition \eqref{eq:cosine}. Numerical error and convergence rate. $N=7$. $T=1.5$. $d=2$.
	}
	%\vspace{2 mm}
	\centering
	\begin{tabular}{|c| c c c c|}
		\hline
		$\varepsilon$ &  DOF& $L^2$ error &  $R_{\text{DOF}}$& $R_\varepsilon$ \\
		\hline
		
	&	\multicolumn{4}{|c|}{$ k=1$}   \\
		\hline

1E-03	&	244	&	1.52E-02	&		&		\\
5E-04	&	372	&	7.94E-03	&	1.53	&	0.93	\\
1E-04	&	945	&	1.25E-03	&	1.98	&	1.15	\\
5E-05	&	1248	&	1.00E-03	&	8.01	&	0.32	\\
1E-05	&	2608	&	1.84E-04	&	2.30	&	1.05	\\
5E-06	&	3508	&	9.96E-05	&	2.07	&	0.89	\\
1E-06	&	5596	&	3.81E-05	&	2.06	&	0.60	\\
		\hline
		
	&	\multicolumn{4}{|c|}{$ k=2$}   \\
		\hline
5E-05	&	1143	&	5.41E-04	&		&		\\
1E-05	&	2043	&	1.15E-04	&	2.67	&	0.96	\\
5E-06	&	2736	&	6.91E-05	&	1.74	&	0.73	\\
1E-06	&	4842	&	1.24E-05	&	3.00	&	1.07	\\
5E-07	&	5994	&	8.29E-06	&	1.90	&	0.59	\\
1E-07	&	9045	&	1.74E-06	&	3.79	&	0.97	\\
5E-08	&	11142	&	1.08E-06	&	2.28	&	0.69	\\
		\hline
		
	&	\multicolumn{4}{|c|}{$ k=3$}   \\
		\hline
5E-05	&	1056	&	3.65E-04	&		&		\\
1E-05	&	2048	&	8.85E-05	&	2.14	&	0.88	\\
5E-06	&	2320	&	5.41E-05	&	3.94	&	0.71	\\
1E-06	&	3904	&	1.45E-05	&	2.53	&	0.82	\\
5E-07	&	4480	&	6.32E-06	&	6.02	&	1.20	\\
1E-07	&	6224	&	1.30E-06	&	4.80	&	0.98	\\
5E-08	&	7680	&	5.84E-07	&	3.82	&	1.16	\\
		\hline
		
	\end{tabular}
	\label{table:defor_d2}
\end{table}

\begin{figure}[htp]
	\begin{center}
		\subfigure[]{\includegraphics[width=.42\textwidth]{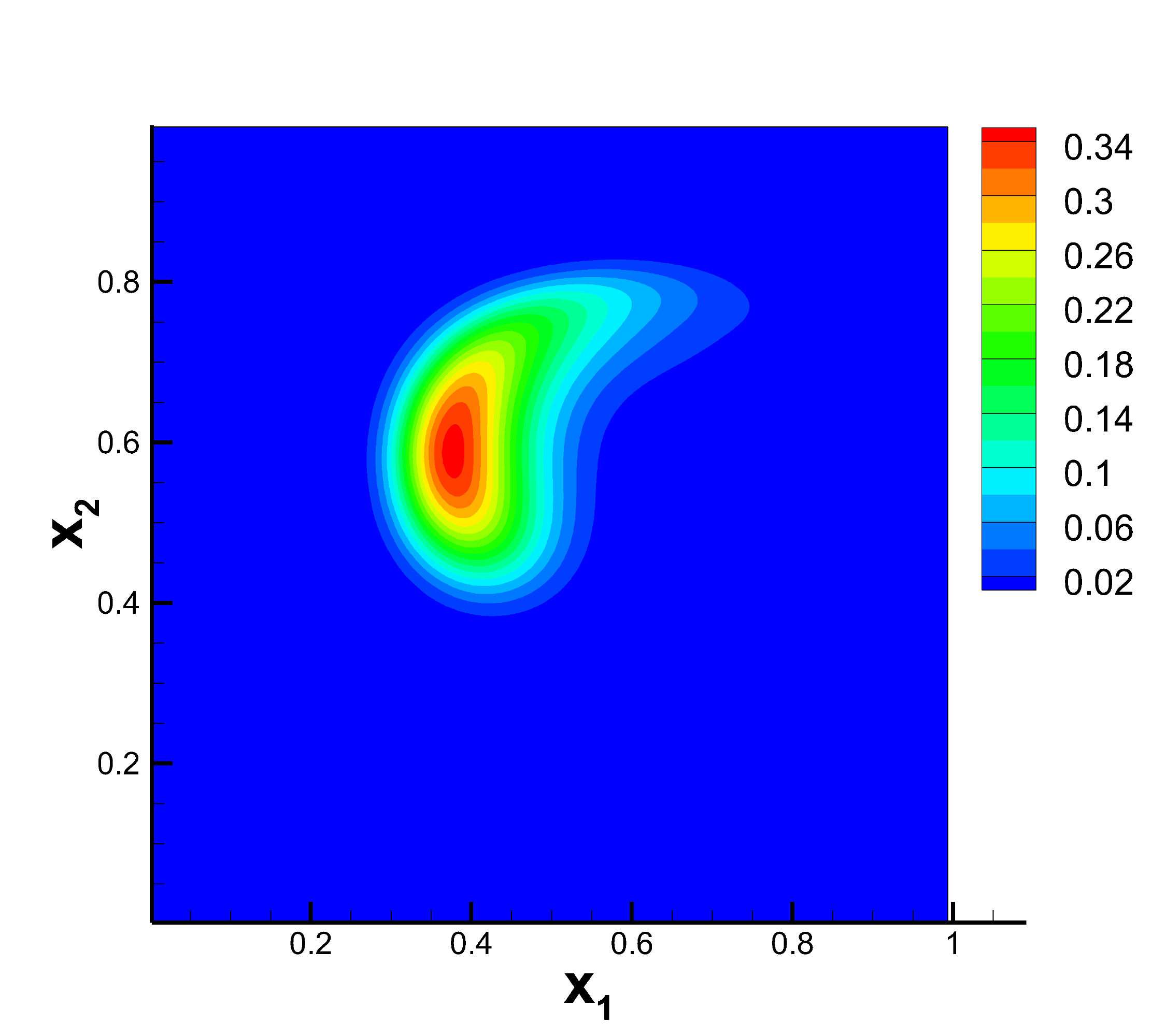}}
		\subfigure[]{\includegraphics[width=.42\textwidth]{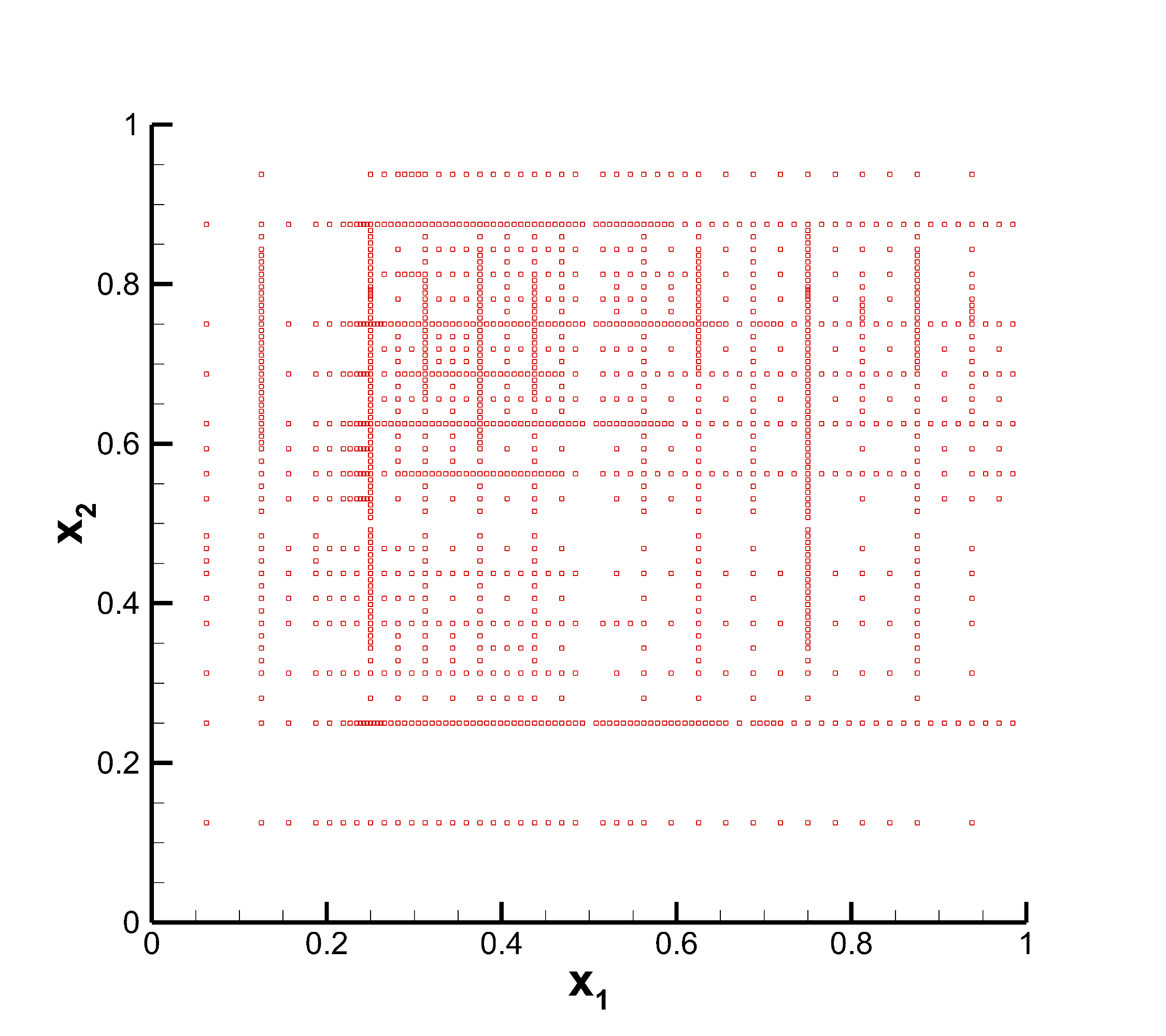}}\\
		\subfigure[]{\includegraphics[width=.42\textwidth]{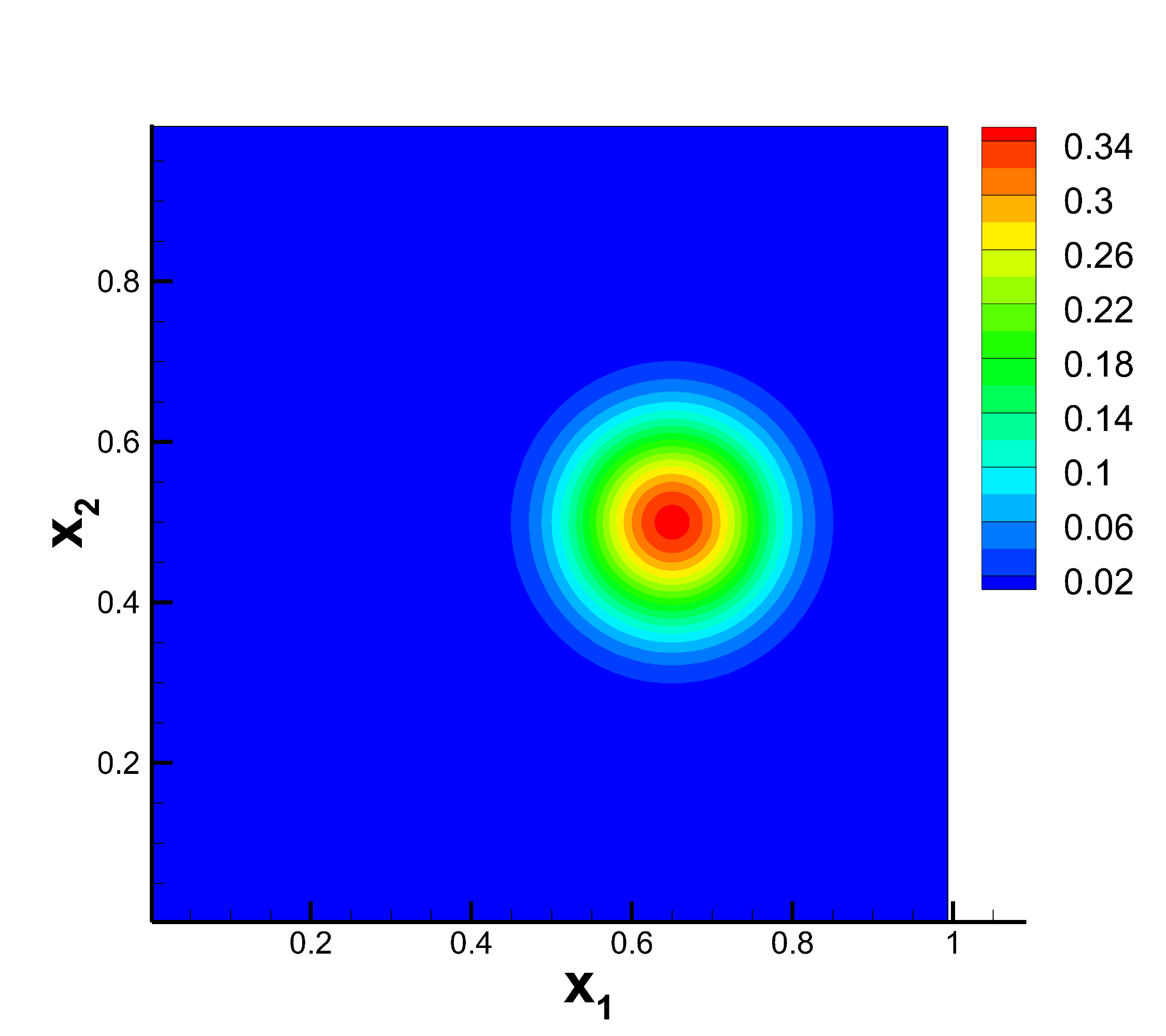}}
		\subfigure[]{\includegraphics[width=.42\textwidth]{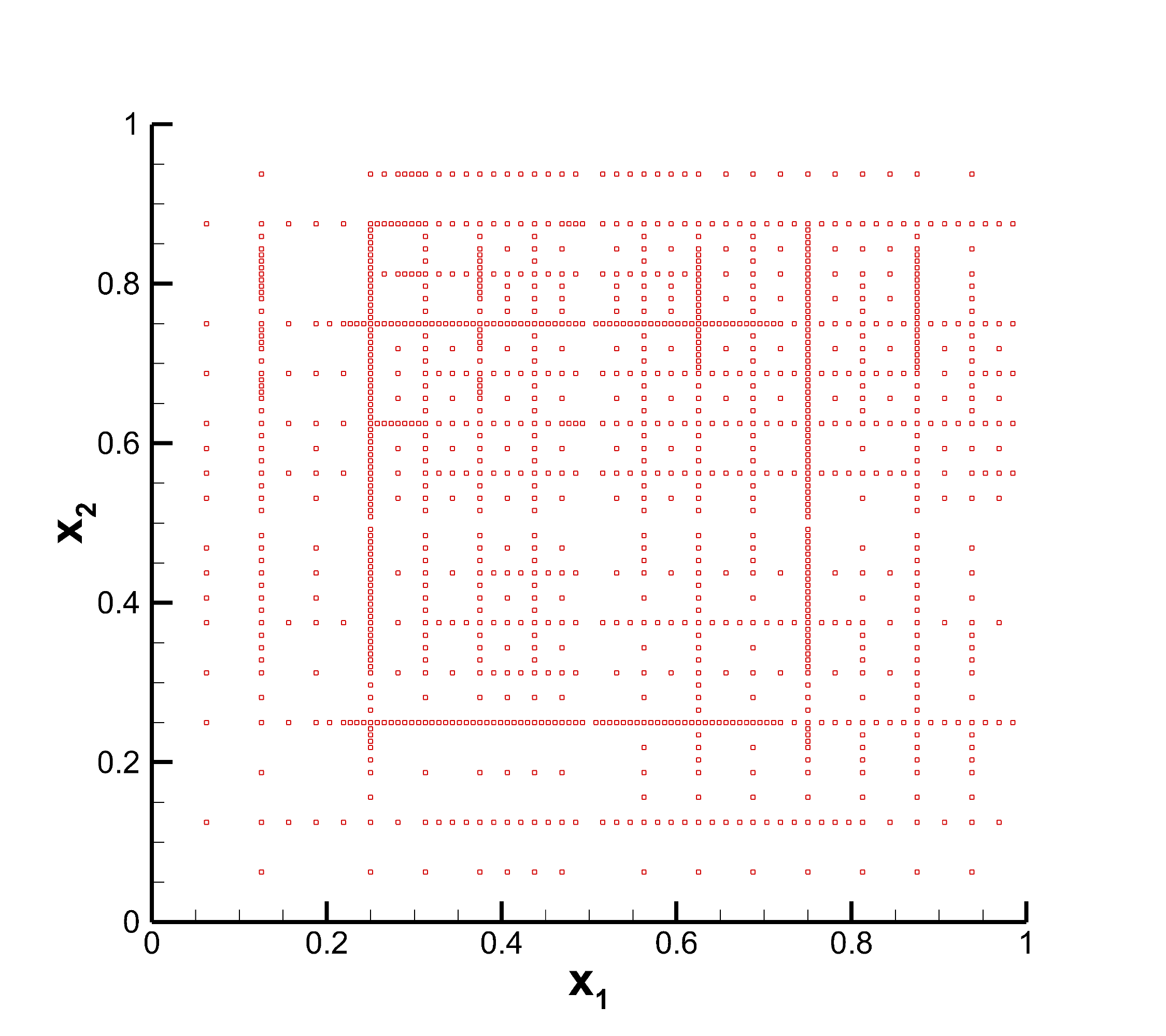}}
	\end{center}
	\caption{Example \ref{ex:deformational} with initial condition \eqref{eq:cosine}. $N=7$. $k=3$. $\varepsilon=10^{-7}$.  (a-b) $t=T/2$. (c-d) $t=T.$ }
	\label{fig:defo_n7}
\end{figure}

In Figure \ref{fig:defo_n7}, we present the contour plots and the associated active elements of the numerical solutions computed with   $N=7,\,k=3,\,\varepsilon=10^{-7}$ at $t=T/2$ when the shape of the bell is severely deformed, and $t=T$ when the solution is recovered into its initial state.    The elements  tend to cluster where the solution deforms as expected, and the shape of the cosine bell is well recovered at $t=T.$

\begin{figure}[htp]
	\begin{center}
		\subfigure[]{\includegraphics[width=.42\textwidth]{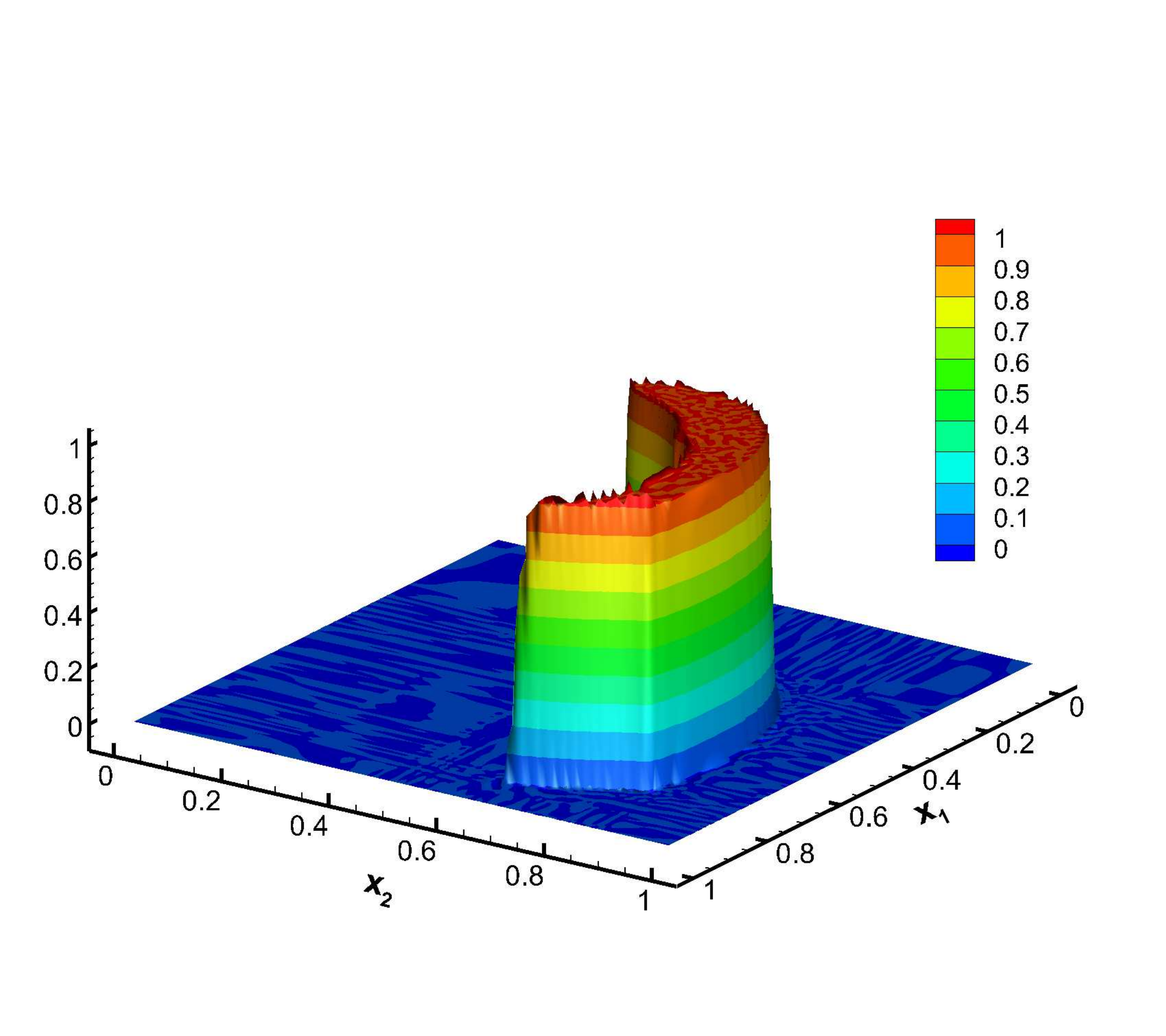}}
		\subfigure[]{\includegraphics[width=.42\textwidth]{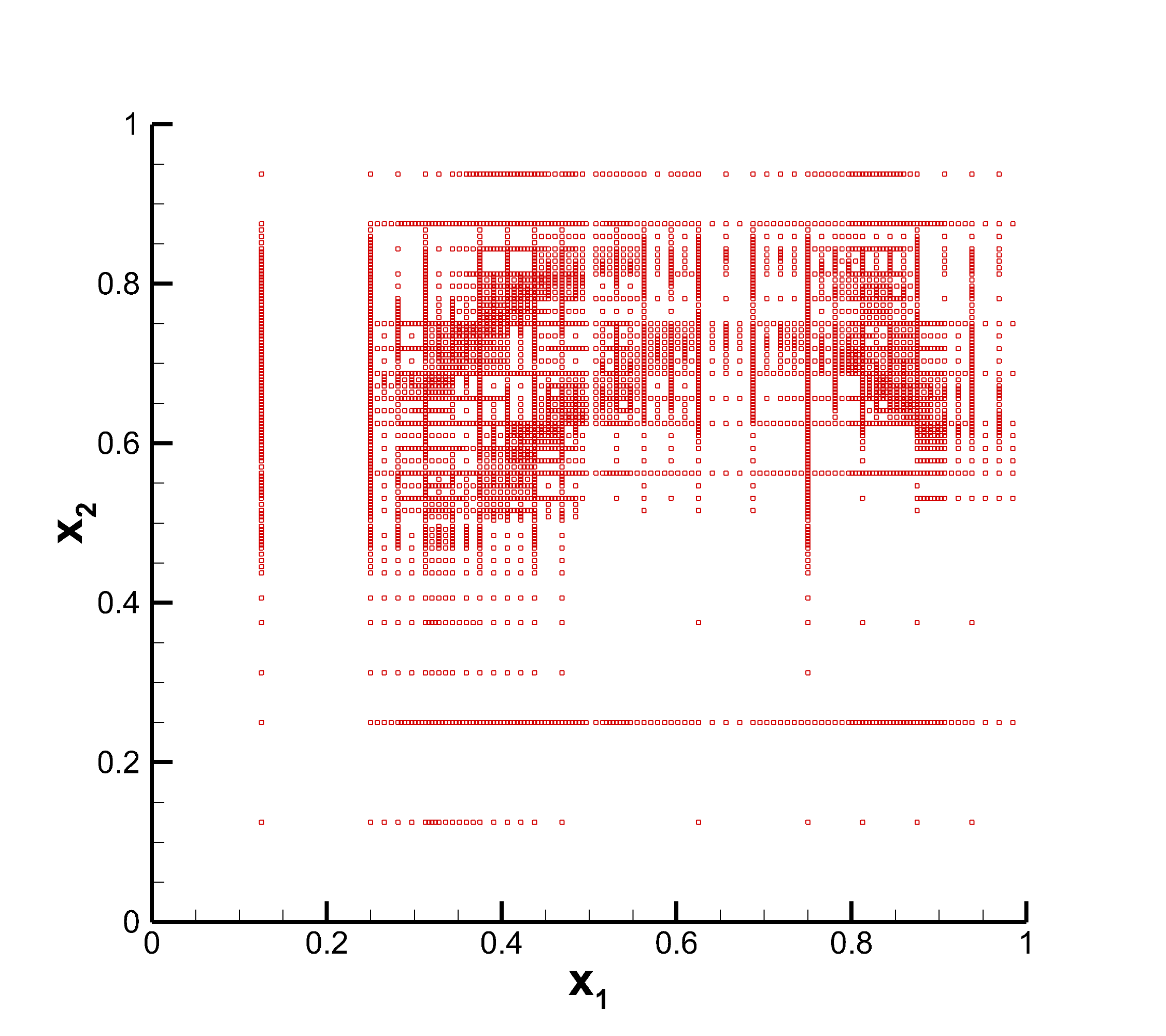}}\\
		\subfigure[]{\includegraphics[width=.42\textwidth]{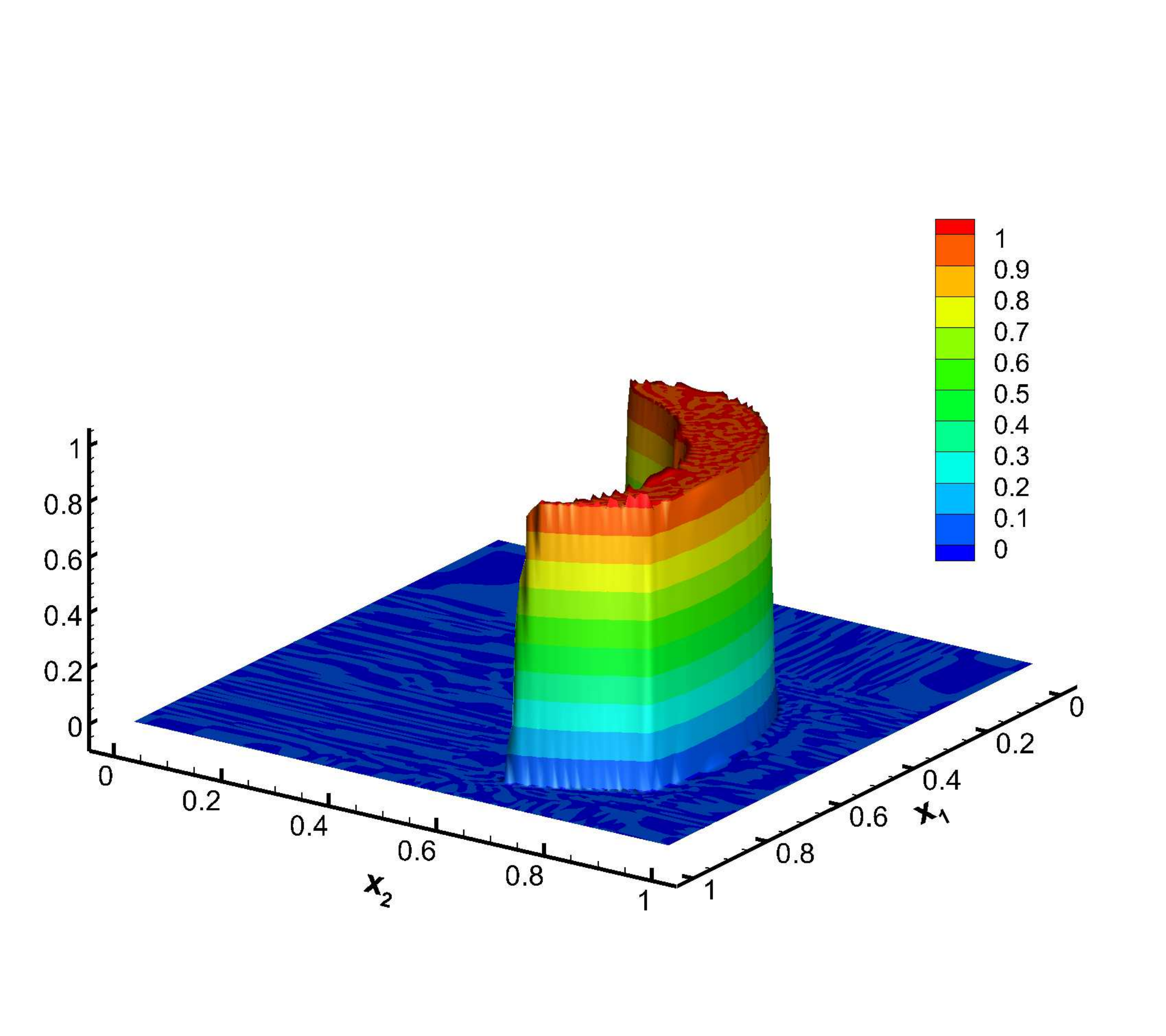}}
		\subfigure[]{\includegraphics[width=.42\textwidth]{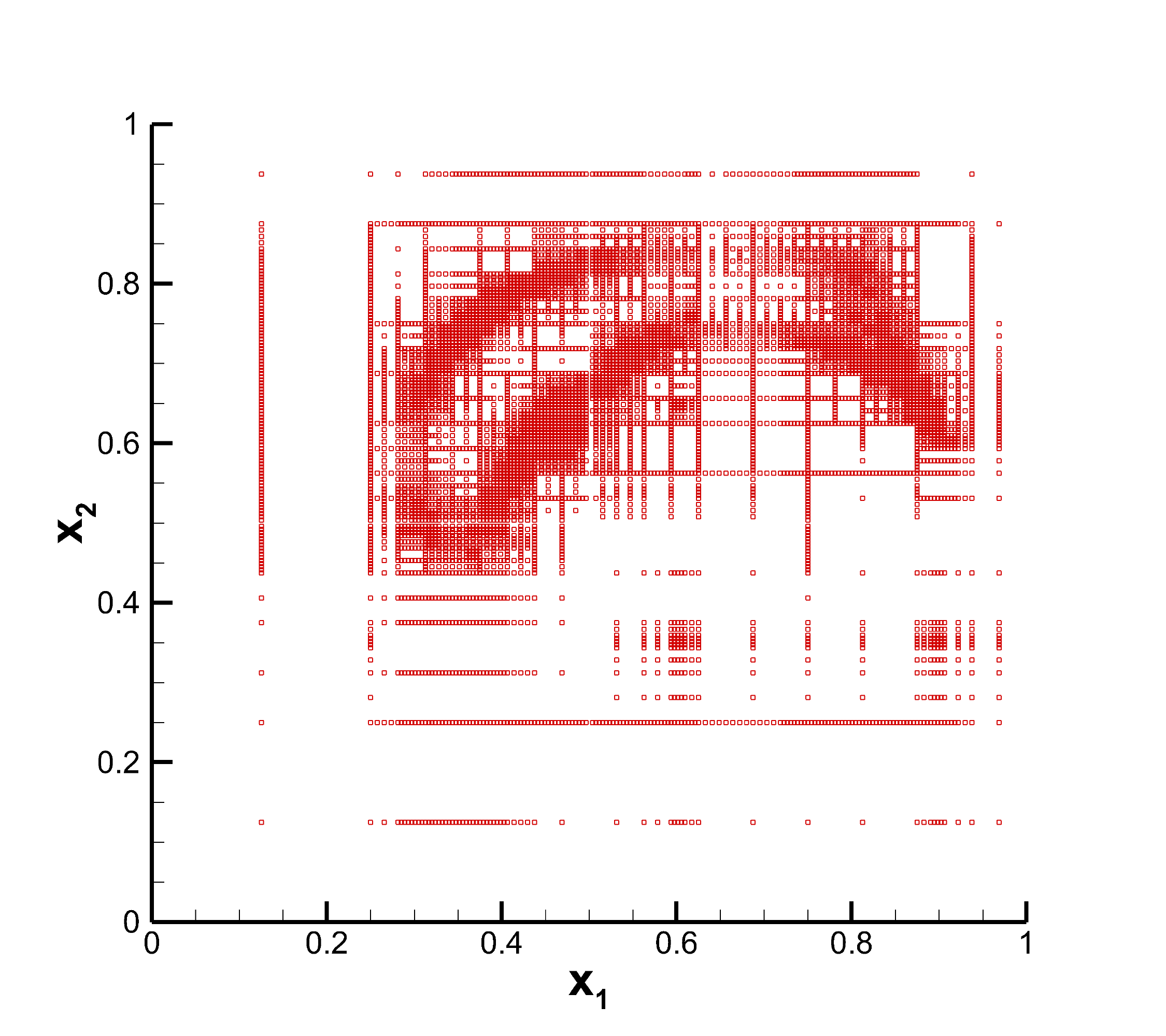}}
	\end{center}
	\caption{Example \ref{ex:deformational} with initial condition \eqref{eq:sbr_discontinuous}. $t=T/2$. $N=7$. $k=3$. $\varepsilon=10^{-5}$. (a-b) $L^1$ norm based criteria.  (c-d) $L^2$ norm based criteria.}
	\label{fig:defo_dis_t0}
\end{figure}

\begin{figure}[htp]
	\begin{center}
		\subfigure[]{\includegraphics[width=.42\textwidth]{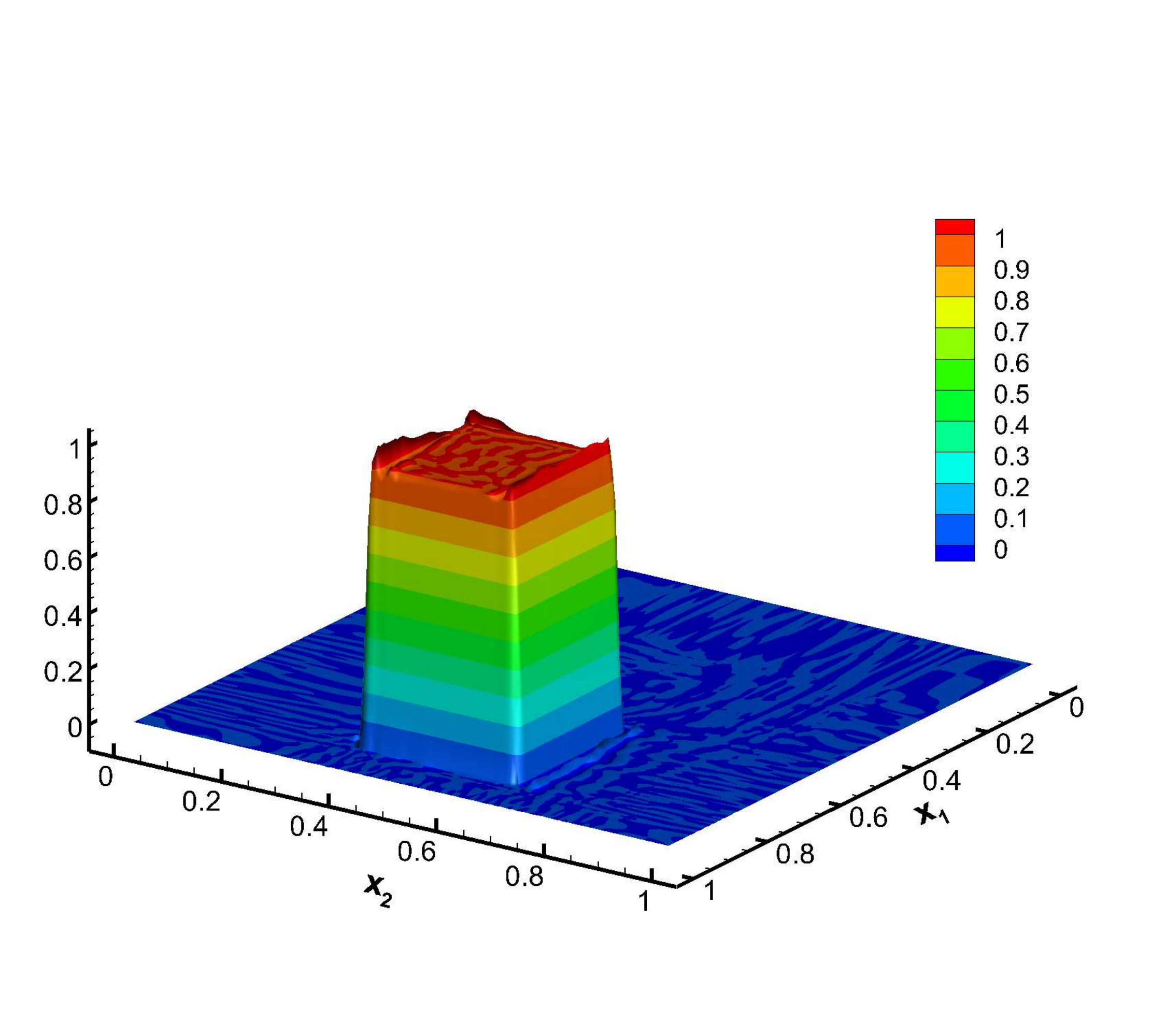}}
		\subfigure[]{\includegraphics[width=.42\textwidth]{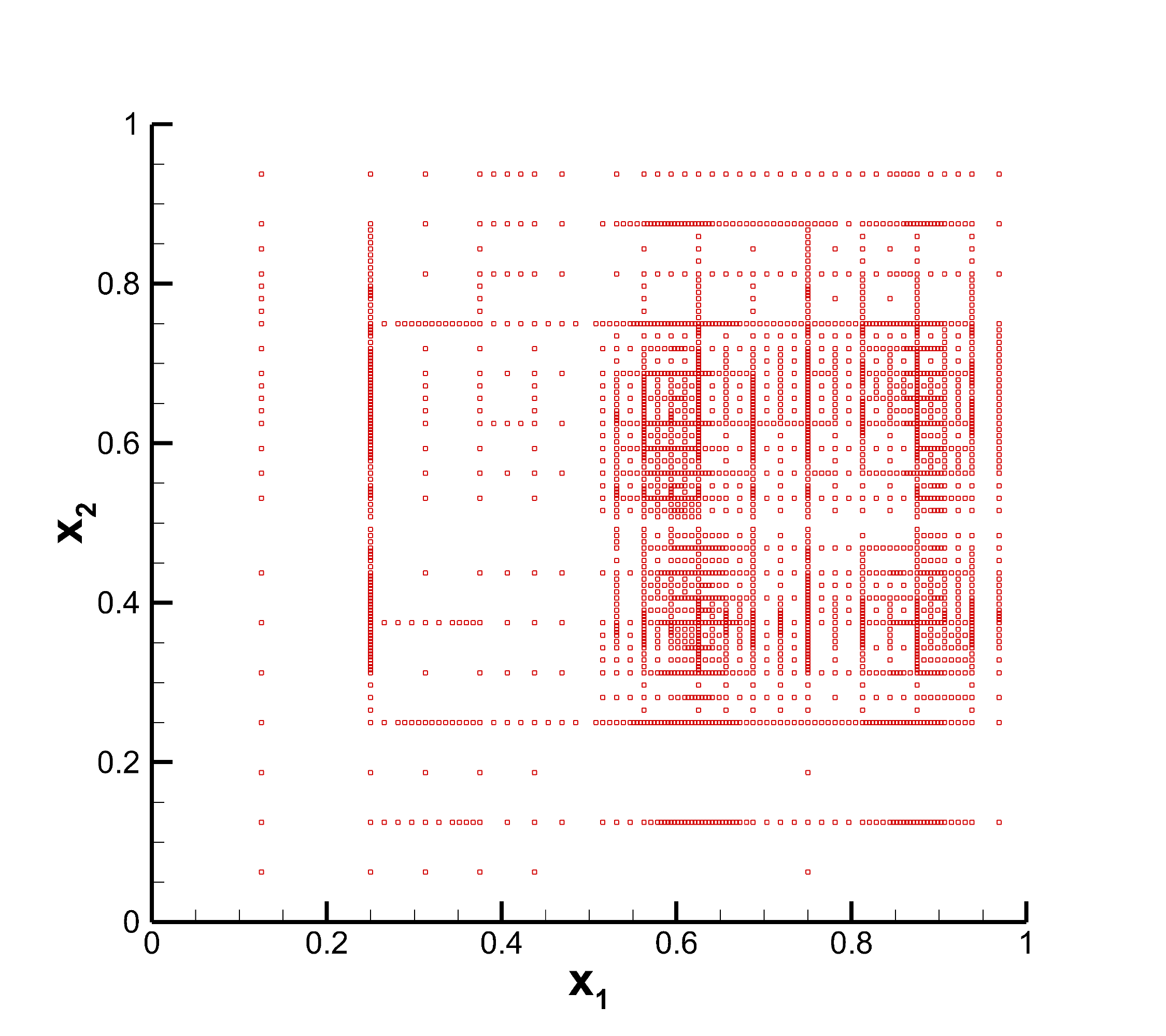}}\\
		\subfigure[]{\includegraphics[width=.42\textwidth]{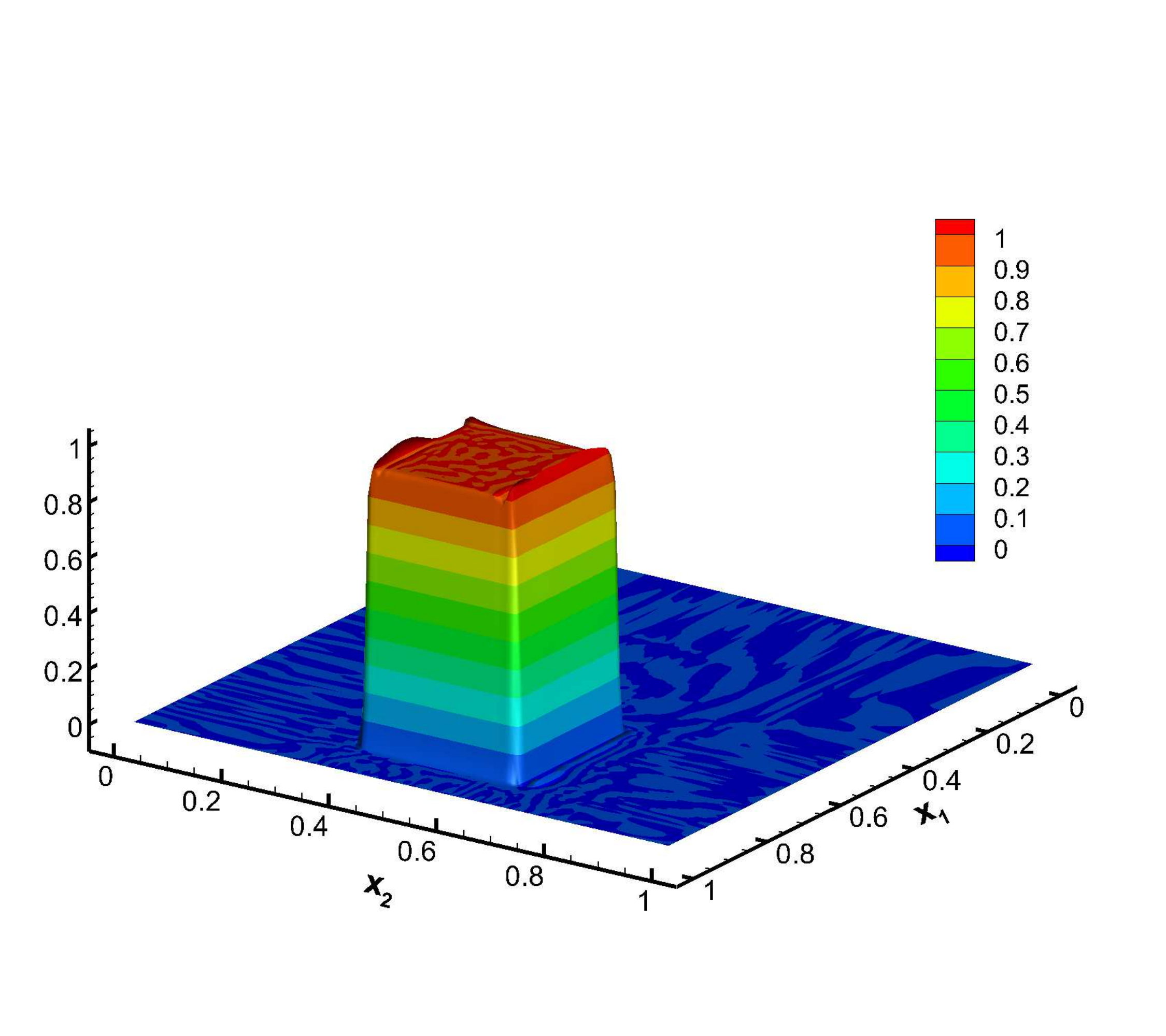}}
		\subfigure[]{\includegraphics[width=.42\textwidth]{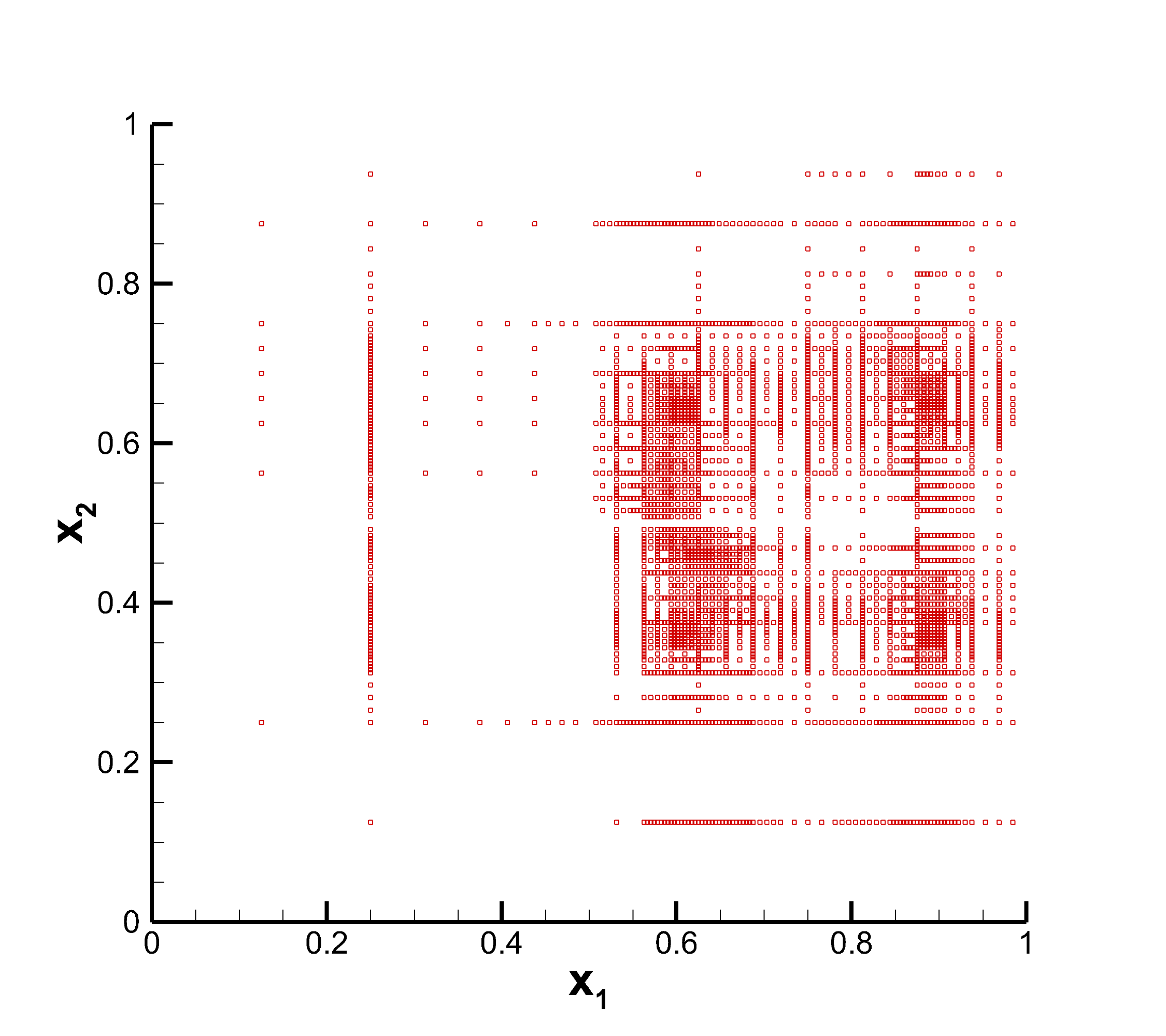}}
	\end{center}
	\caption{Example \ref{ex:deformational} with initial condition \eqref{eq:sbr_discontinuous}. $t=T$. $N=7$. $k=3$. $\varepsilon=10^{-5}$. (a-b) $L^1$ norm based criteria.  (c-d) $L^2$ norm based criteria.}
	\label{fig:defo_dis_t1}
\end{figure}

 We also consider the discontinuous initial condition \eqref{eq:sbr_discontinuous}, and use both $L^1$ and $L^2$ based refinement/coarsening criteria with $N=7$, $k=3$ and $\varepsilon=10^{-5}$. The numerical solutions and the associated active elements at $t=T/2$ and $t=T$ are plotted in Figures \ref{fig:defo_dis_t0} and \ref{fig:defo_dis_t1}, respectively. For this challenging test, the numerical solution tends to be more oscillatory, again because no limiting mechanism is present in the scheme. The $L^2$ norm based criteria generate less oscillatory profiles but use more elements for computation.

\section{Vlasov-Poisson simulations}
\label{sec:kinetic}

In this section, we apply the adaptive multiresolution DG methods to solve the kinetic transport equation. Here, we consider   the VP system, which is a fundamental model in plasma simulation. The solution is known to develop filamentation (fine structures) in the phase space. Therefore, it is a good test problem for the adaptive algorithm. For simplicity, we restrict our attention to two dimensional cases, but comment that the algorithm can be readily generalized to higher dimensions and other types of kinetic models.

%In the literature, adaptive semi-Lagrangian type wavelet method \cite{besse2003adaptive,gutnic_vlasov_2004,besse_wavelet-mra-based_2008}   and adaptive RKDG method \cite{zhuqiu} have been constructed. 

%The Vlasov-Poisson (VP) system. Besides the famous ``curse of dimensionality'', another major challenge for deterministic Vlasov simulations is the solution may develop filamentations in the phase space. It is therefore a good test bed for adaptive algorithms.  In the context of adaptive computation for Vlasov equations, we mention semi-Lagrangian type wavelet method \cite{besse2003adaptive,gutnic_vlasov_2004,besse_wavelet-mra-based_2008} and the $h$-adaptive RKDG method \cite{zhuqiu}. 

%\subsection{The collisionless Vlasov-Poisson system} 

\begin{exa}
\label{exa:vp}
We first consider the   non-dimensionalized single-species nonlinear VP system for plasma simulations in the zero-magnetic limit
\begin{eqnarray}
&&f_t + \bv\cdot\nabla_\bx f + \bE(t,\bx)\cdot\nabla_\bv f=0, \label{eq:V}\\
&&-\Delta_\bx\Phi(\bx) = \rho - 1,\quad \bE(\bx)=-\nabla_\bx\Phi \label{eq:poisson}
\end{eqnarray}
where $f(t,\bx,\bv)$ denotes the probability distribution function of electrons. $\bE(t,\bx)$ is the self-consistent electrostatic field given by Poisson's equation\eqref{eq:poisson} and $\rho(t,\bx)=\int_{\bv} f(t,\bx,\bv)  d\bv$ denotes the electron density. Ions are assumed to form a neutralizing background. 
\end{exa}

Periodic boundary condition is imposed in $x$-space. As a standard practice, the computational domain in $v$ is truncated to $  [-V_c,V_c]$, where $V_c$ is a constant chosen large enough  to impose zero boundary condition in the $v$-direction  $f_h(t,x,\pm V_c) = 0.$
The following set of initial conditions will be considered as classical benchmark numerical tests.
\begin{itemize}
	\item Landau damping:
	\begin{equation}
	f(0,x,v) = f_M(v)(1+A\cos(kx)),\quad x\in[0,L],\,v\in[-V_c,V_c],
	\end{equation}
	where $A=0.5$, $k=0.5$, $L=4\pi$, $V_c=2\pi$, and $f_M(v)=\frac{1}{\sqrt{2\pi}}e^{-v^2/2}$.
	\item Bump-on-tail instability:
	\begin{equation}
		f(0,x,v) = f_{BT}(v)(1+A\cos(kx)),\quad x\in[0,L],\,v\in[-V_c,V_c],
	\end{equation}
	where $A=0.04$, $k=0.3$, $L=20\pi/3$, $V_c=13$, and $$f_{BT}(v)=n_p\exp\left(-\frac{v^2}{2}\right)+n_b\exp\left(-\frac{|v-u|^2}{2v_t^2}\right),$$ 
	where $n_p = \frac{9}{10\sqrt{10\pi}},\,n_b = \frac{2}{10\sqrt{10\pi}},\,u=4.5,\,v_t=0.5.$
	\item Two-stream instability I:
	\begin{equation}
	f(0,x,v) = f_{TSI}(v)(1+A\cos(kx)),\quad x\in[0,L],\,v\in[-V_c,V_c],
	\end{equation}
	where $A=0.05$, $k=0.5$, $L=4\pi$, $V_c=2\pi$, and $f_{TSII}(v)=\frac{1}{\sqrt{2\pi}}v^2e^{-v^2/2}$.
	\item Two-stream instability II:
		\begin{equation}
		f(0,x,v) = f_{TSII}(v)(1+A\cos(kx)),\quad x\in[0,L],\,v\in[-V_c,V_c],
		\end{equation}
		where $A=0.05$, $k=2/13$, $L=13\pi$, $V_c=5$, and $$f_{TSII}(v)=\frac{1}{2v_t\sqrt{2\pi}}\left(\exp\left(-\frac{|u+v|^2}{2v_t^2}\right)+\exp\left(-\frac{|u-v|^2}{2v_t^2}\right)\right),$$ 
		where $u=0.99,\,v_t=0.3.$
\end{itemize}

In the literature,  RKDG schemes for the VP system  \cite{Ayuso2009, heath2012discontinuous,  cheng_vp} have been extensively studied. They are shown to have superior performance in conservation. Our previous work on sparse grid DG method \cite{guo_sparsedg} focused on the closely related Vlasov-Amp\`{e}re  (VA) system. The solver in \cite{guo_sparsedg} successfully reduced the DOFs of the equations while maintaining key conservation properties. However, when $t$ gets large and filamentation becomes severe, the sparse grid method has difficulties resolving the fine structures in the phase space. It is therefore to our interest to investigate if the adaptive multiresolution scheme can achieve a good balance between computational cost and numerical resolution.

\begin{figure}[htp]
	\begin{center}
		\subfigure[]{\includegraphics[width=.42\textwidth]{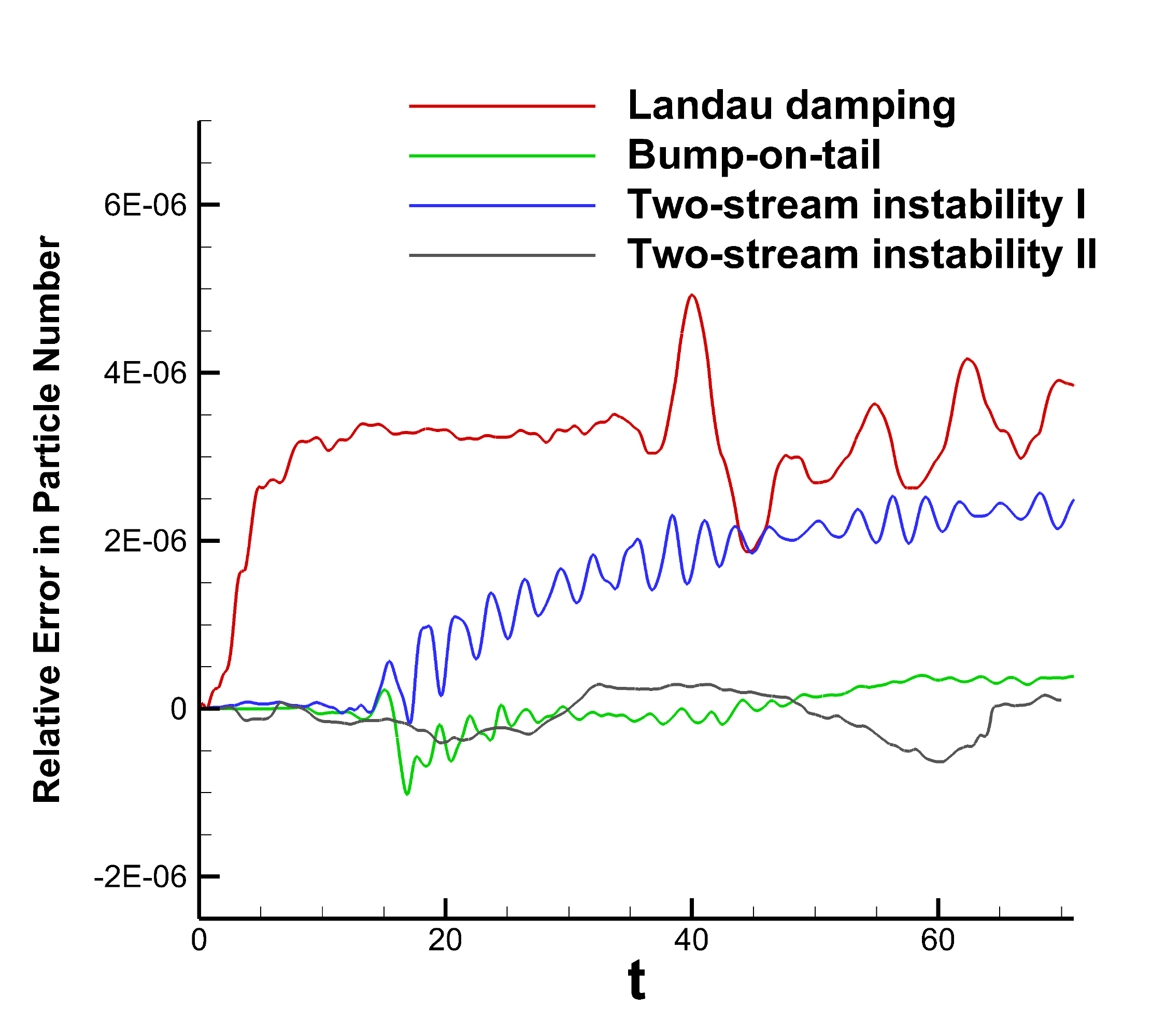}}
		\subfigure[]{\includegraphics[width=.42\textwidth]{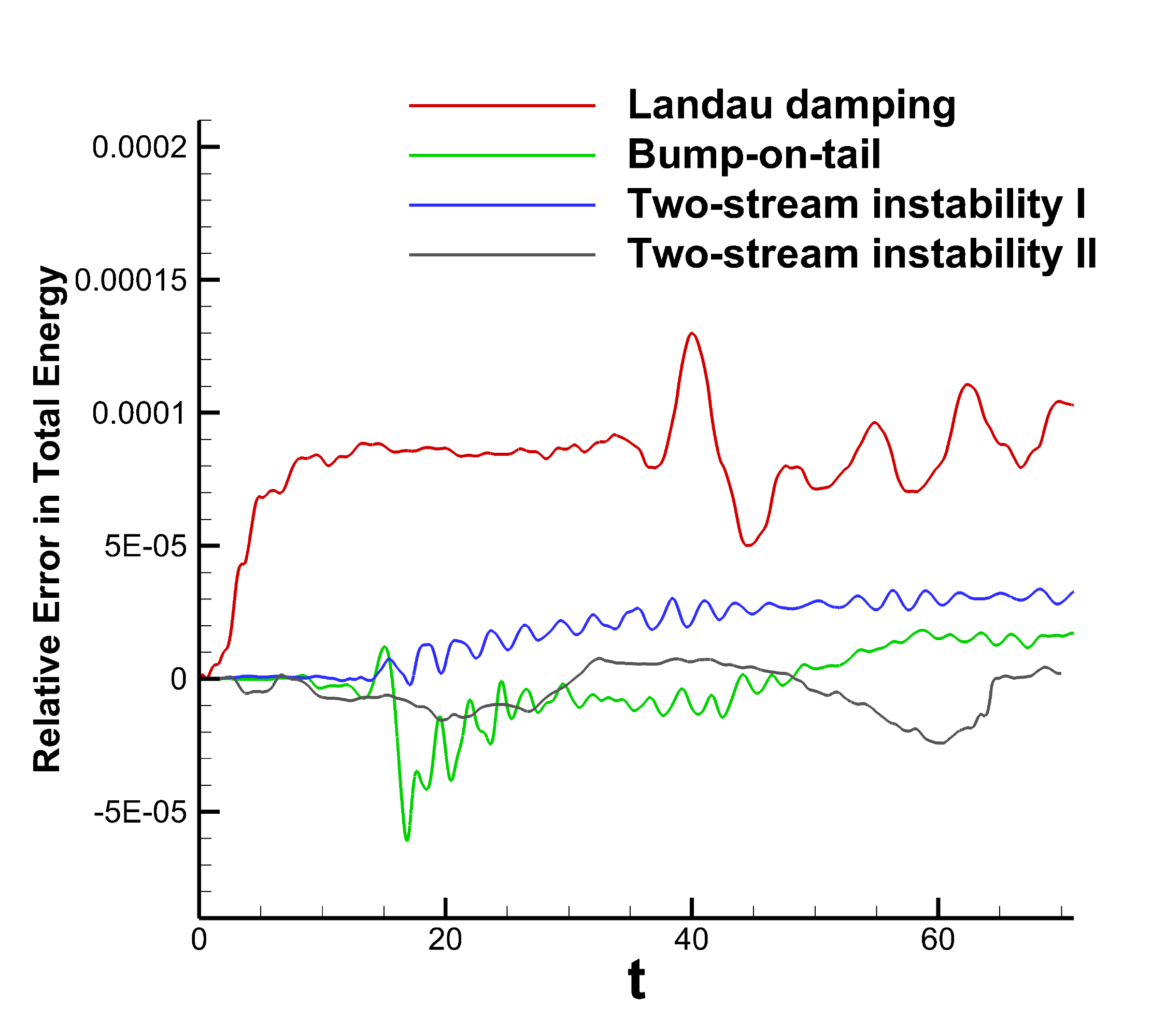}}\\
		\subfigure[]{\includegraphics[width=.42\textwidth]{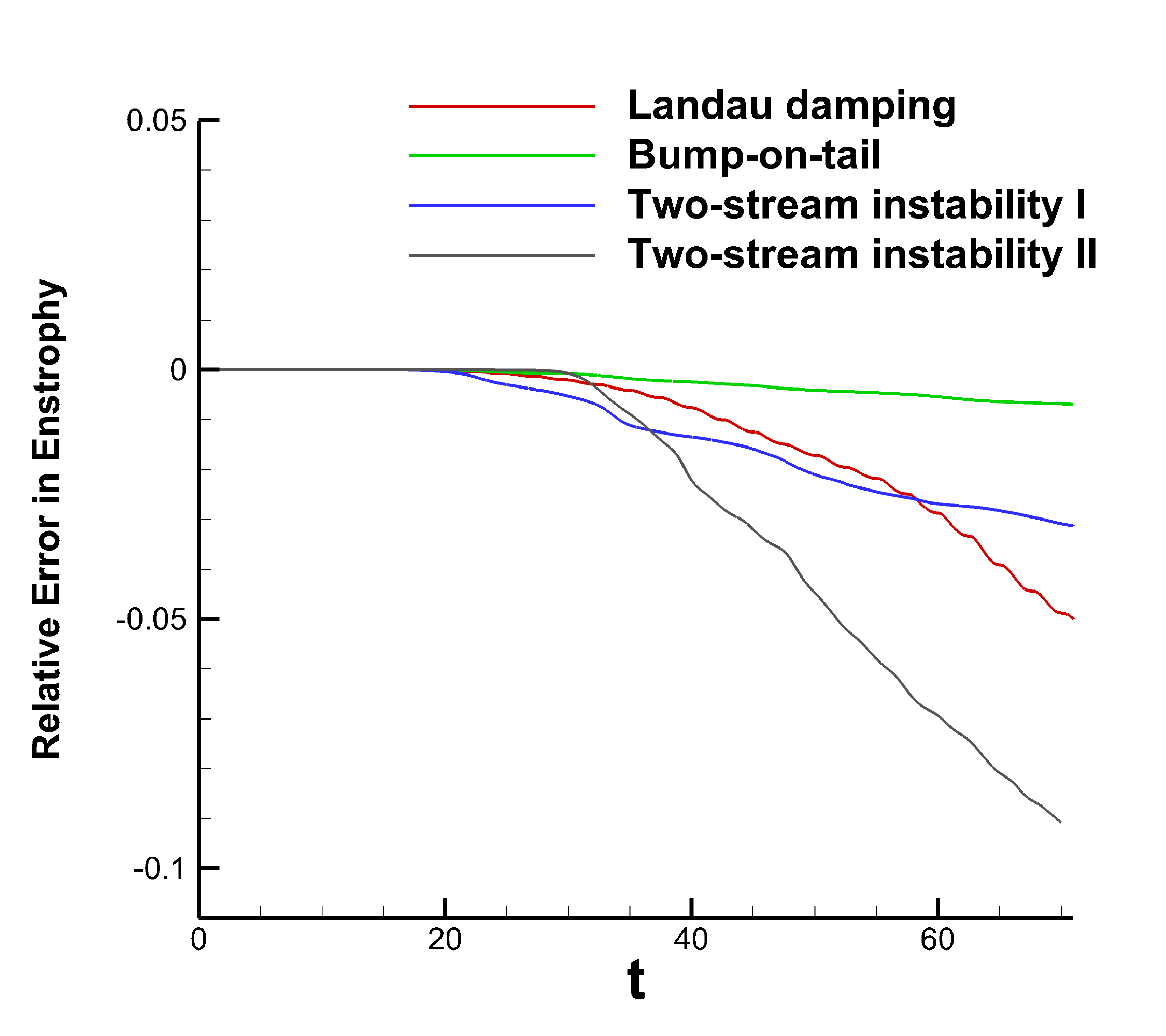}}
		\subfigure[]{\includegraphics[width=.42\textwidth]{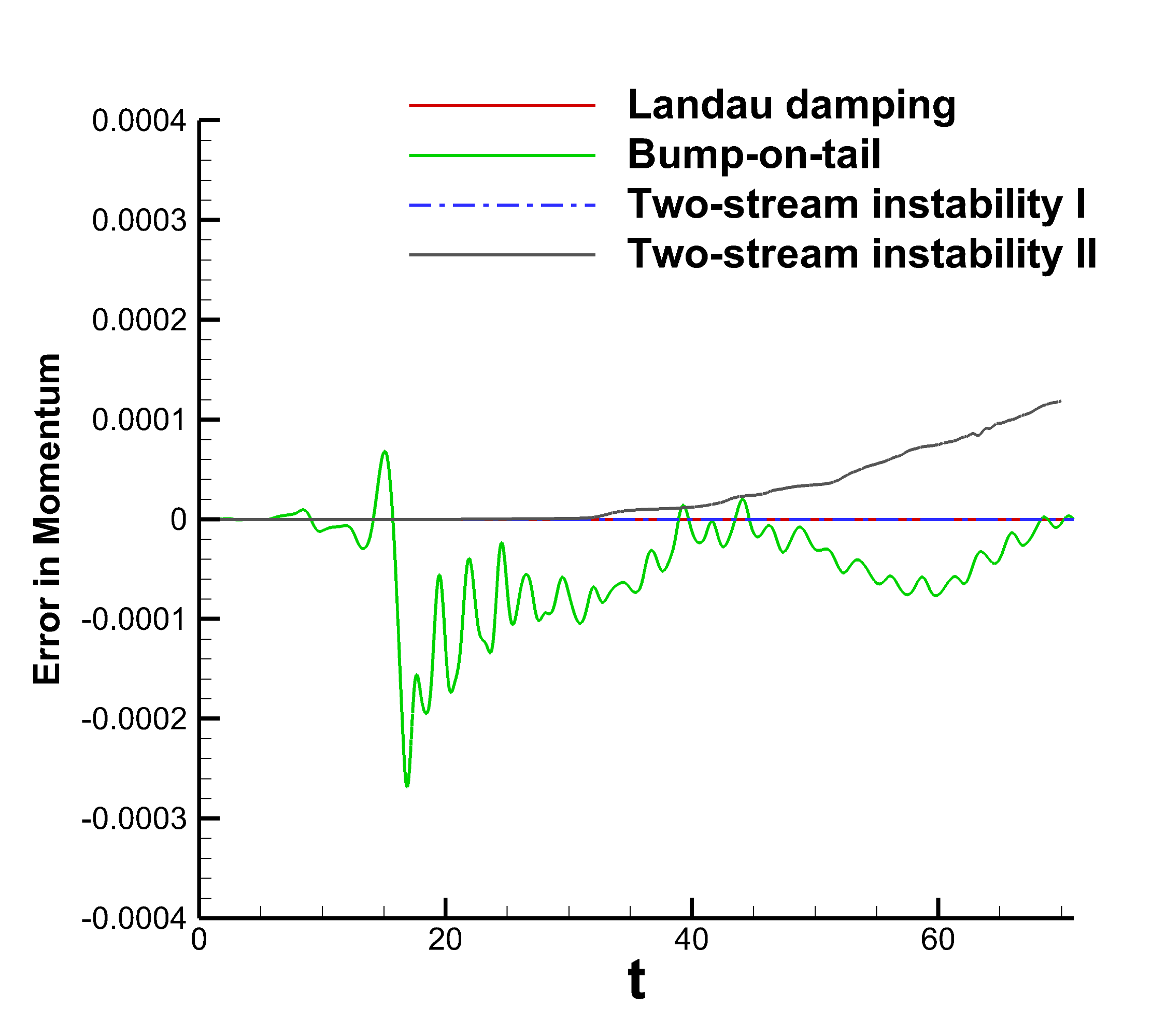}}\\
	\end{center}
	\caption{Example \ref{exa:vp}. Evolution of the relative errors in total particle number (a), total energy (b), enstrophy (c), and evolution of error in momentum (d) with the indicated initial conditions.  $N=7$. $k=3$. $\varepsilon=10^{-5}$.}
	\label{fig:evo_vp}
\end{figure}

We apply the adaptive algorithm to the Vlasov equation as outlined in Section \ref{sec:method}. The Poisson equation is solved by a standard local DG method \cite{Arnold_2002_SIAM_DG} on the finest level mesh in the $x$-direction.
 In the simulations, we use $N=7$, $\varepsilon=10^{-5}$ and $k=3$. 
 First we investigate the conservative properties of the   scheme.  
 The VP system is known to preserve many physical invariants, including 
 the particle number  $\int_\bx\int_\bv f(t,\bx,\bv)\,d\bx d\bv,$ momentum $\int_\bx\int_\bv \bv f(t,\bx,\bv)\,d\bx d\bv,$ enstrophy $\int_\bx\int_\bv |f(t,\bx,\bv)|^2\,d\bx d\bv,$ and total energy $\frac12\int_\bx\int_\bv f(t,\bx,\bv)|\bv|^2\,d\bx d\bv +\frac12 \int_\bx |\bE(t,\bx)|^2\,d\bx.$ Generally speaking, it is difficult for a numerical method to preserve all those invariants.
 By careful design, DG methods have been designed to preserve the particle number and the energy of the system \cite{Ayuso2009,cheng2014energy}. %Here, we will investigate the time evolution of those numerical quantities to benchmark the adaptive schemes.
 For our scheme, in Figure \ref{fig:evo_vp}, we report the time evolution of the relative errors in total particle number, total energy, enstrophy, and evolution of error in momentum. It is observed that the total particle number is conserved up to the magnitude of $\varepsilon$. This is not as well conserved as  a traditional RKDG method. However, it is expected because the adaptive algorithm only keeps  elements   above the error threshold in the hash table and causing   the truncation errors at the velocity boundary to be on the same magnitude of $\varepsilon$,   contributing to the numerical errors in particle numbers. However, we do comment  that  the addition and removal of elements other than level $\bl=\mathbf{0}$ in the refinement and coarsening steps will not change the numerical mass because  the basis functions are orthogonal. Similarly, the total energy and momentum also show visible and slightly larger errors than the standard RKDG method. The enstrophy exhibits the most visible decay because of the choice of upwind numerical flux. %More detailed discussions of numerical conservation of RKDG method can be found in \cite{Ayuso2009,cheng_vp, cheng2014energy}.
 
 %Also note that even though the total energy is not  conserved, the relative error is still on a  scale of $\varepsilon$, again because of the truncation error at the boundary. Furthermore, he momentum is well conserved and  the visible decay in enstrophy is expected because of  the dissipative upwind flux.

%We further present some numerical data to benchmark the proposed adaptive scheme. We first consider the Log Fourier modes of the electric field $E(t,x)$ as  functions of time, which are defined as 
%$$\log FM_n(t)=\log_{10}\left(\frac1L\sqrt{\left| \int_0^L E(t,x)\,\sin(knx)\,dx\right|^2+\left| \int_0^L E(t,x)\,\cos(knx)\,dx\right|^2}\right).$$ 
%In Figure \ref{fig:fou_lan}, we show the time evolution of the first four Log Fourier modes when simulating Landau damping. Again, the configuration of $N=7$, $k=3$, and $\varepsilon=10^{-5}$ is used in our computation. We only report the results for Landau damping for brevity, and they agree with other calculations in the literature.
%
%\begin{figure}[htp]
%	\begin{center}
%		\subfigure[]{\includegraphics[width=.42\textwidth]{figures/fm1-eps-converted-to}}
%		\subfigure[]{\includegraphics[width=.42\textwidth]{figures/fm2-eps-converted-to}}\\
%		\subfigure[]{\includegraphics[width=.42\textwidth]{figures/fm3-eps-converted-to}}
%		\subfigure[]{\includegraphics[width=.42\textwidth]{figures/fm4-eps-converted-to}}\\
%	\end{center}
%	\caption{The first four log Fourier modes of Landau damping. $k=3$. $N=7$. $\varepsilon=10^{-5}$.}
%	\label{fig:fou_lan}
%\end{figure}

\begin{figure}[htp]
	\begin{center}
		\subfigure[]{\includegraphics[width=.42\textwidth]{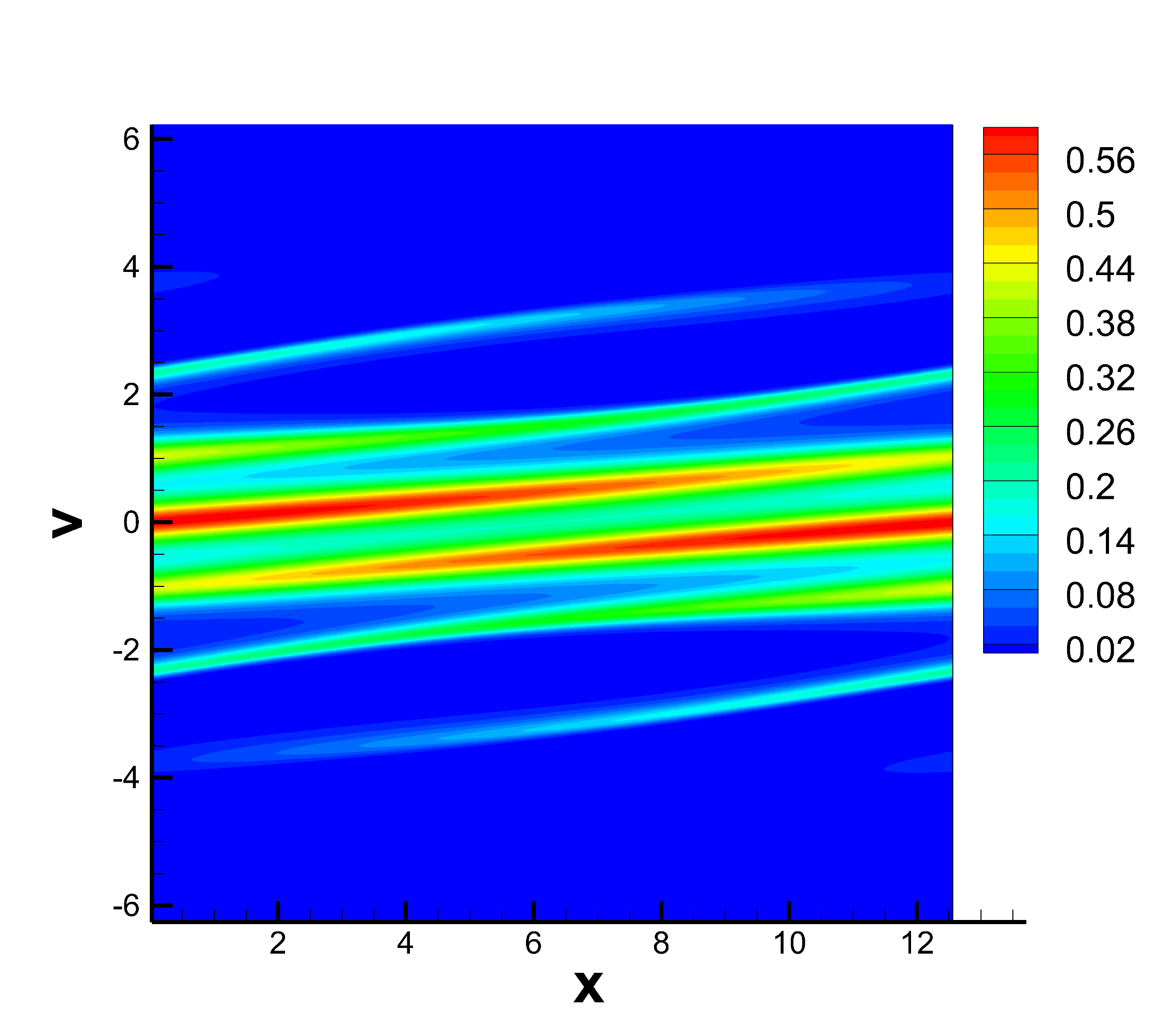}}
		\subfigure[]{\includegraphics[width=.42\textwidth]{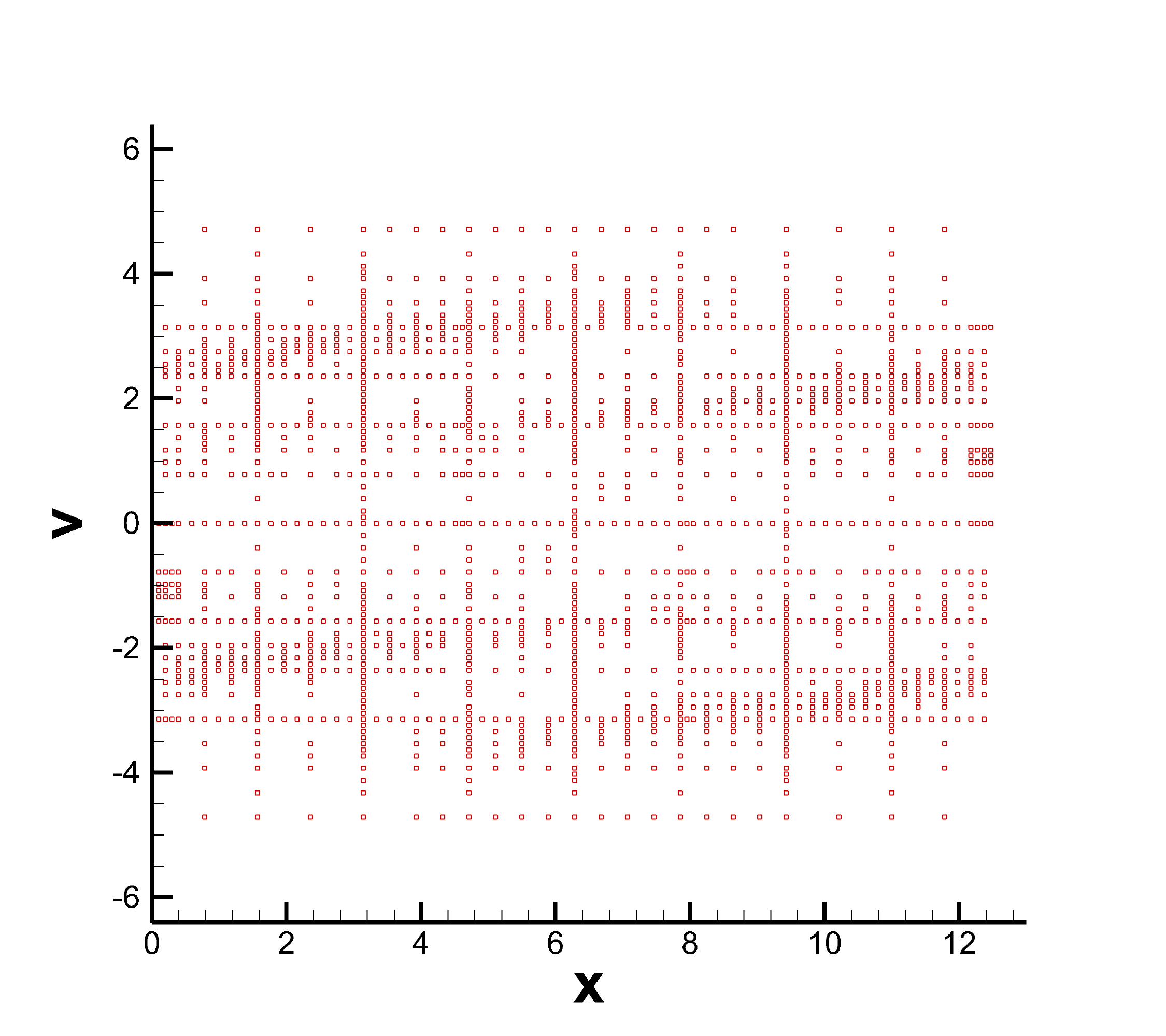}}\\
		\subfigure[]{\includegraphics[width=.42\textwidth]{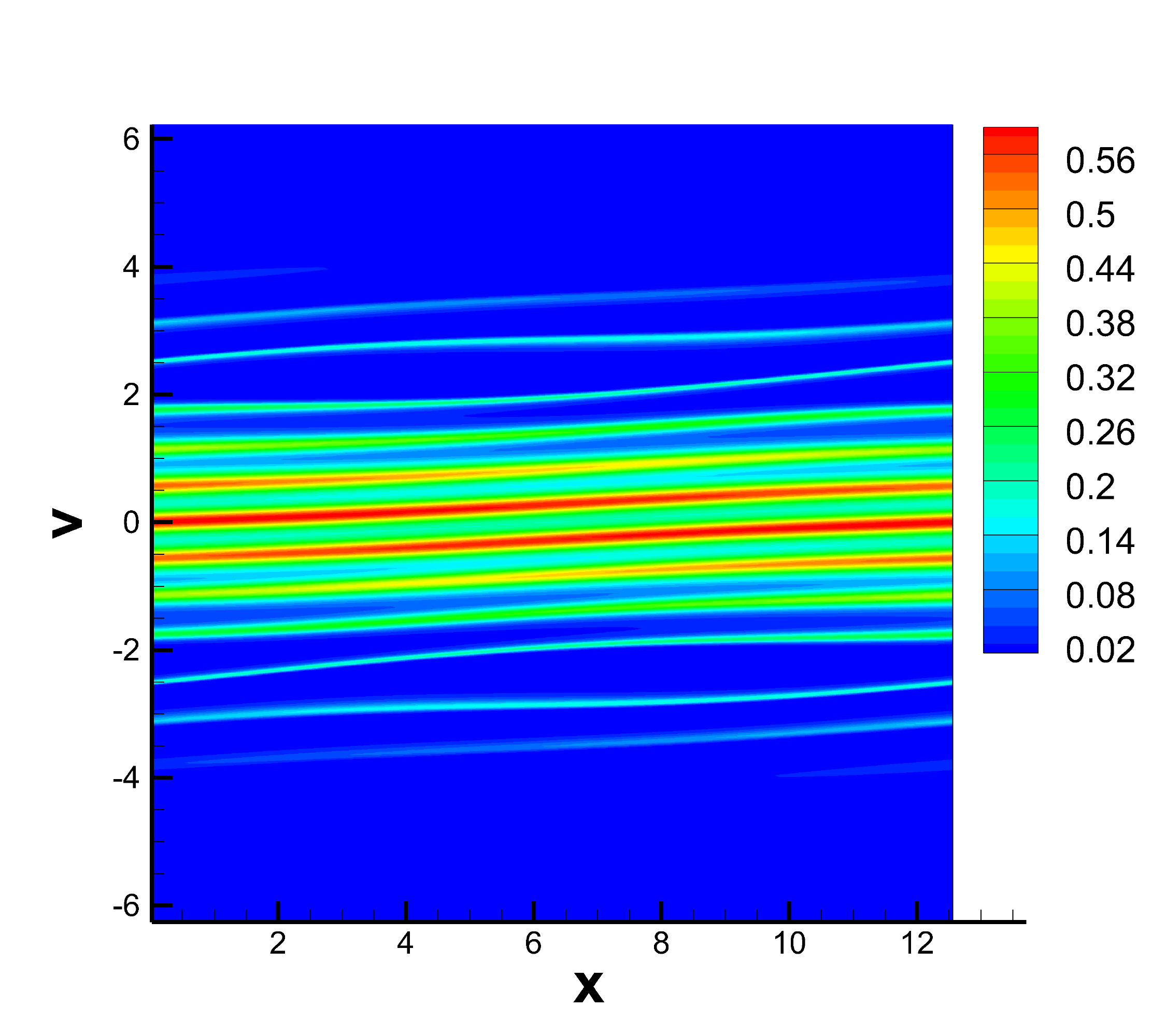}}
		\subfigure[]{\includegraphics[width=.42\textwidth]{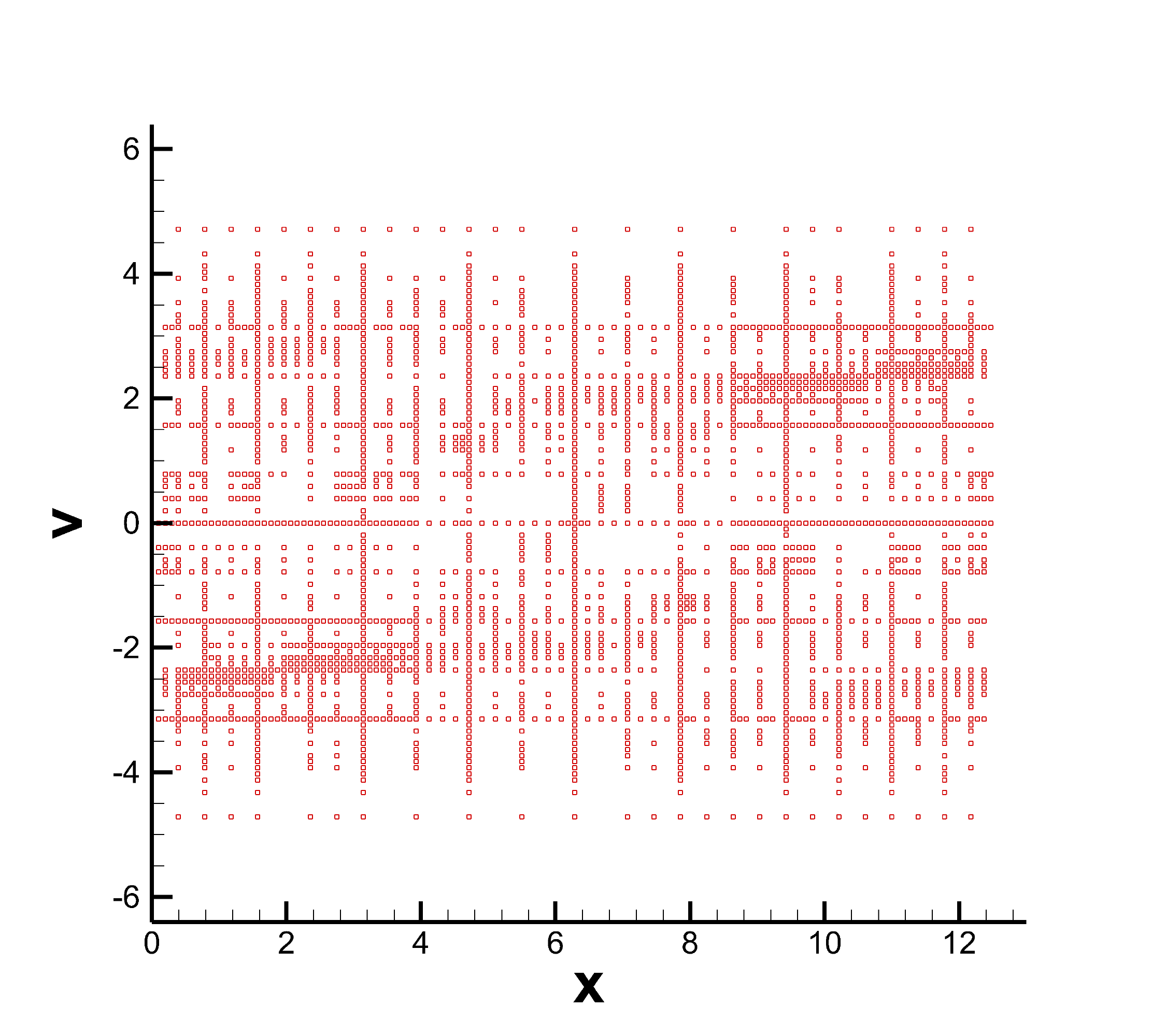}}\\
				\subfigure[]{\includegraphics[width=.42\textwidth]{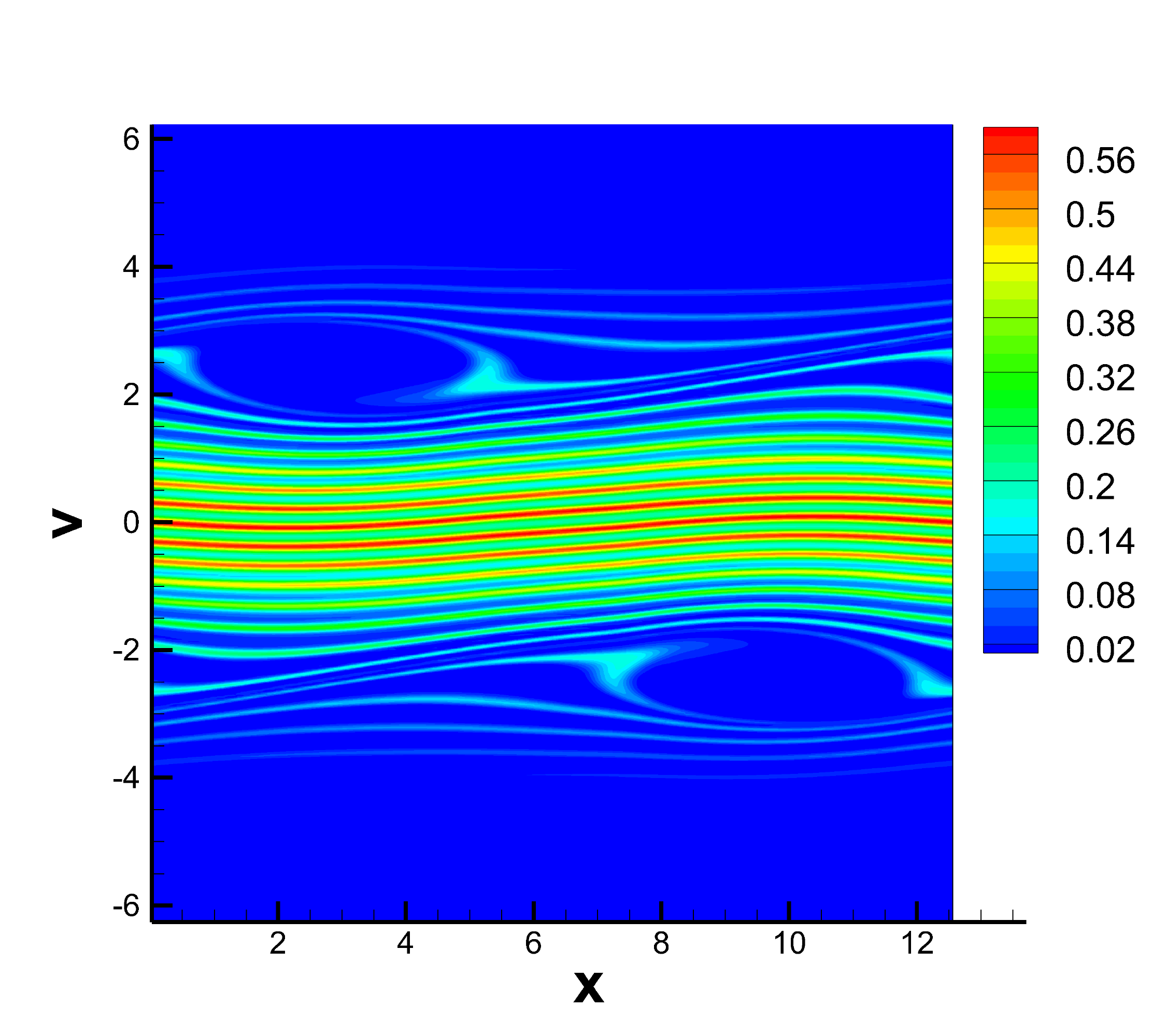}}
				\subfigure[]{\includegraphics[width=.42\textwidth]{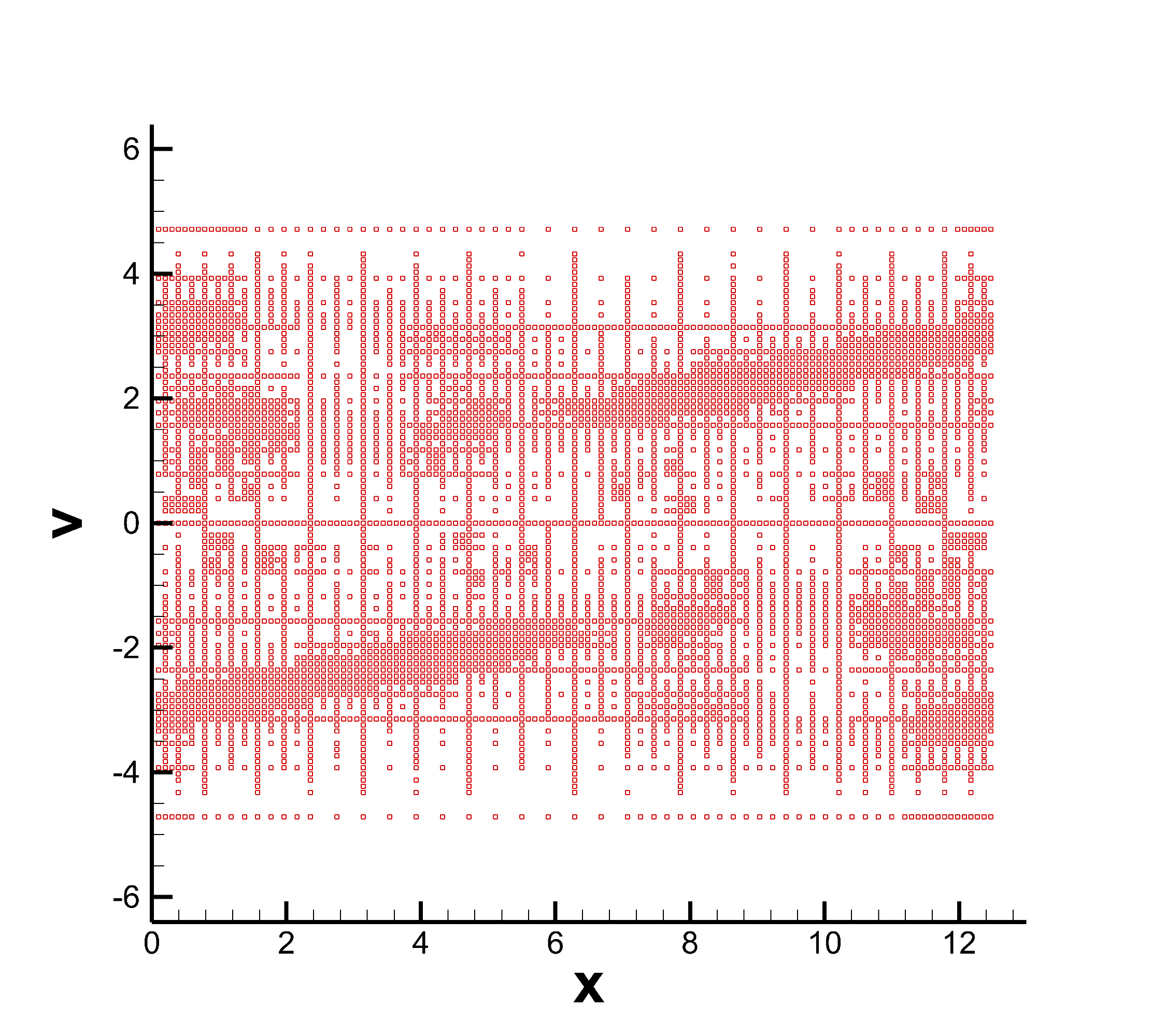}}\\
		%\subfigure[]{\includegraphics[width=.42\textwidth]{figures/con_lan_t60_p3}}
		%\subfigure[]{\includegraphics[width=.42\textwidth]{figures/mesh_lan_t60_p3}
		%\\
	\end{center}
	\caption{Example \ref{exa:vp}. Landau damping. Phase space contour plots and the associated active elements at   $t=10$ (a-b),  $t=20$ (c-d), $t=40$ (e-f). $N=7$. $k=3$.  $\varepsilon=10^{-5}$.}
	\label{fig:con_lan}
\end{figure}

\begin{figure}[htp]
	\begin{center}
		\subfigure[]{\includegraphics[width=.42\textwidth]{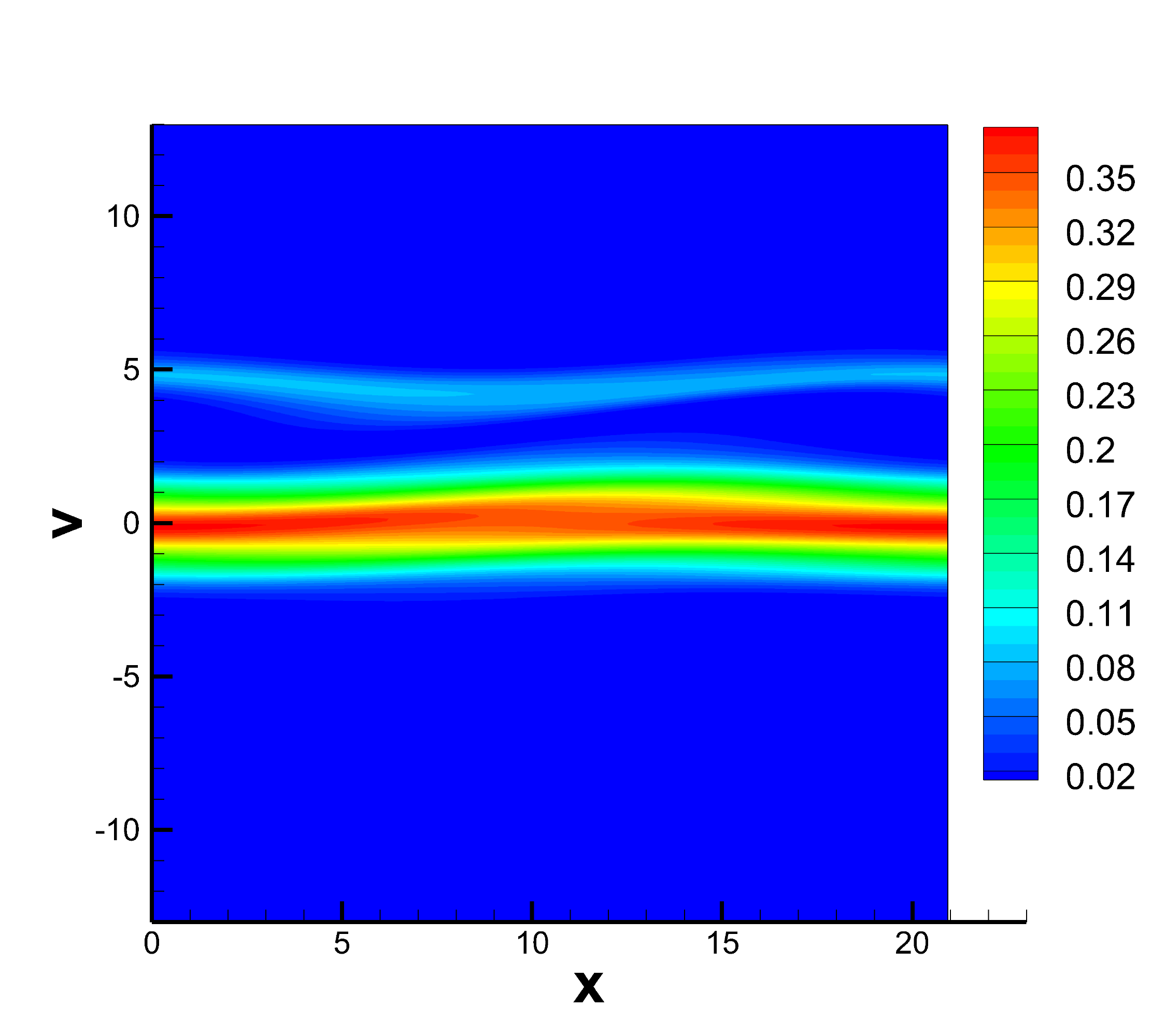}}
		\subfigure[]{\includegraphics[width=.42\textwidth]{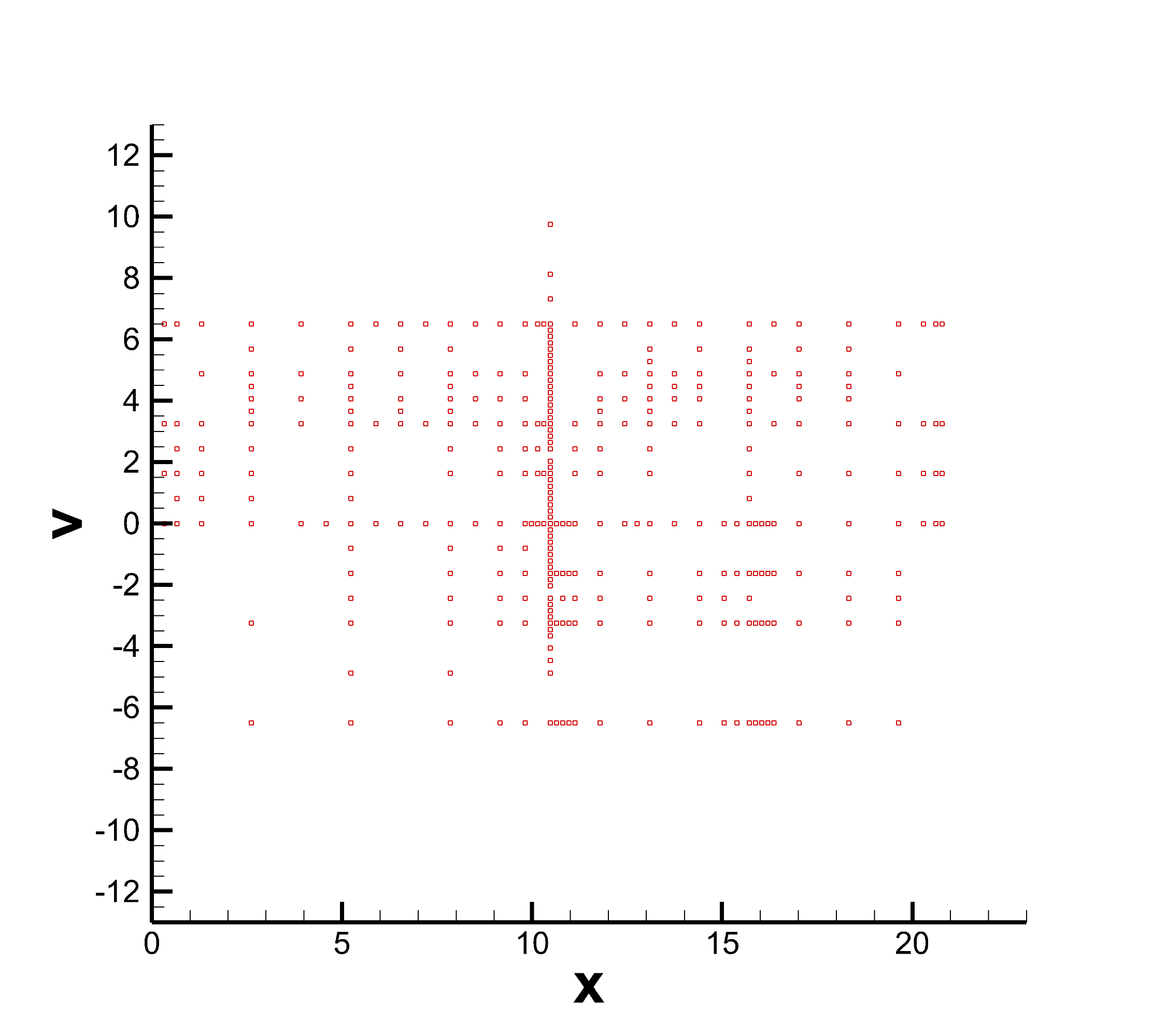}}\\
		\subfigure[]{\includegraphics[width=.42\textwidth]{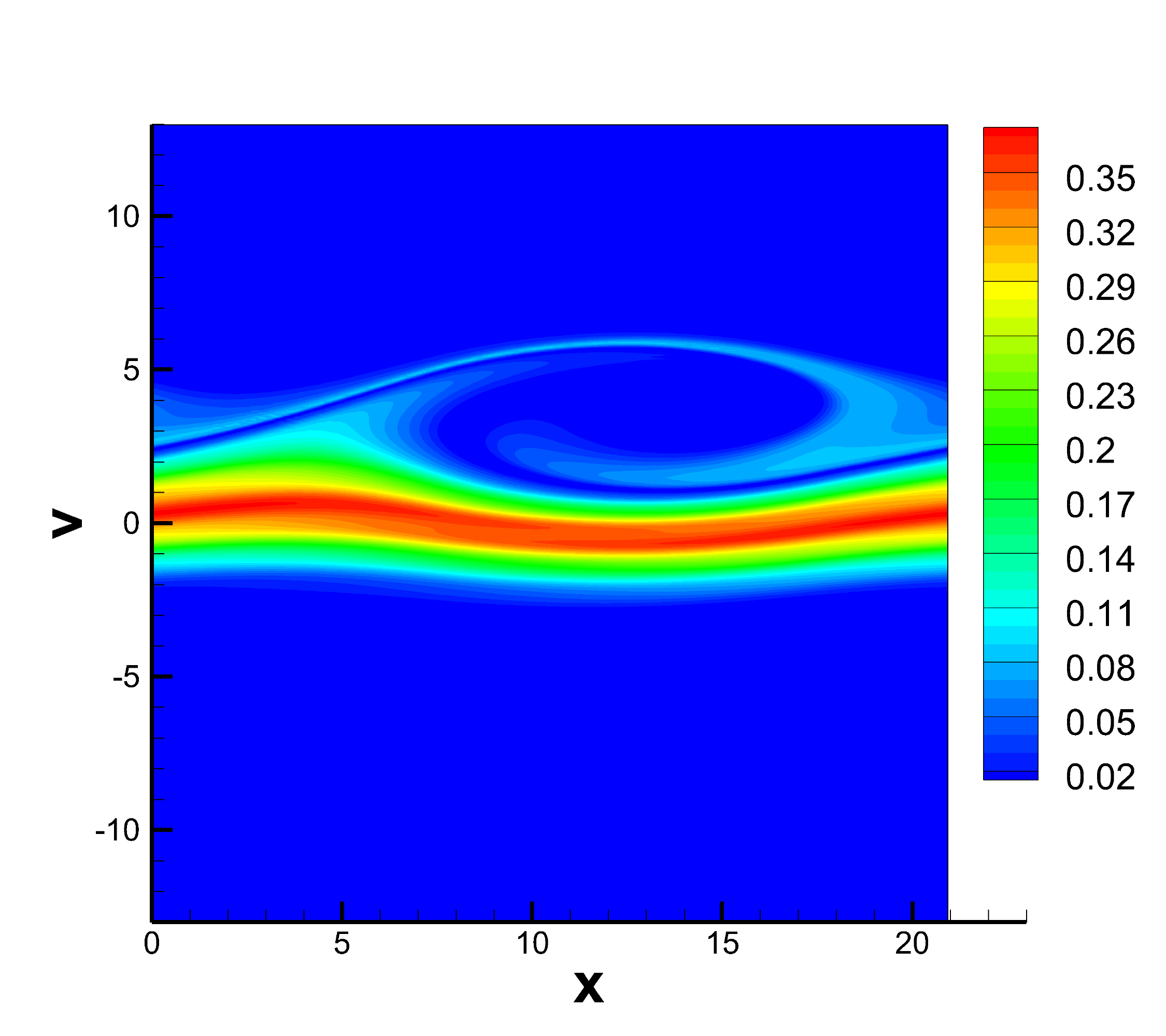}}
		\subfigure[]{\includegraphics[width=.42\textwidth]{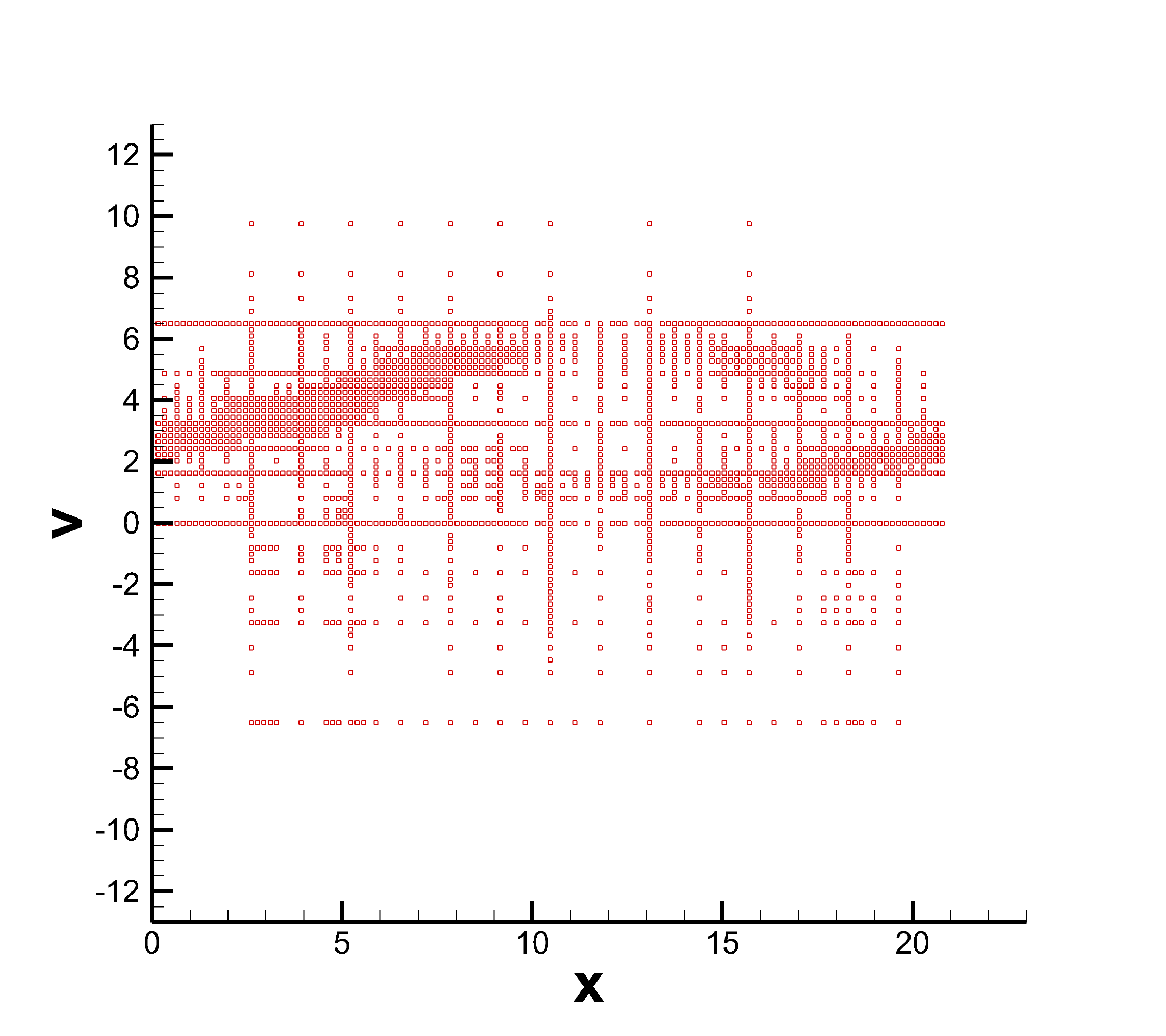}}\\
				\subfigure[]{\includegraphics[width=.42\textwidth]{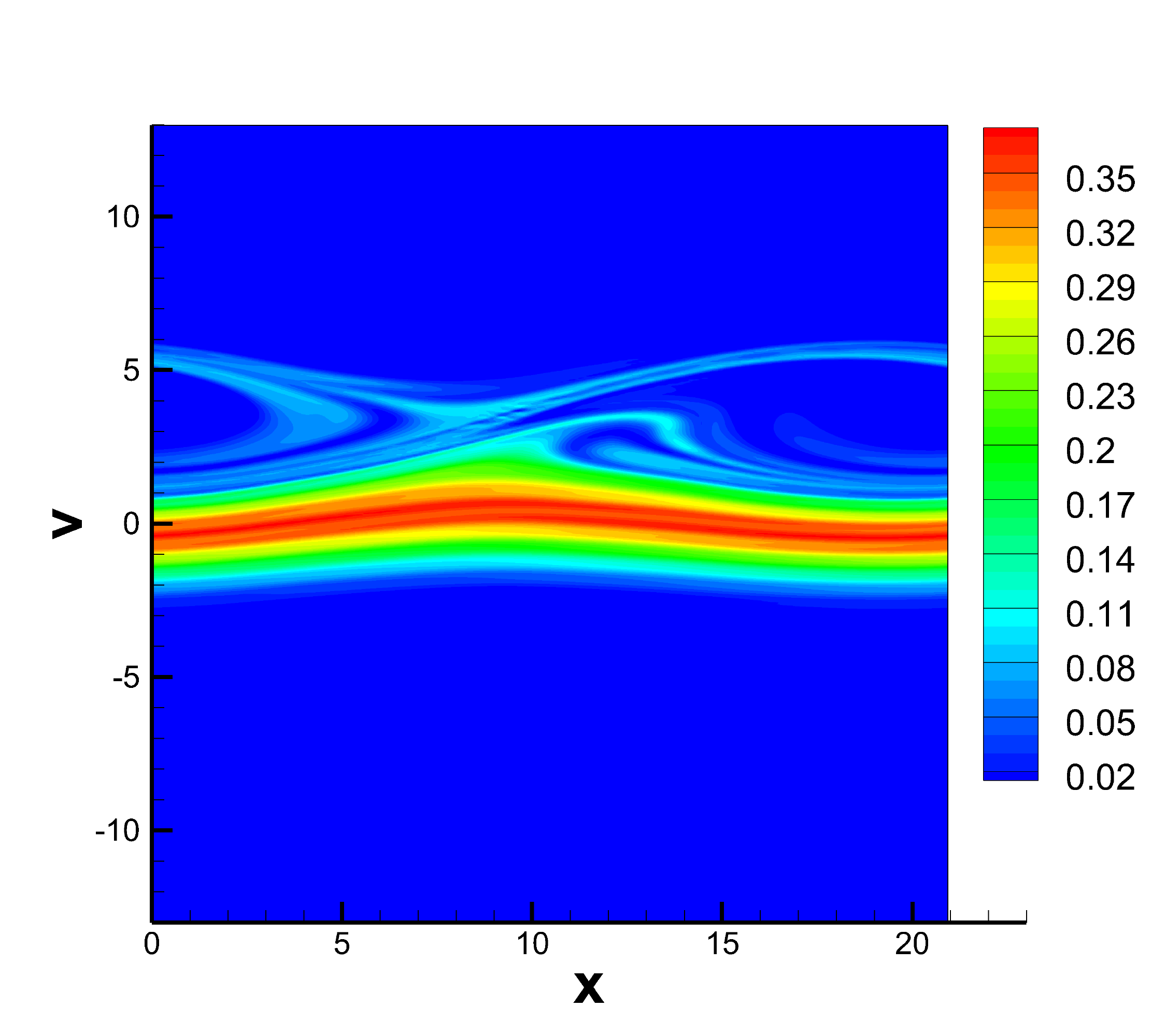}}
				\subfigure[]{\includegraphics[width=.42\textwidth]{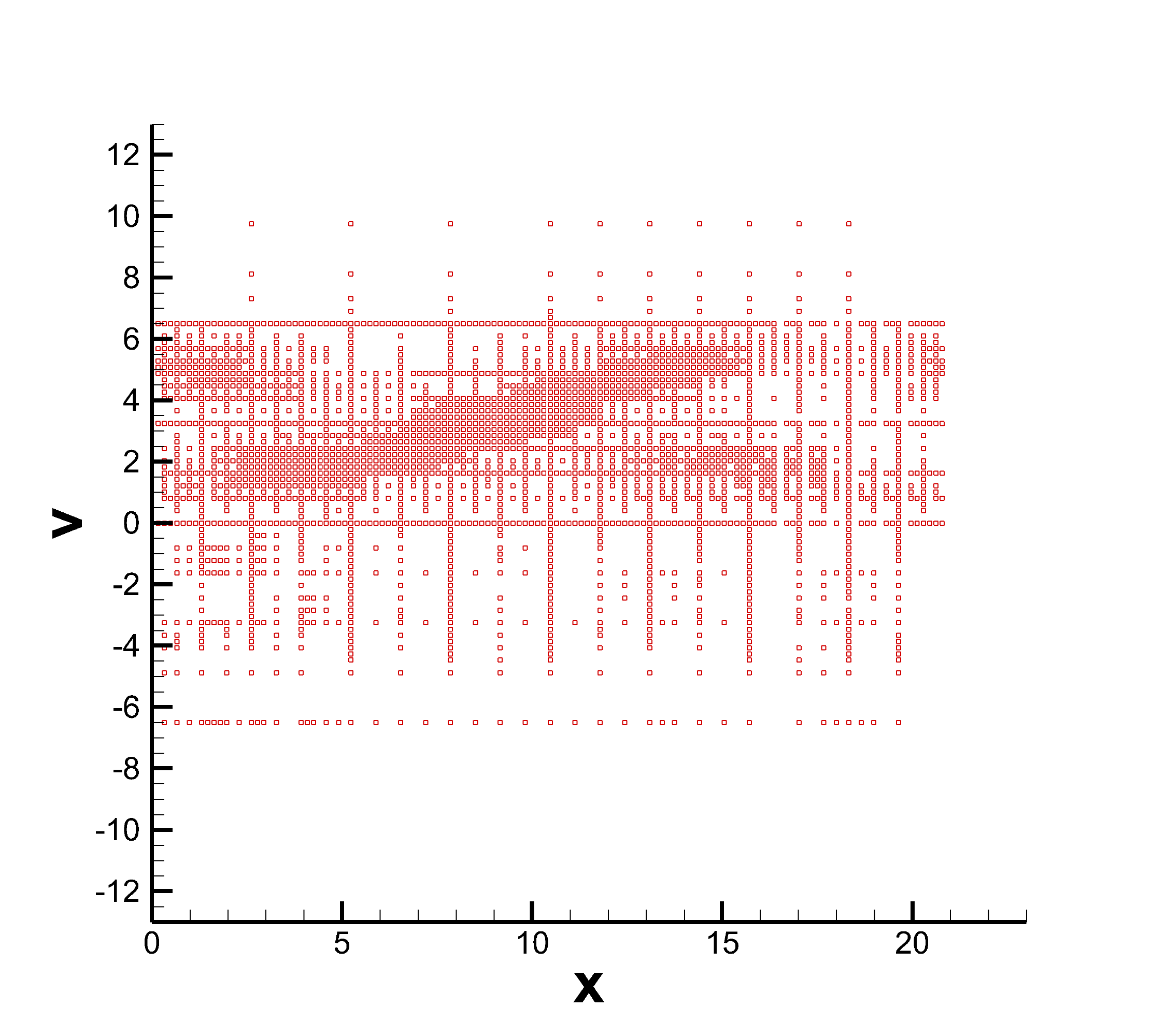}}\\
%		\subfigure[]{\includegraphics[width=.42\textwidth]{figures/con_bump_t60_p3}}
%		\subfigure[]{\includegraphics[width=.42\textwidth]{figures/mesh_bump_t60_p3}}\\
	\end{center}
	\caption{Example \ref{exa:vp}. Bump-on-tail instability. Phase space contour plots and the associated active elements at   $t=10$ (a-b),  $t=20$ (c-d), $t=40$ (e-f). $N=7$. $k=3$.  $\varepsilon=10^{-5}$.}
	\label{fig:con_bump}
\end{figure}

\begin{figure}[htp]
	\begin{center}
		\subfigure[]{\includegraphics[width=.42\textwidth]{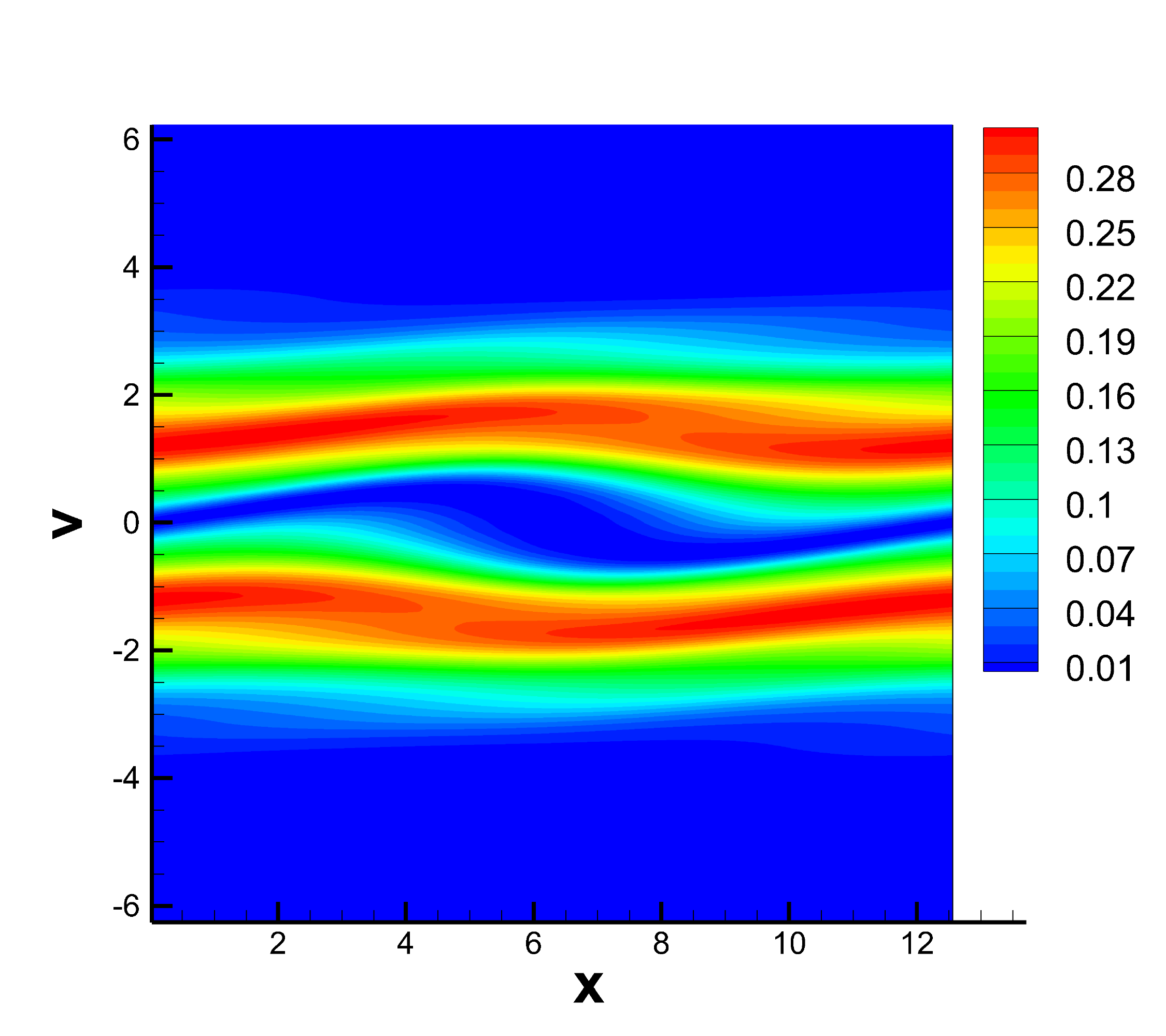}}
		\subfigure[]{\includegraphics[width=.42\textwidth]{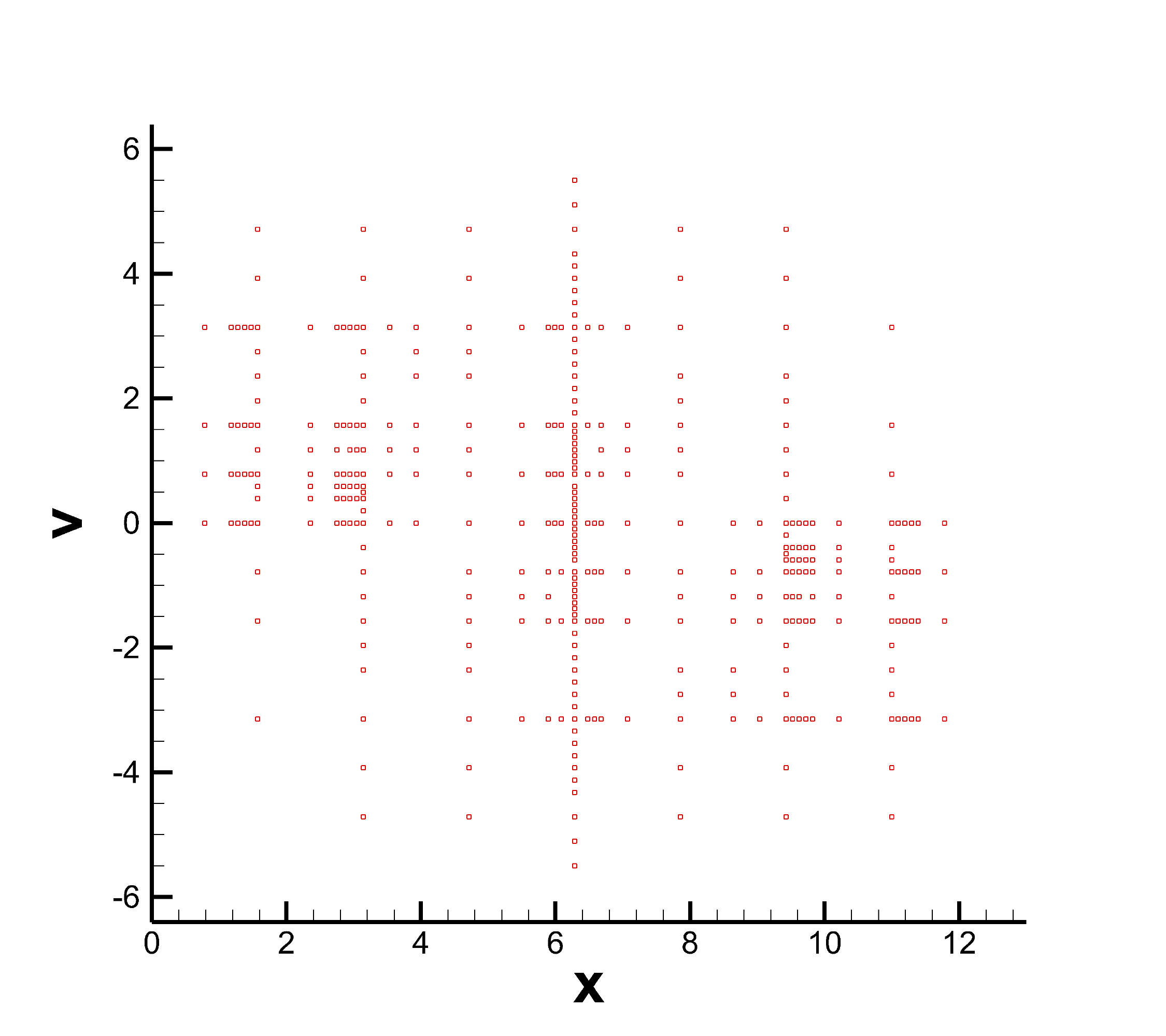}}\\
		\subfigure[]{\includegraphics[width=.42\textwidth]{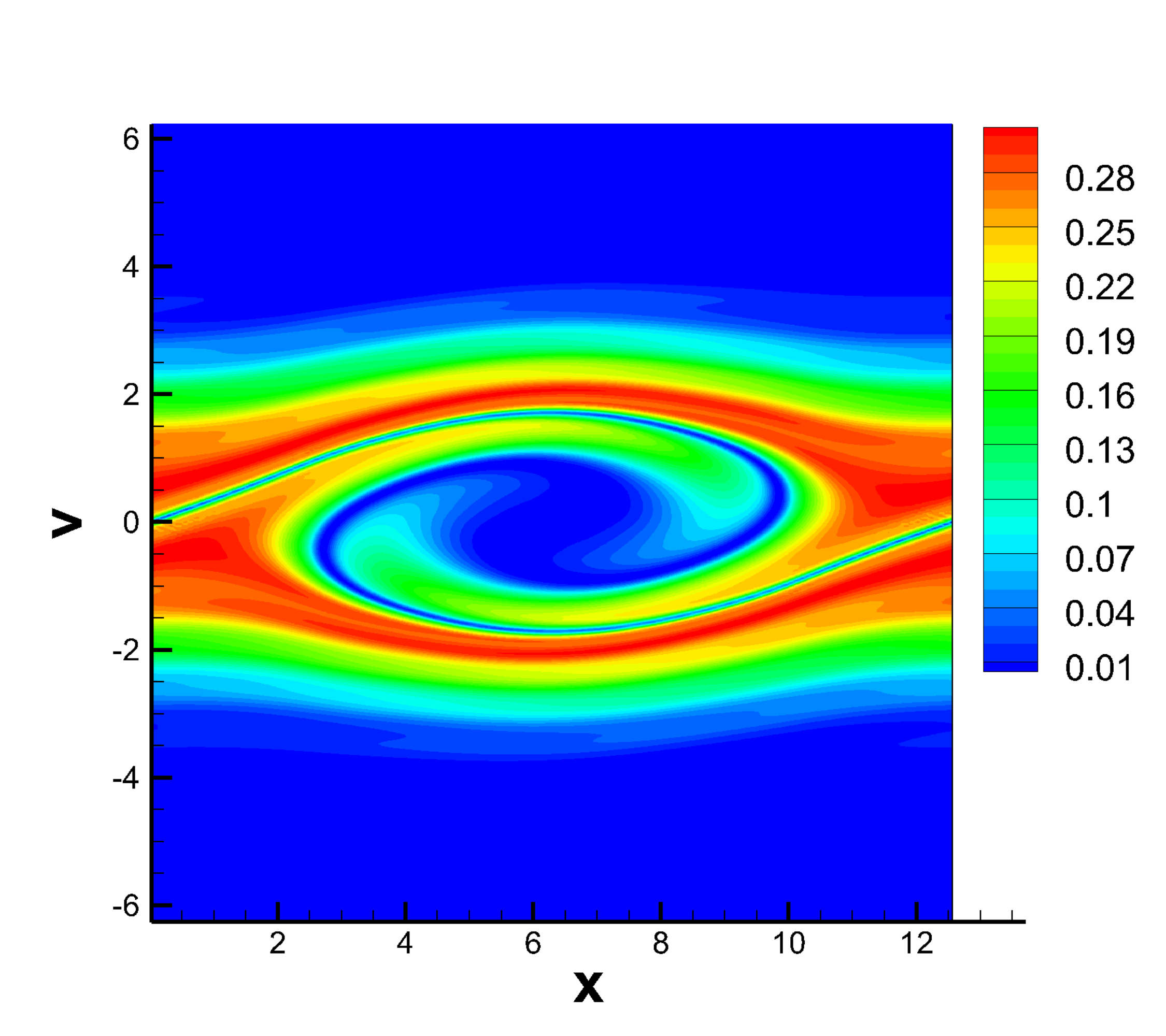}}
		\subfigure[]{\includegraphics[width=.42\textwidth]{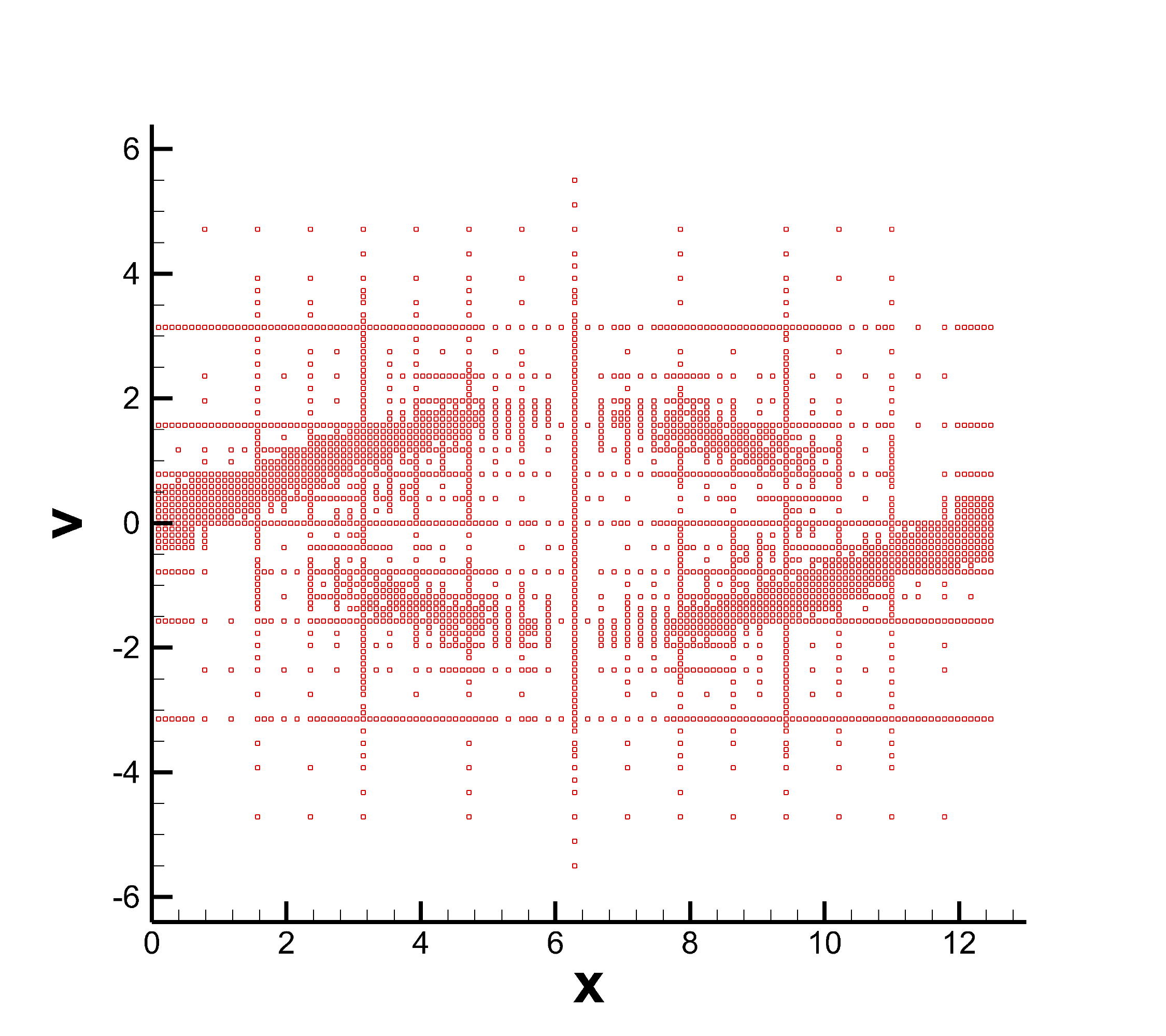}}\\
				\subfigure[]{\includegraphics[width=.42\textwidth]{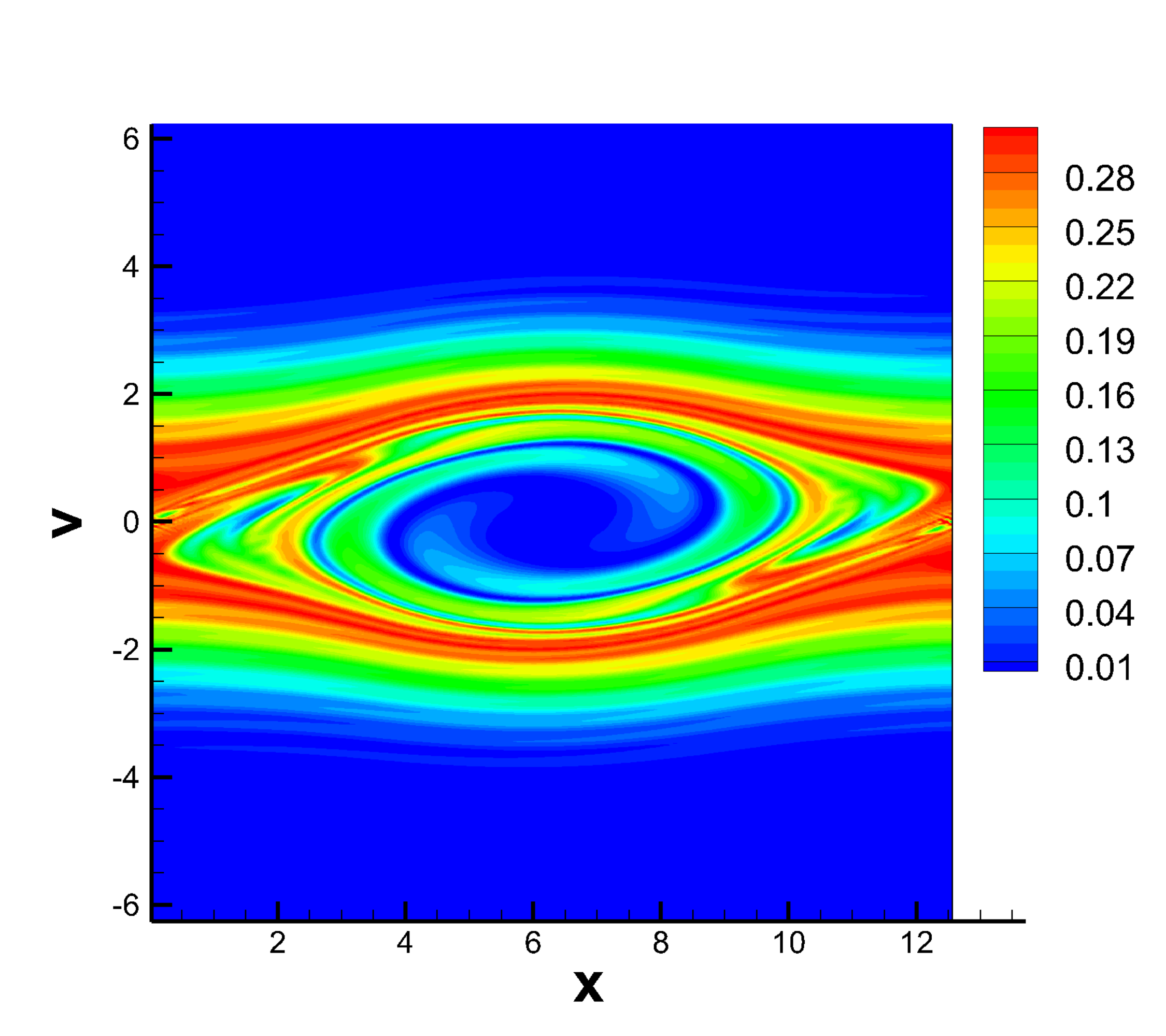}}
				\subfigure[]{\includegraphics[width=.42\textwidth]{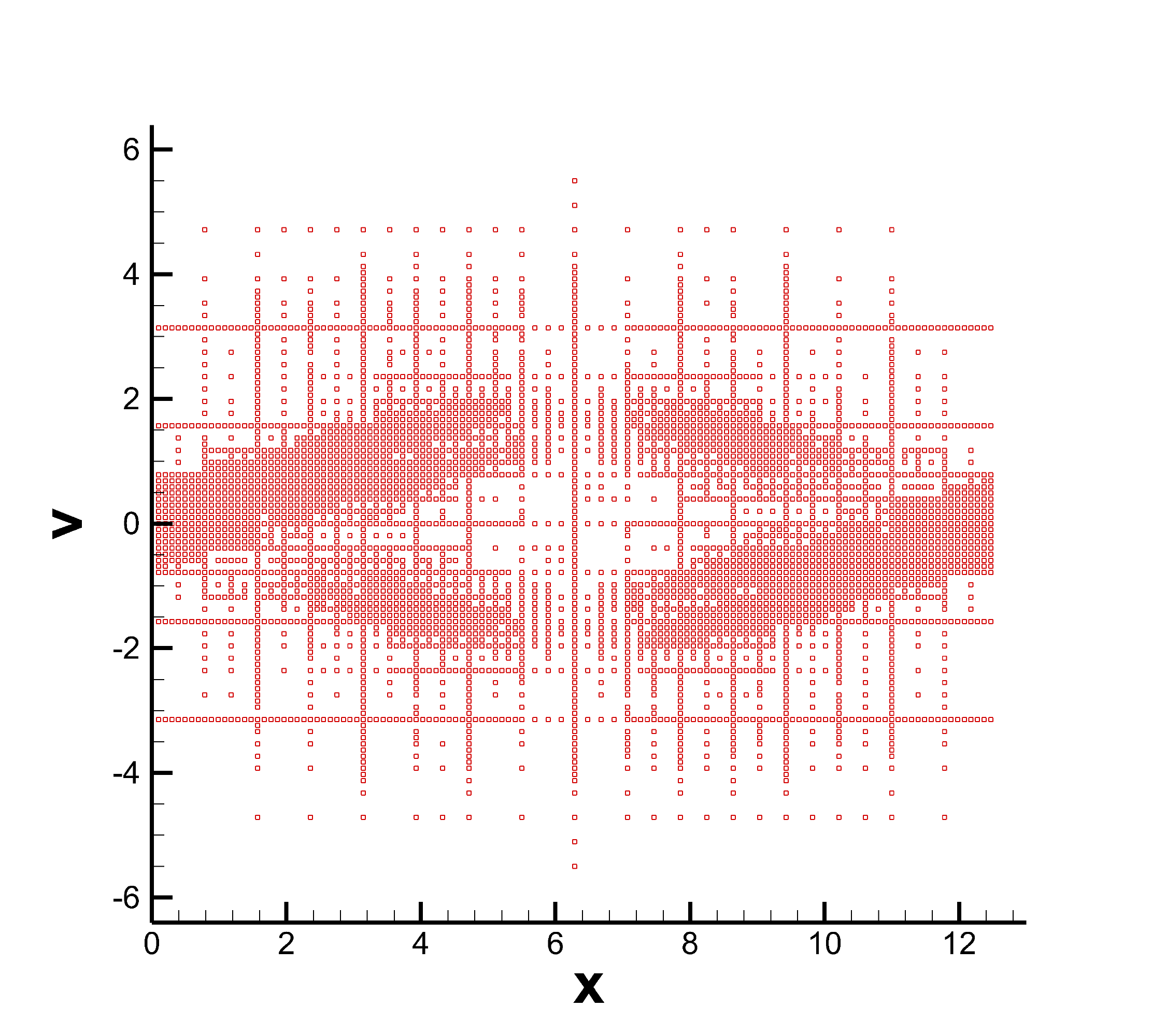}}\\
%		\subfigure[]{\includegraphics[width=.42\textwidth]{figures/con_two_t60_p3}}
%		\subfigure[]{\includegraphics[width=.42\textwidth]{figures/mesh_two_t60_p3}}\\
	\end{center}
	\caption{Example \ref{exa:vp}. Two-stream instability I. Phase space contour plots and the associated active elements at   $t=10$ (a-b),  $t=20$ (c-d), $t=40$ (e-f). $N=7$. $k=3$.  $\varepsilon=10^{-5}$.}
	\label{fig:con_two}
\end{figure}

\begin{figure}[htp]
	\begin{center}
		\subfigure[]{\includegraphics[width=.42\textwidth]{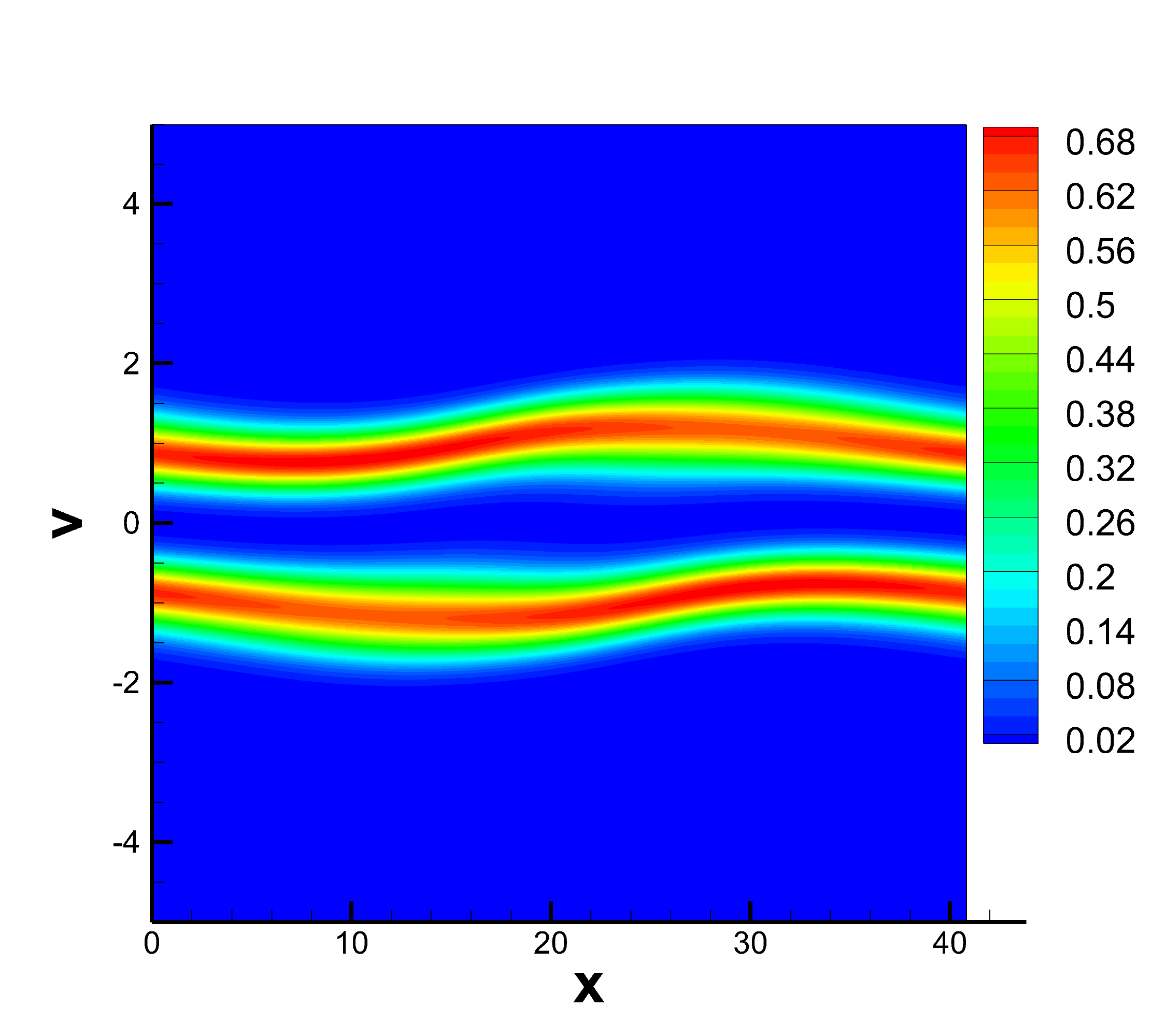}}
		\subfigure[]{\includegraphics[width=.42\textwidth]{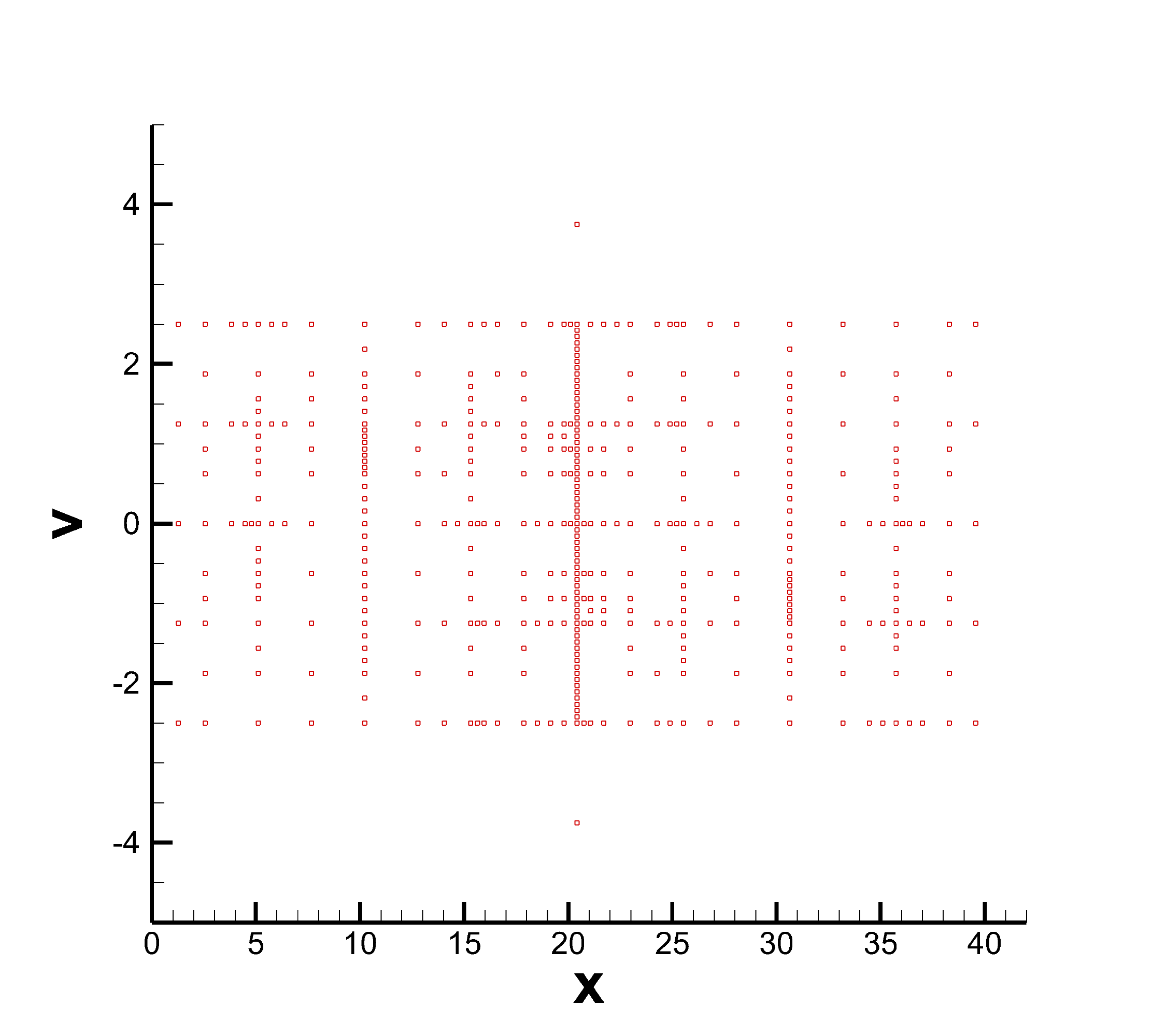}}\\
				\subfigure[]{\includegraphics[width=.42\textwidth]{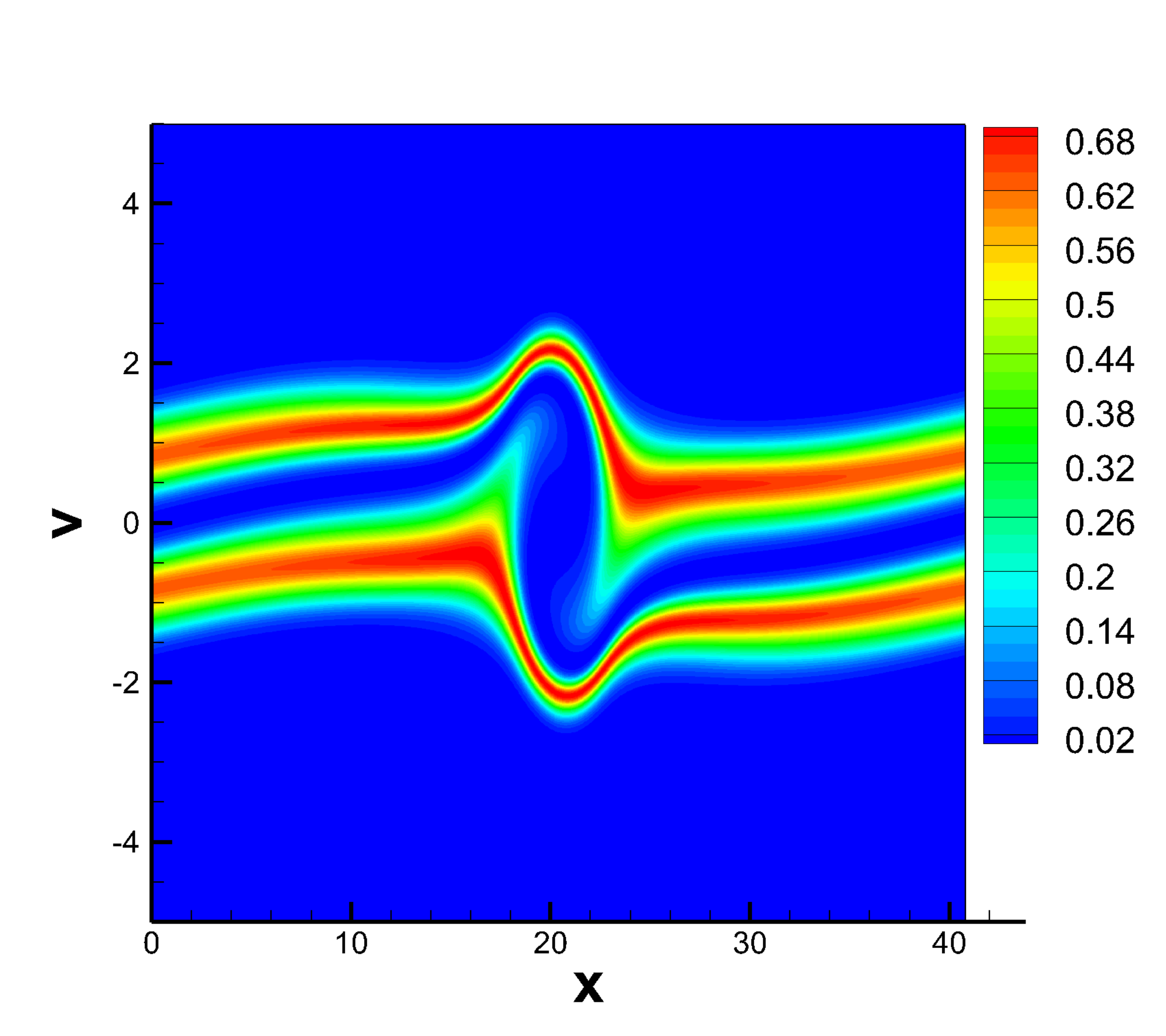}}
				\subfigure[]{\includegraphics[width=.42\textwidth]{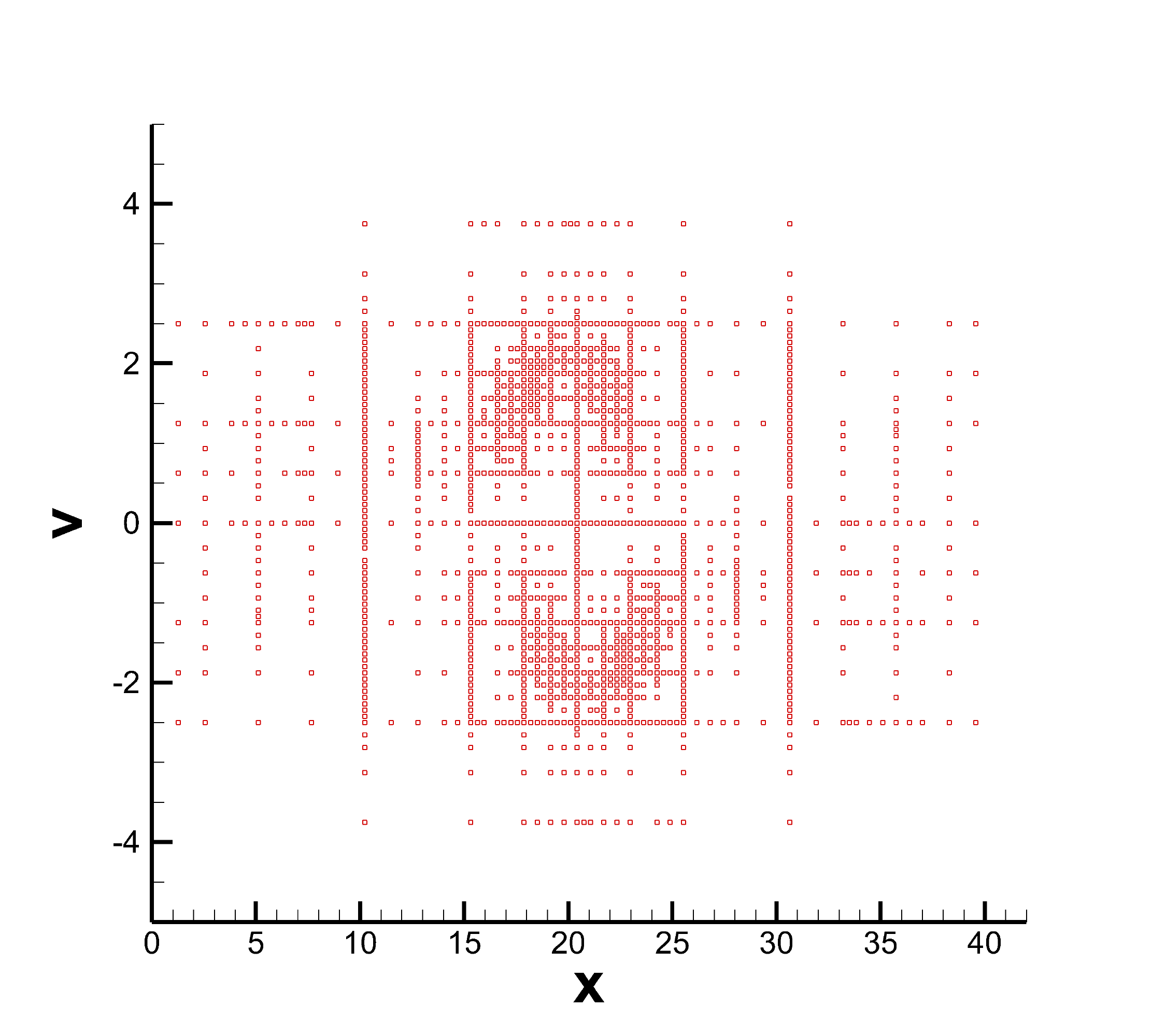}}\\
		\subfigure[]{\includegraphics[width=.42\textwidth]{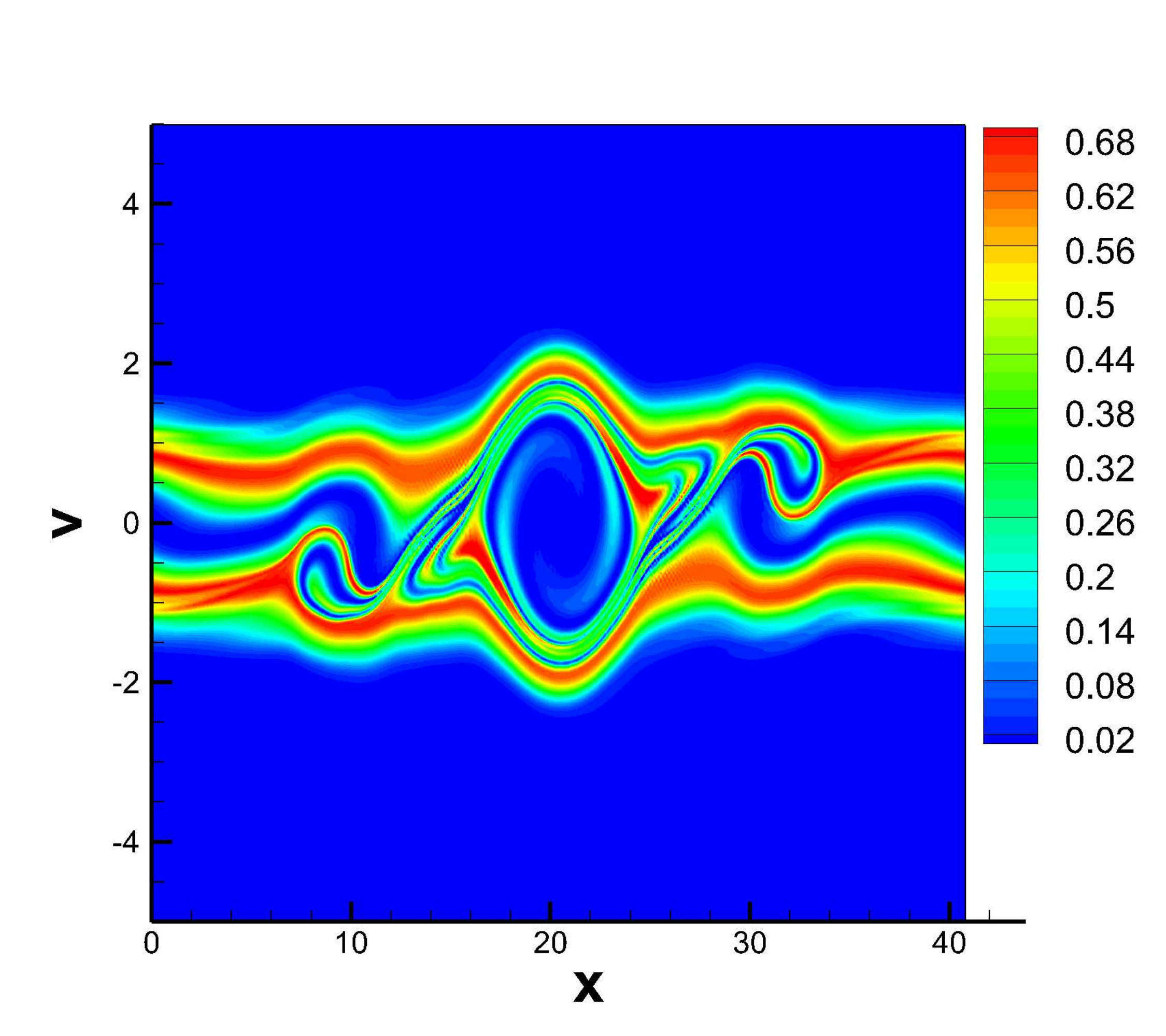}}
		\subfigure[]{\includegraphics[width=.42\textwidth]{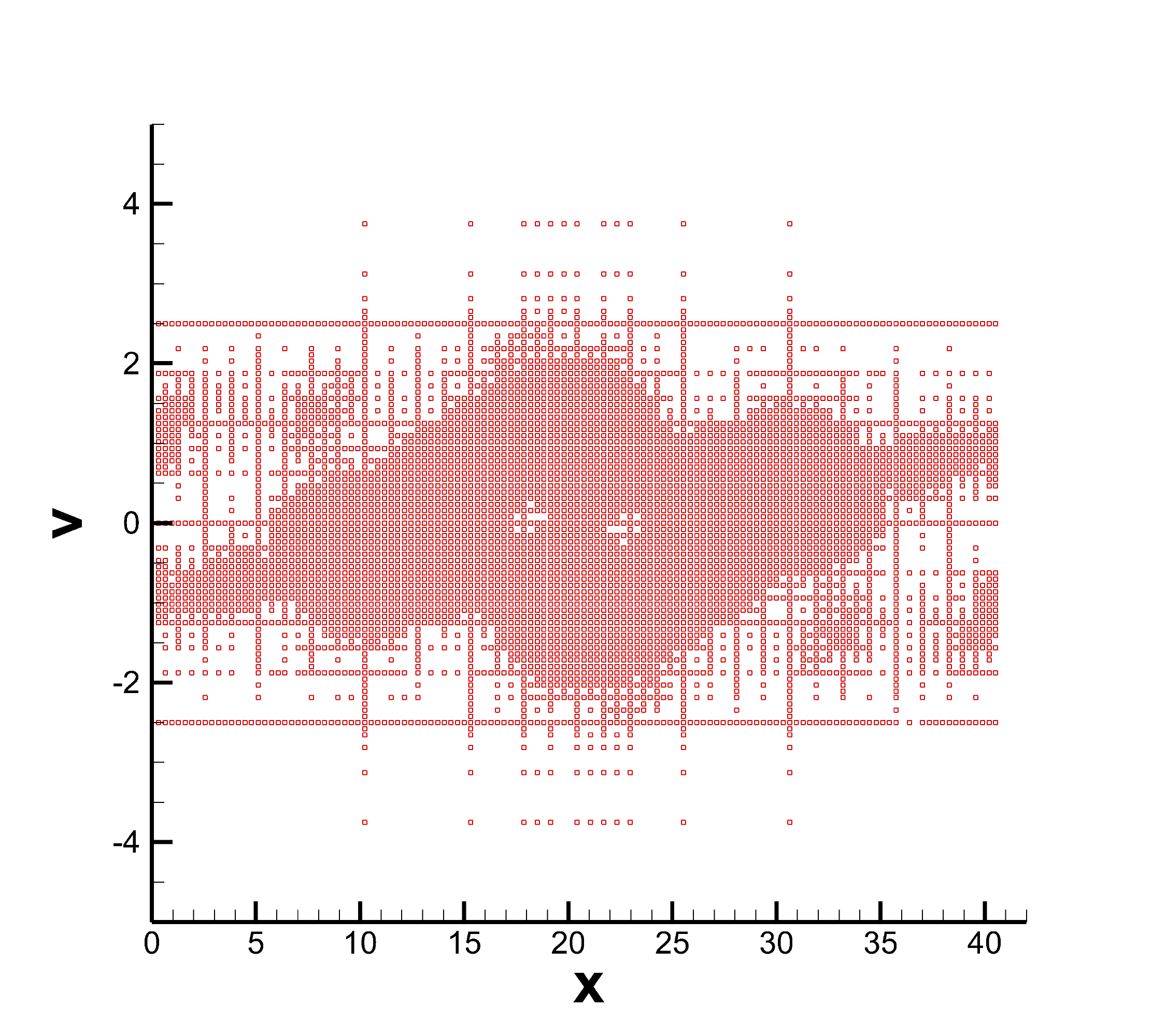}}\\
%		\subfigure[]{\includegraphics[width=.42\textwidth]{figures/con_two2_p3_t70}}
%		\subfigure[]{\includegraphics[width=.42\textwidth]{figures/mesh_two2_p3_t70}}\\
	\end{center}
	\caption{Example \ref{exa:vp}. Two-stream instability II. Phase space contour plots and the associated active elements at   $t=10$ (a-b),  $t=20$ (c-d), $t=40$ (e-f). $N=7$. $k=3$.  $\varepsilon=10^{-5}$.}
	\label{fig:con_two2}
\end{figure}

\begin{figure}[htp]
	\begin{center}
		\subfigure[]{\includegraphics[width=.42\textwidth]{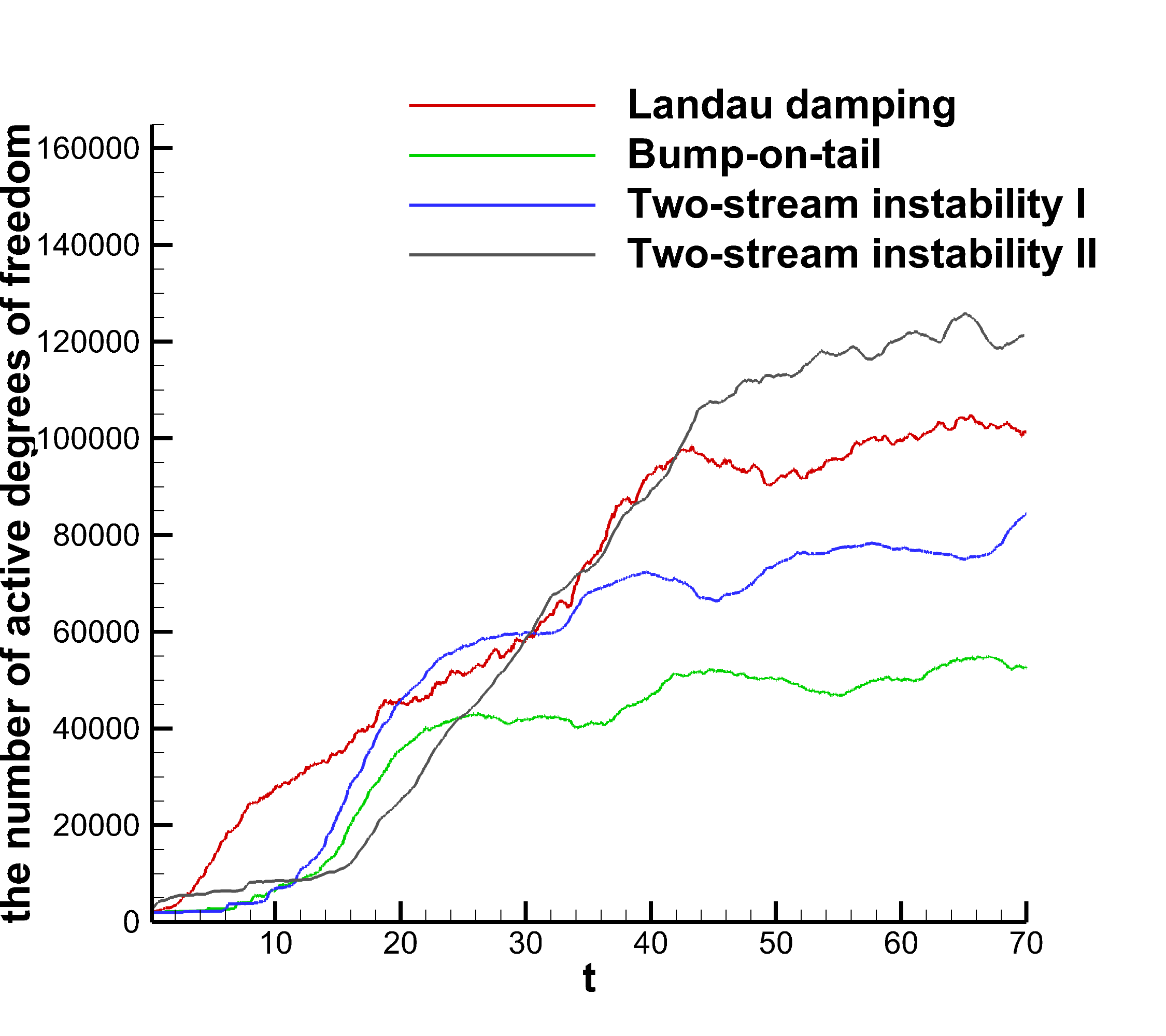}}
	\end{center}
	\caption{Example \ref{exa:vp}. Time histories of the number of active degrees of freedom for Landau damping, bump-on-tail, and two types of two-stream instabilities. $k=3$. $N=7$. $\varepsilon=10^{-5}$.}
	\label{fig:vp_elements}
\end{figure}

In Figures \ref{fig:con_lan}-\ref{fig:con_two2}, we present the phase space contour plots and the associated active elements at several instances of time for all four initial conditions. In Figure \ref{fig:vp_elements}, the time evolution of the numbers of active degrees of freedom are plotted. It is observed that when the solution has not developed rich filamentation structures, only a small amount of degrees of freedom are used. As time evolves, thiner and thiner  filaments are generated because of  phase mixing. The adaptive scheme can automatically add degrees of freedom   to adequately resolve the fine structures. We remark that the quality of the numerical results are quite comparable to those computed by the more expensive full grid DG method with similar mesh size, e.g., see \cite{cheng_vp}, while much less degrees of freedom are needed, leading to computational savings.

%height=.30\textwidth

\begin{figure}[htp]
	\begin{center}
		\subfigure[]{\includegraphics[width=.42\textwidth]{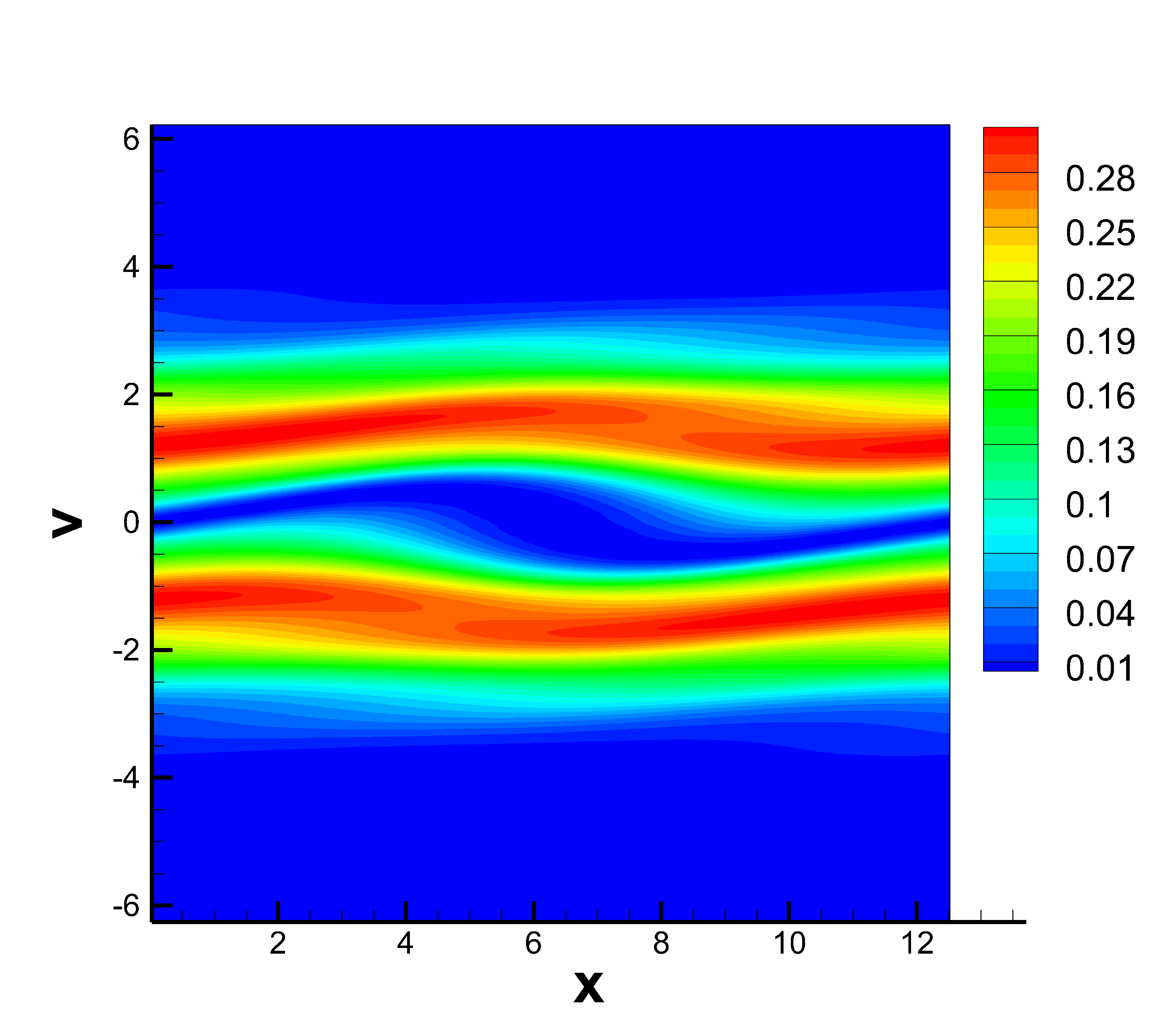}}
			\subfigure[]{\includegraphics[width=.42\textwidth]{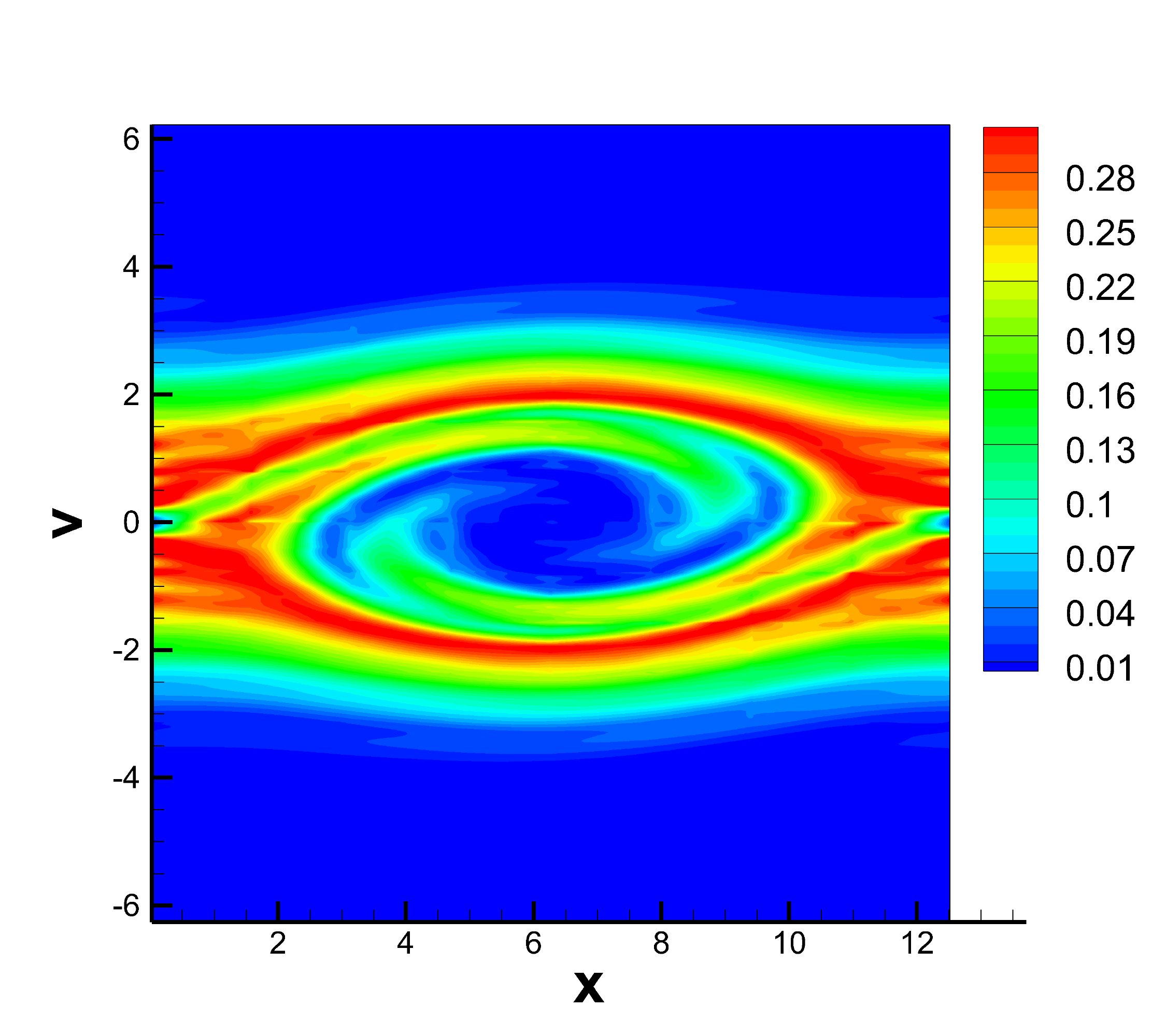}}
	\end{center}
	\caption{Example \ref{exa:vp}. Two-stream instability I. Phase space contour plots of the numerical solution by the sparse grid DG scheme \cite{guo_sparsedg}. $N=7$.  $k=3$. (a) $t=10$. (b) $t=20$. }
	\label{fig:vp_sparse}
\end{figure}

\begin{figure}[htp]
	\begin{center}
		\subfigure[]{\includegraphics[width=.42\textwidth]{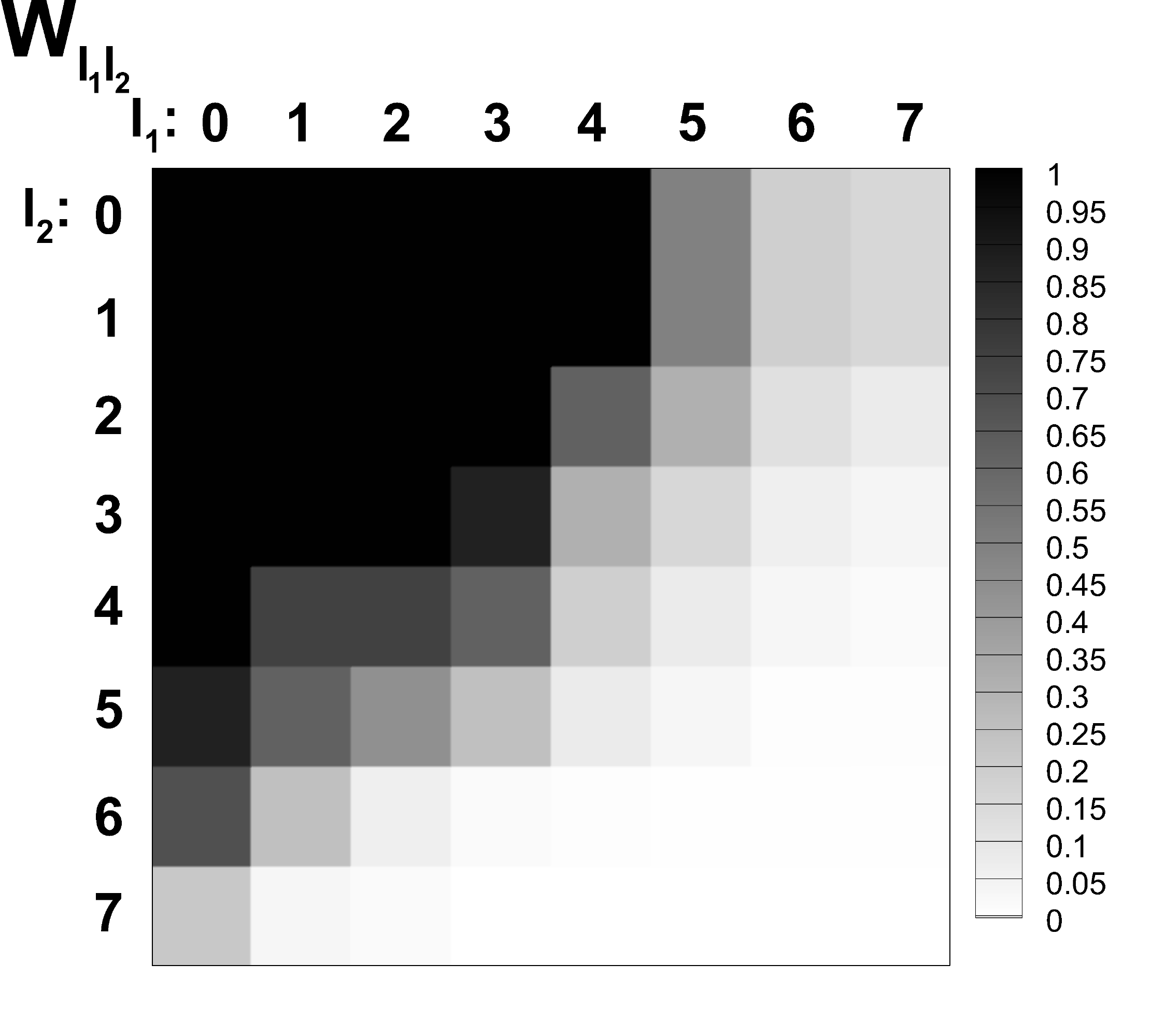}}
		\subfigure[]{\includegraphics[width=.42\textwidth]{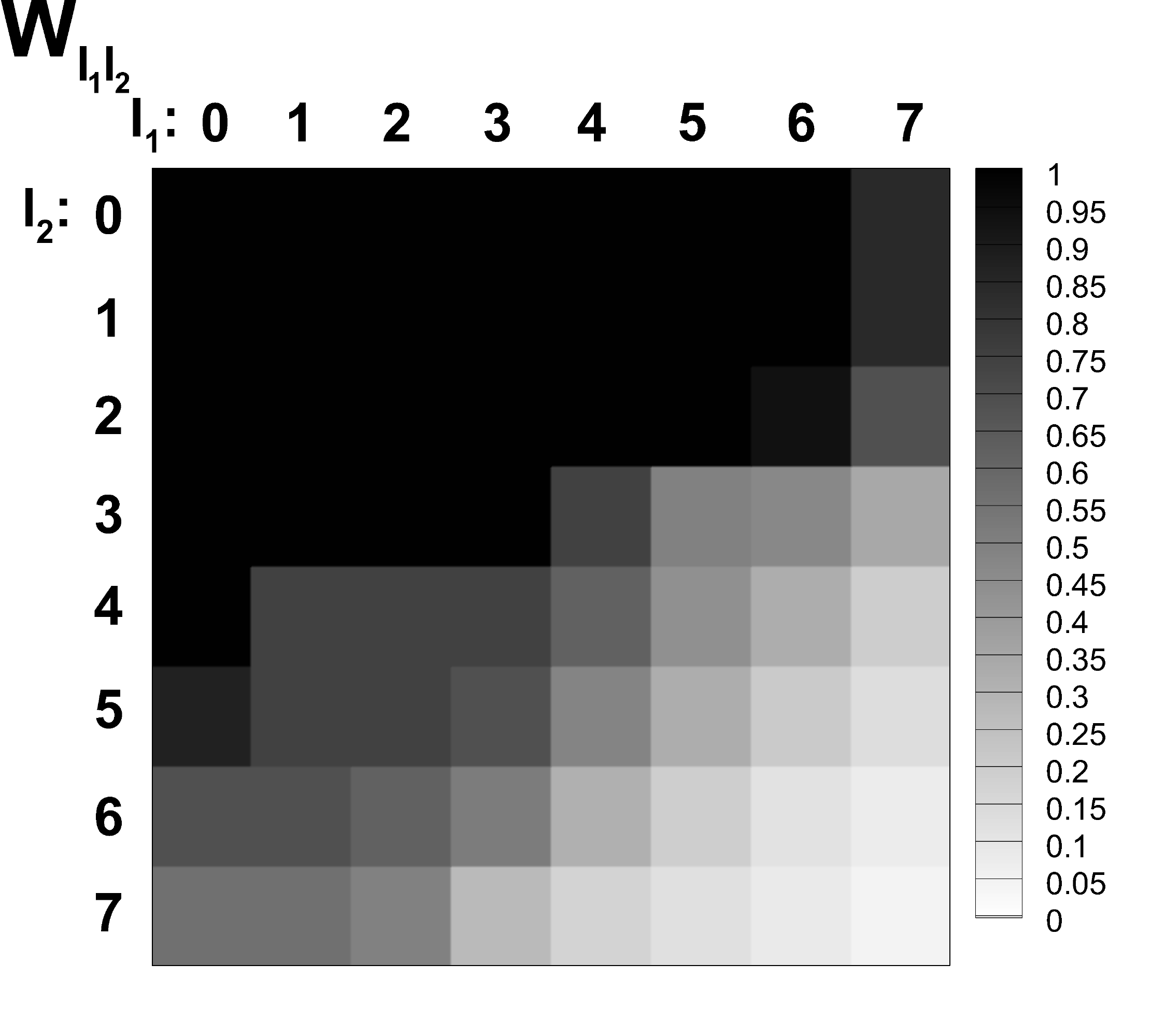}}
	\end{center}
	\caption{Example \ref{exa:vp}. Two-stream instability I. The percentage of used elements for each incremental space $\bW_\bl$, $\bl=\{l_1,l_2\}$ and $|\bl|_\infty\leq N$. $N=7$. $k=3$. $\varepsilon=10^{-5}$. (a) $t=10$. (b) $t=20$. }
	\label{fig:vp_density}
\end{figure}

Lastly, we present a numerical comparison between the sparse grid DG method \cite{guo_sparsedg} and the adaptive scheme in this paper. We consider 
 two-stream instability I  with  parameter choice $N=7$, $k=3$ for both methods solving the VP system. %computed by the sparse grid DG method with approximation space $\hat{V}_7^3$ to demonstrate the advantages of the proposed adaptive MRA DG method. 
The phase space contour plots of the sparse grid DG method at $t=10$ when the solution is very smooth and at $t=20$ when the solution has developed filamentations are provided in Figure \ref{fig:vp_sparse}, to be compared with the results of the adaptive scheme in Figure \ref{fig:con_two}.   In Figure \ref{fig:vp_density}, we plot the percentage of used elements for each incremental space $\bW_\bl$ at at time $t=10$ and $20$ for the adaptive method. While both schemes actually provide similar accurate description of  the macroscopic moments, there is a qualitative difference when the numerical resolution of $f$ is concerned. As expected, when the solution is smooth at $t=10$, both the sparse grid DG method and the adaptive method can generate reliable results with comparable degrees of freedom, see Figure \ref{fig:con_two}(a) and Figure \ref{fig:vp_density}(a). At $t=20$, the sparse grid DG method
 does not resolve all the fine structures when compared with the adaptive  method (see Figure \ref{fig:con_two}(c) versus Figure \ref{fig:vp_sparse}(b)). At this time, the adaptive method uses more degrees of freedom than the sparse grid method (see Figure \ref{fig:vp_density}(b)), but  the DOFs are still much less than the full grid method. %However, it does provide comparable results
 
% In comparison, when the solution is smooth, the adaptive scheme is also capable of generating good results with comparable  degrees of freedom to the sparse DG method; while, the adaptive scheme is able to resolve localized structures by automatically adding more degrees of freedom where filamentations occur, and provide similar results to those by the full grid counterpart \cite{cheng2014energy}, see Figure \ref{fig:con_two}. In Figure \ref{fig:vp_density}, we plot the percentage of used elements for each incremental space $\bW_\bl$ at at time $t=10$ and $20$. Again, It is observed the pattern looks similar to a sparse grid, i.e., the upper part of incremental spaces are mainly used in the simulation, when the solution is smooth ($t=10$). Furthermore, when filamentations occur in the solution ($t=20$), the adaptive scheme takes more degrees of freedom from the lower part of incremental spaces in order to resolve such localized structures.

\begin{exa}[Oscillatory VP system]
\label{exa:ovp}
We consider the following oscillatory VP system in the polar coordinates:
\begin{eqnarray}
&&f_t + \frac{v}{\epsilon} f_r +  (E(t,r)+E_{ext}(t,r)) f_v=0, \label{eq:osc_v}\\
&&\partial_r(rE(t,r)) = r\rho(t,r), \label{eq:ocs_poisson}
\end{eqnarray}
where the dimensionless parameter  $\epsilon=0.05$ denotes the ratio between the characteristic lengths
in the transverse and the longitudinal directions \cite{crouseilles2016}. $E_{ext}$ is the external electric field specified as 
$$E_{ext}(t,r) = -\frac{r}{\epsilon}+r\cos^2\left(\frac{t}{\epsilon}\right).$$
The initial condition is set to be a discontinuous function
\begin{equation}\label{eq:beam}
f(0,r,v)=\frac{n_0}{v_t\sqrt{2\pi}}\exp\left(-\frac{v^2}{2v_t^2}\right)\chi_{[-r_m,r_m]}(r),\quad (r,v)\in[-3,3]^2,
\end{equation}
where $n_0=4$,\, $v_t=0.1,\, r_m=1.85,$ and 
\begin{equation*}
\chi_{[-r_m,r_m]}(r) =\left\{\begin{array}{ll} 1& \text{if}\quad -r_m\leq r\leq r_m,\\
0&\text{otherwise}.\end{array}
 \right.
\end{equation*}
\end{exa}

This example has been intensively studied in \cite{crouseilles2013asymptotic,Frenod2015169,crouseilles2016}, where several effective schemes have been developed. 
Note that the initial condition \eqref{eq:beam} considered here is discontinuous and represents a semi-Gaussian beam in particle accelerator physics \cite{crouseilles2016}. We impose zero boundary conditions in both $r$ and $v$ directions. An LDG method is used to solve Poisson's equation \eqref{eq:ocs_poisson} and the closure condition $E(0)=0$ is strongly imposed in the formulation. In the simulation, we let $\epsilon=0.05$, and  $k=2$, $N=7$, $\varepsilon=10^{-4}$. We consider both $L^1$ (\eqref{eq:l1} \eqref{eq:l1_c}) and $L^2$ (\eqref{eq:l2} \eqref{eq:l2_c}) norm based criteria as the refining and coarsening indicators.

\begin{figure}[htp]
	\begin{center}
		\subfigure[]{\includegraphics[width=.42\textwidth]{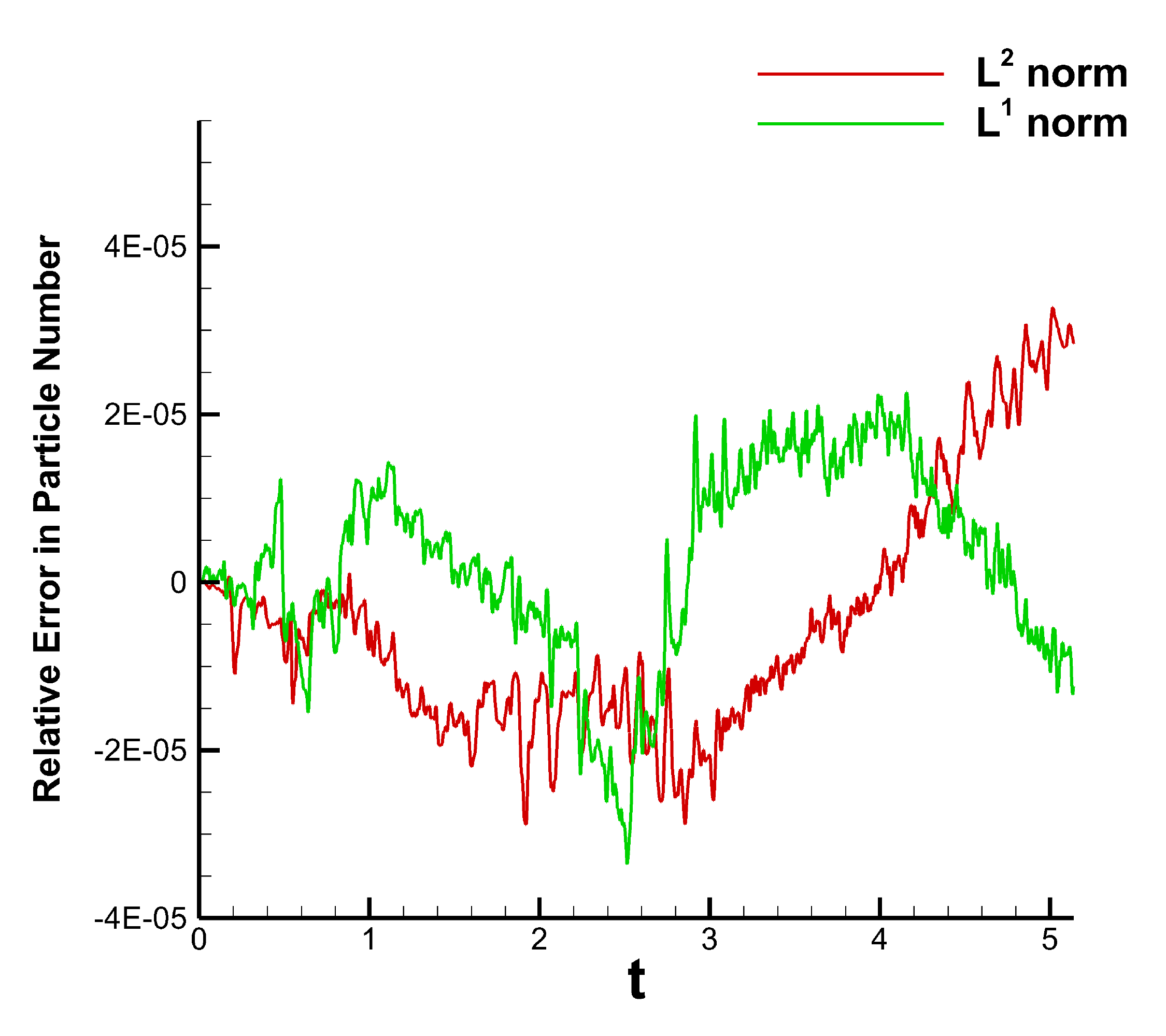}}
		\subfigure[]{\includegraphics[width=.42\textwidth]{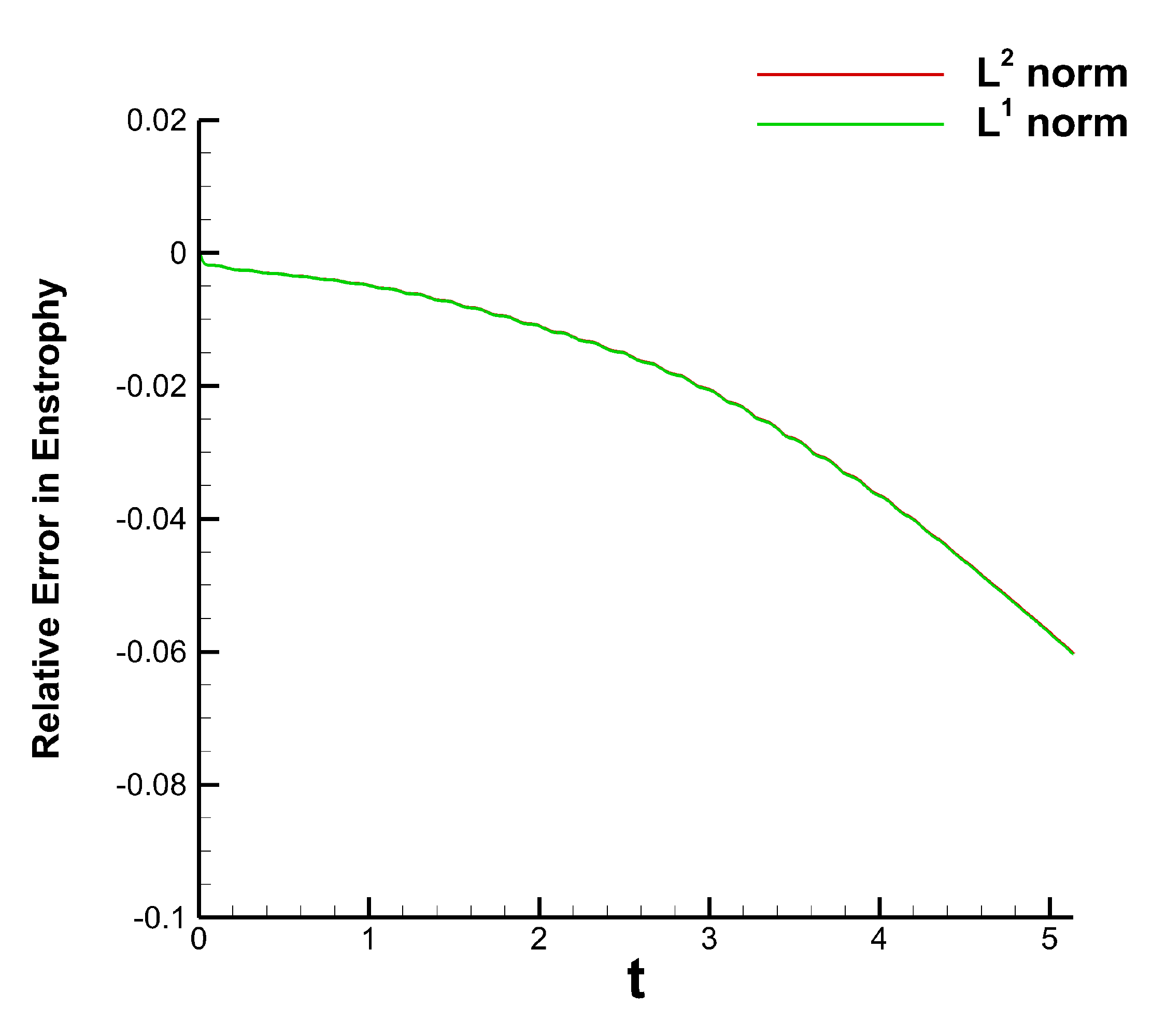}}
	\end{center}
	\caption{Example \ref{exa:ovp}.  Evolution of the relative errors in total particle number (a), enstrophy (b). $N=7$. $k=2$. $\varepsilon=10^{-4}$. }
	\label{fig:evo_vp_conser}
\end{figure}

We first present the time evolution of the relative errors in total particle number and enstrophy in Figure \ref{fig:evo_vp_conser}. Similar to the previous VP system, the scheme with both adaptive indicators is able to conserve the particle number up to the magnitude of $\varepsilon$. The enstrophy decays due to the choice of the numerical flux. 

\begin{figure}[htp]
	\begin{center}
%		\subfigure[]{\includegraphics[width=.42\textwidth]{figures/polar_0_l1}}
%		\subfigure[]{\includegraphics[width=.42\textwidth]{figures/polar_mesh_0_l1}}\\
		\subfigure[]{\includegraphics[width=.42\textwidth]{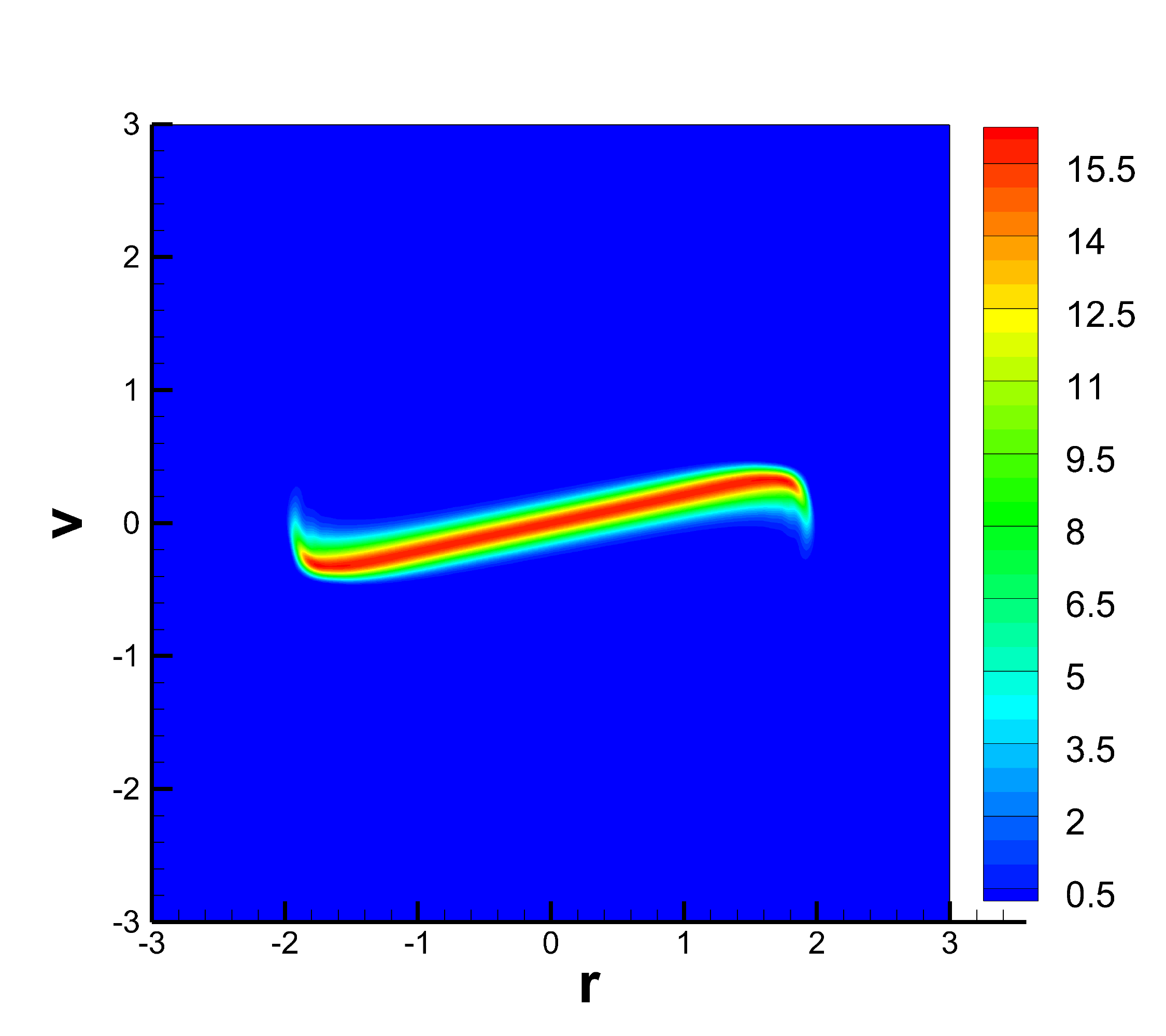}}
		\subfigure[]{\includegraphics[width=.42\textwidth]{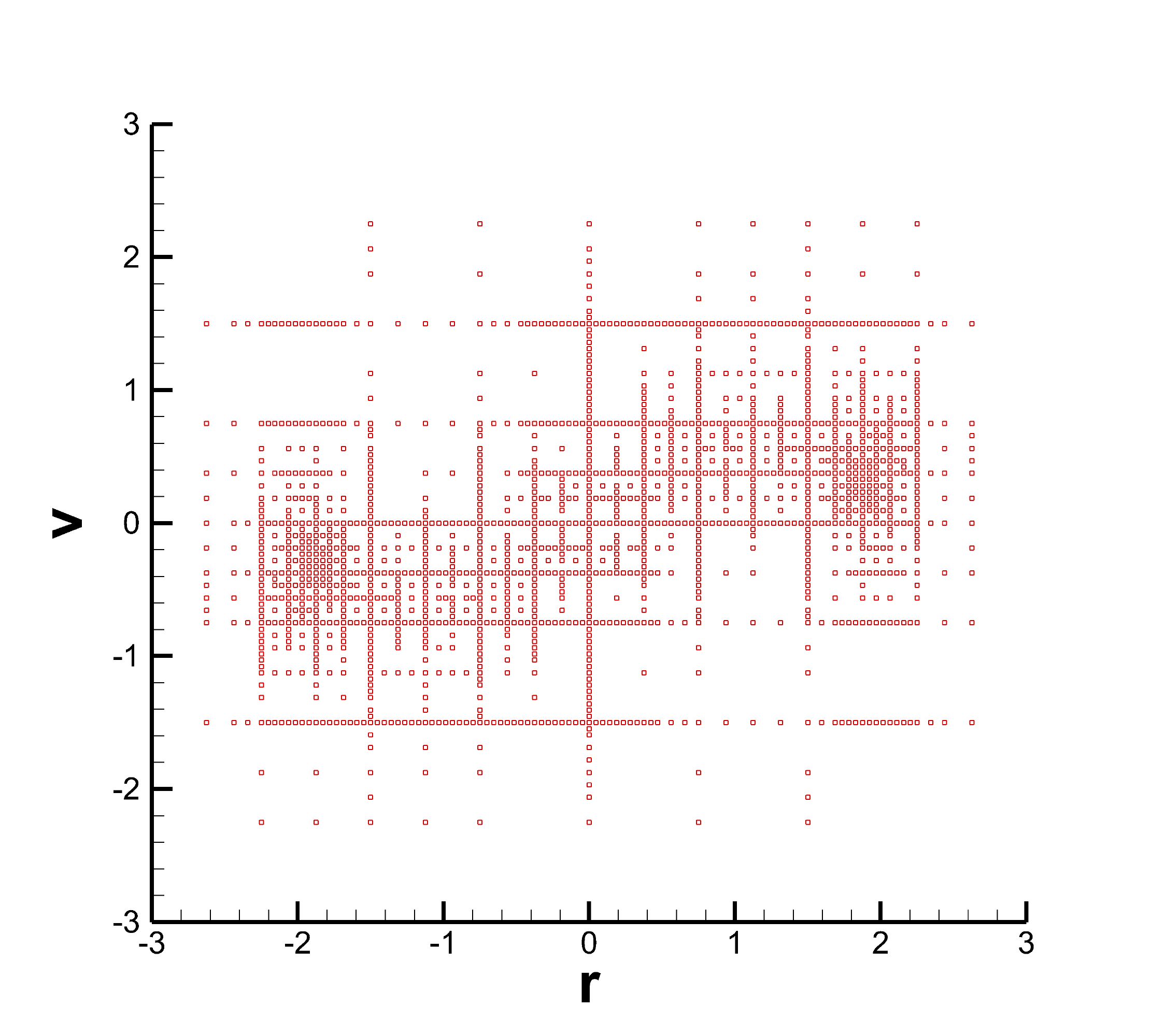}}\\
		\subfigure[]{\includegraphics[width=.42\textwidth]{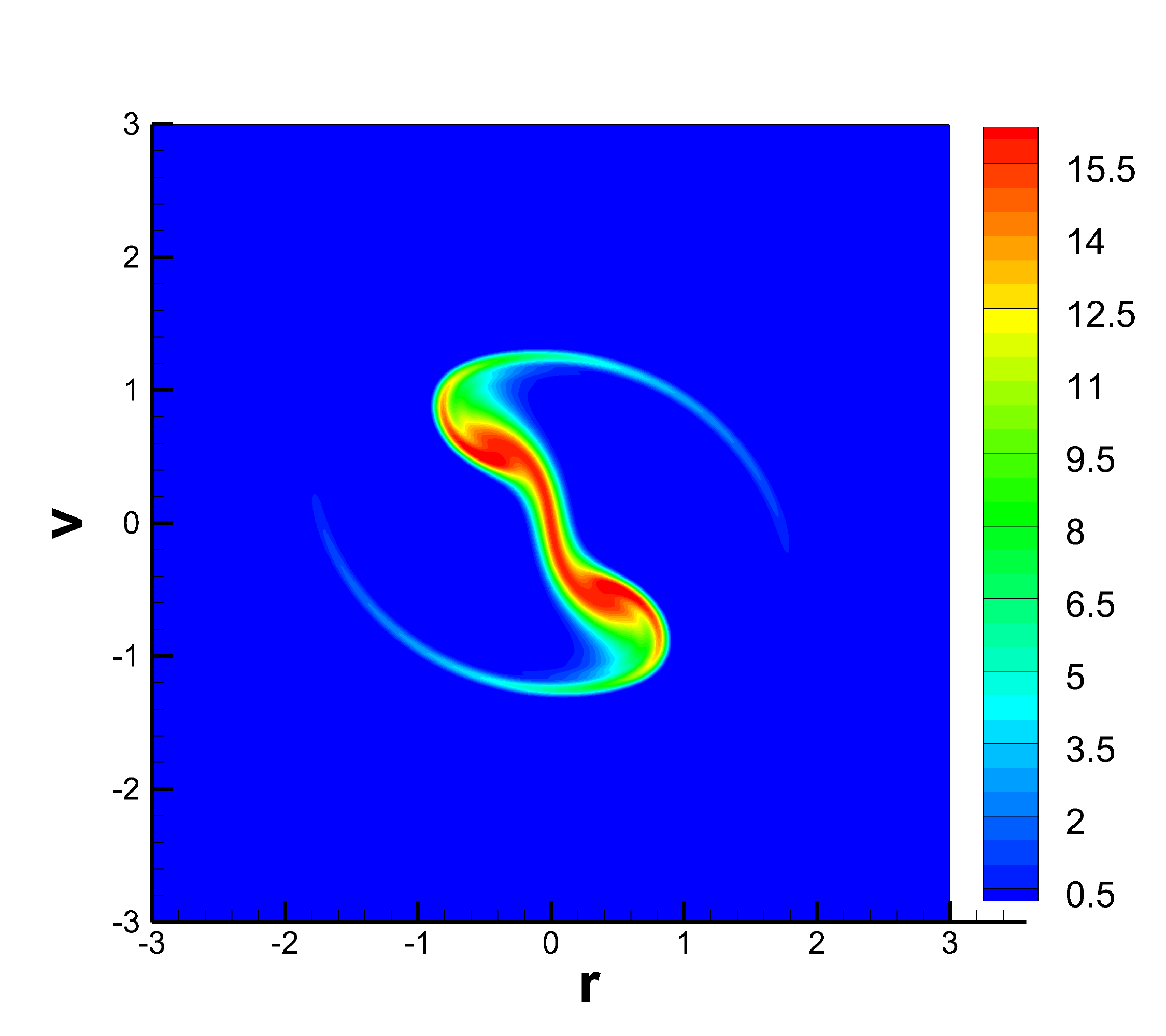}}
		\subfigure[]{\includegraphics[width=.42\textwidth]{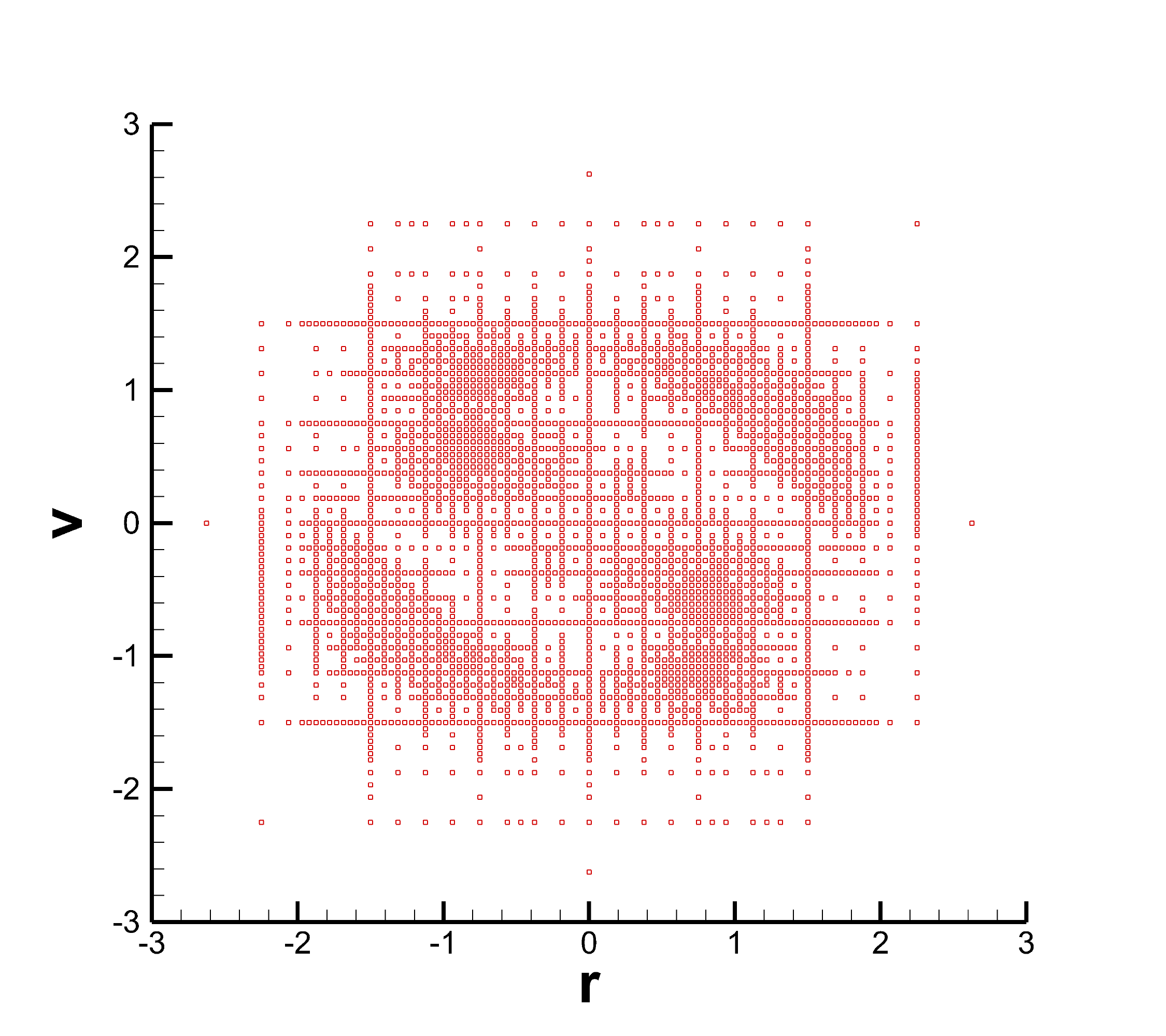}}\\
		\subfigure[]{\includegraphics[width=.42\textwidth]{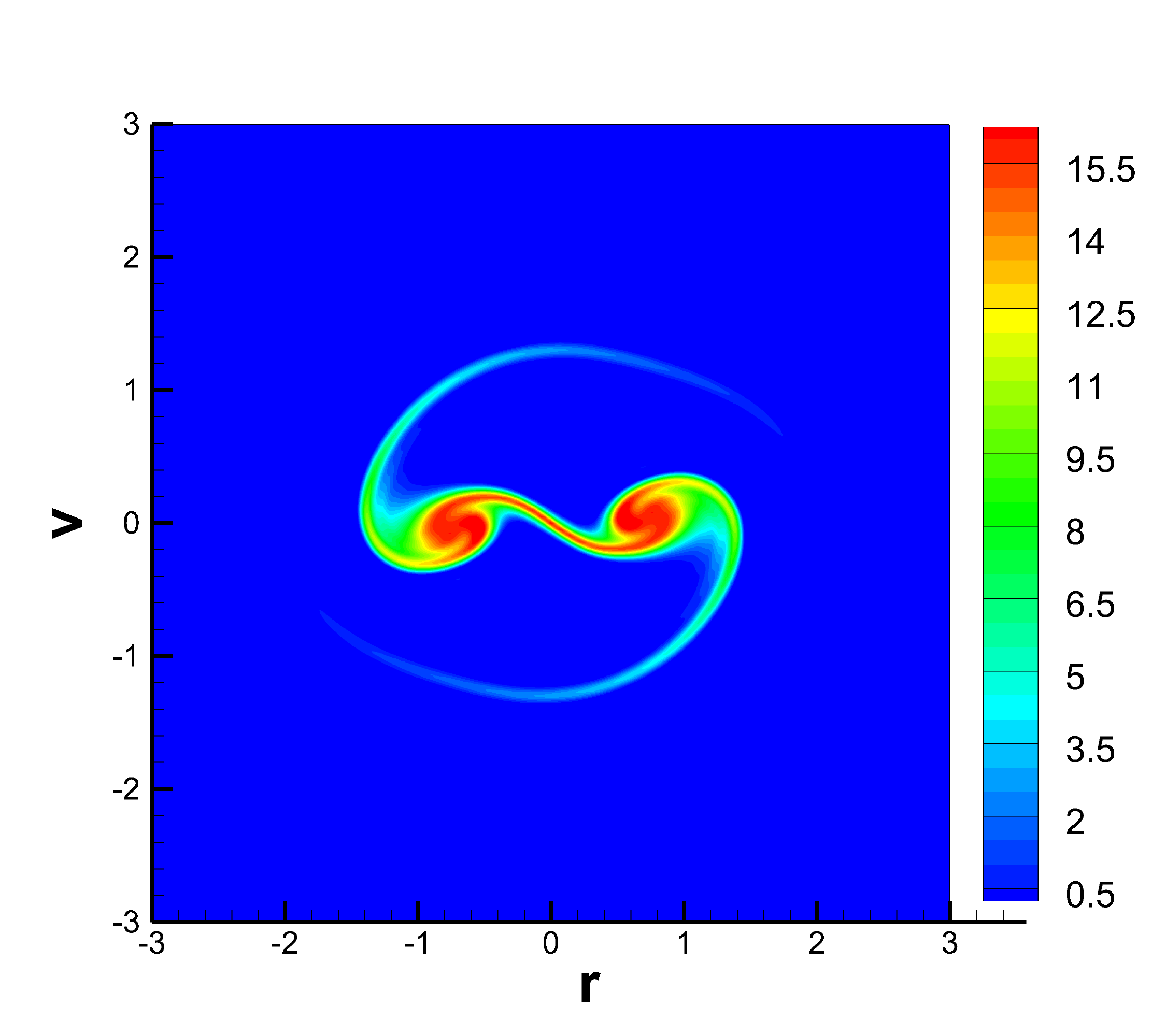}}
		\subfigure[]{\includegraphics[width=.42\textwidth]{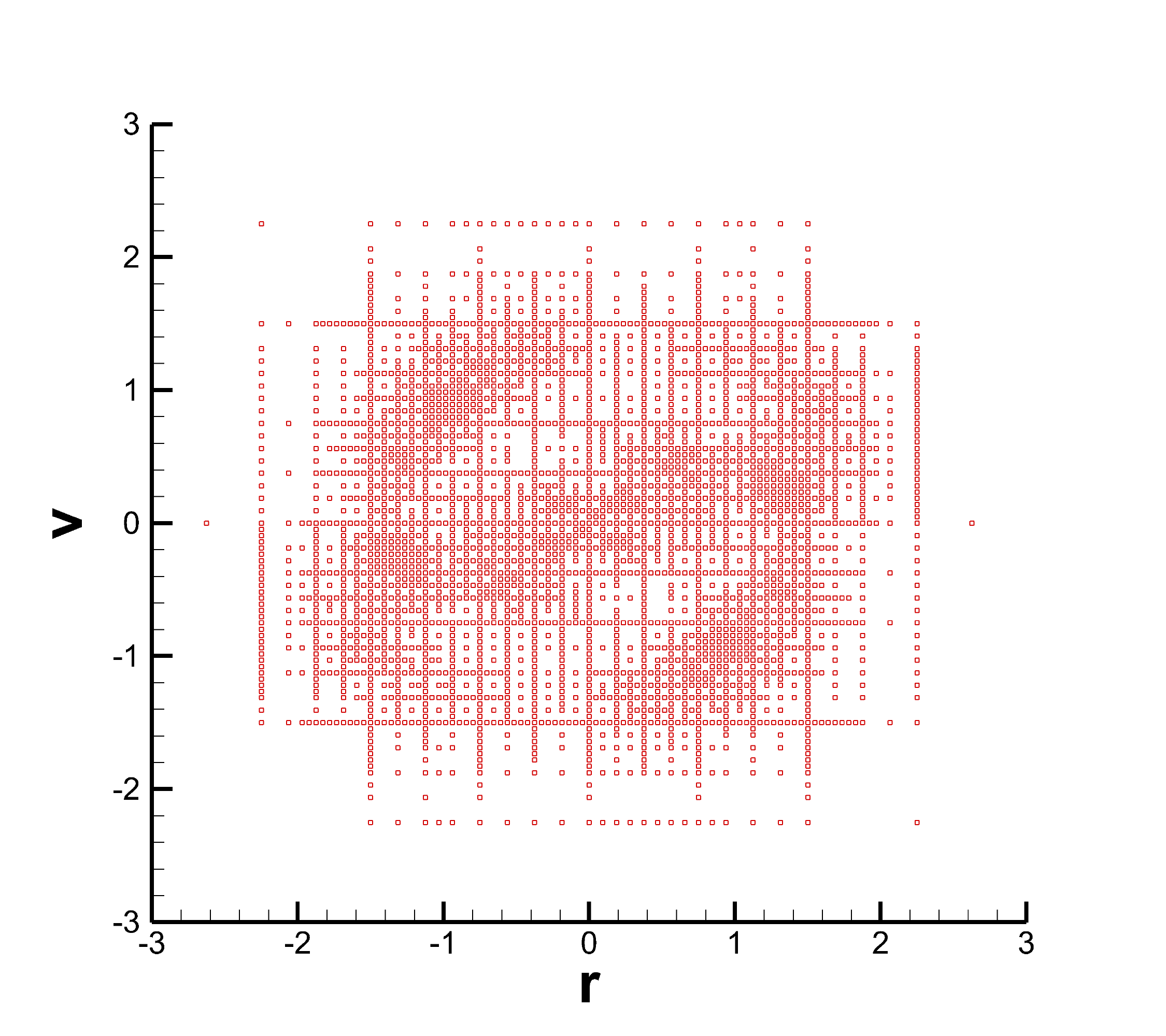}}\\		
	\end{center}
	\caption{Example \ref{exa:ovp}.    Phase space contour plots and the associated active elements at  $t=1.3464$ (a-b), $t=4.3388
	$ (c-d), $t=5.1462$ (e-f). $N=7$. $k=2$. $\varepsilon=10^{-4}$.  The $L^1$ norm based criteria are used as the adaptive indicator.}
	\label{fig:con_osc_l1}
\end{figure}

\begin{figure}[htp]
	\begin{center}
		\subfigure[]{\includegraphics[width=.42\textwidth]{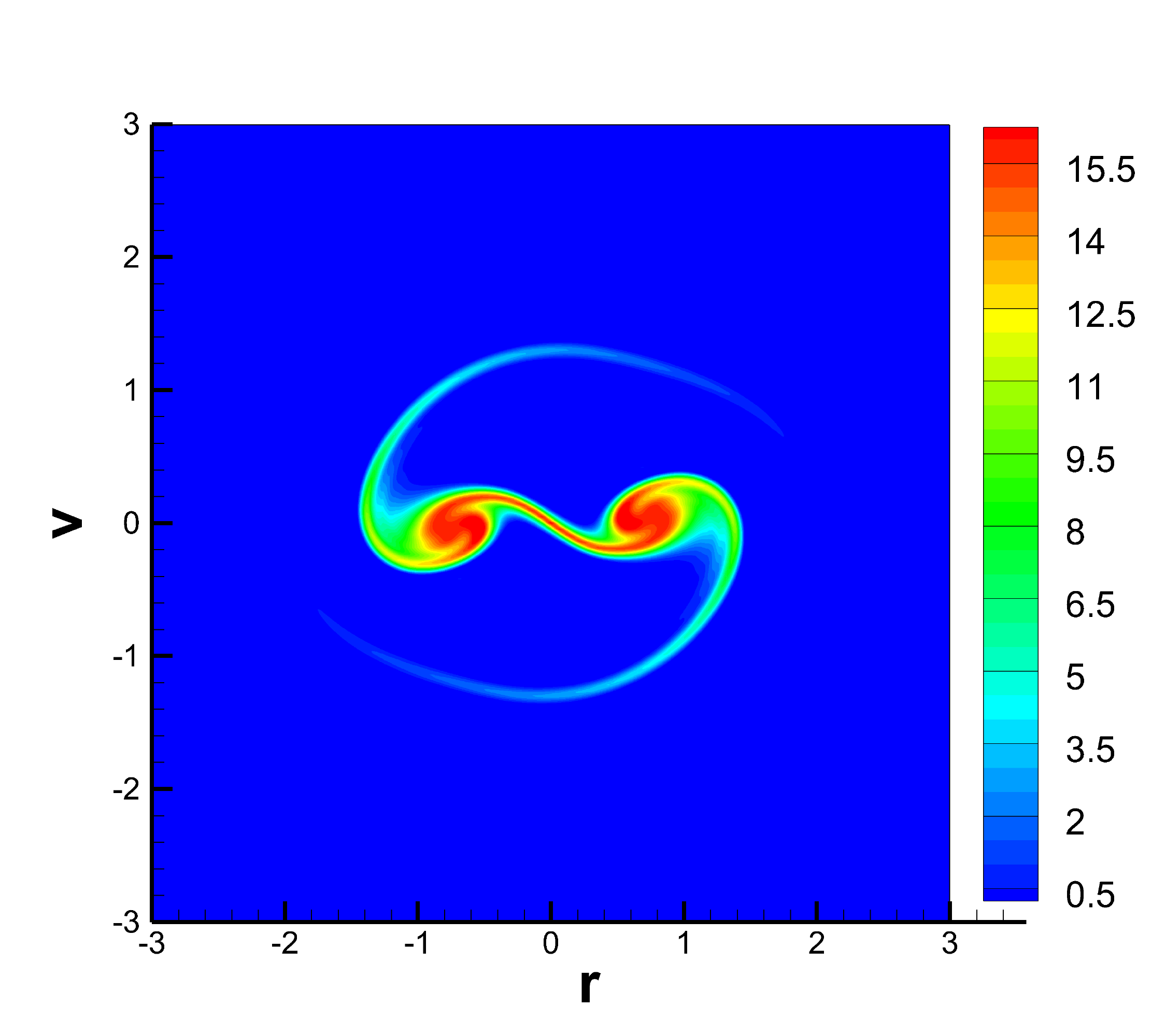}}
		\subfigure[]{\includegraphics[width=.42\textwidth]{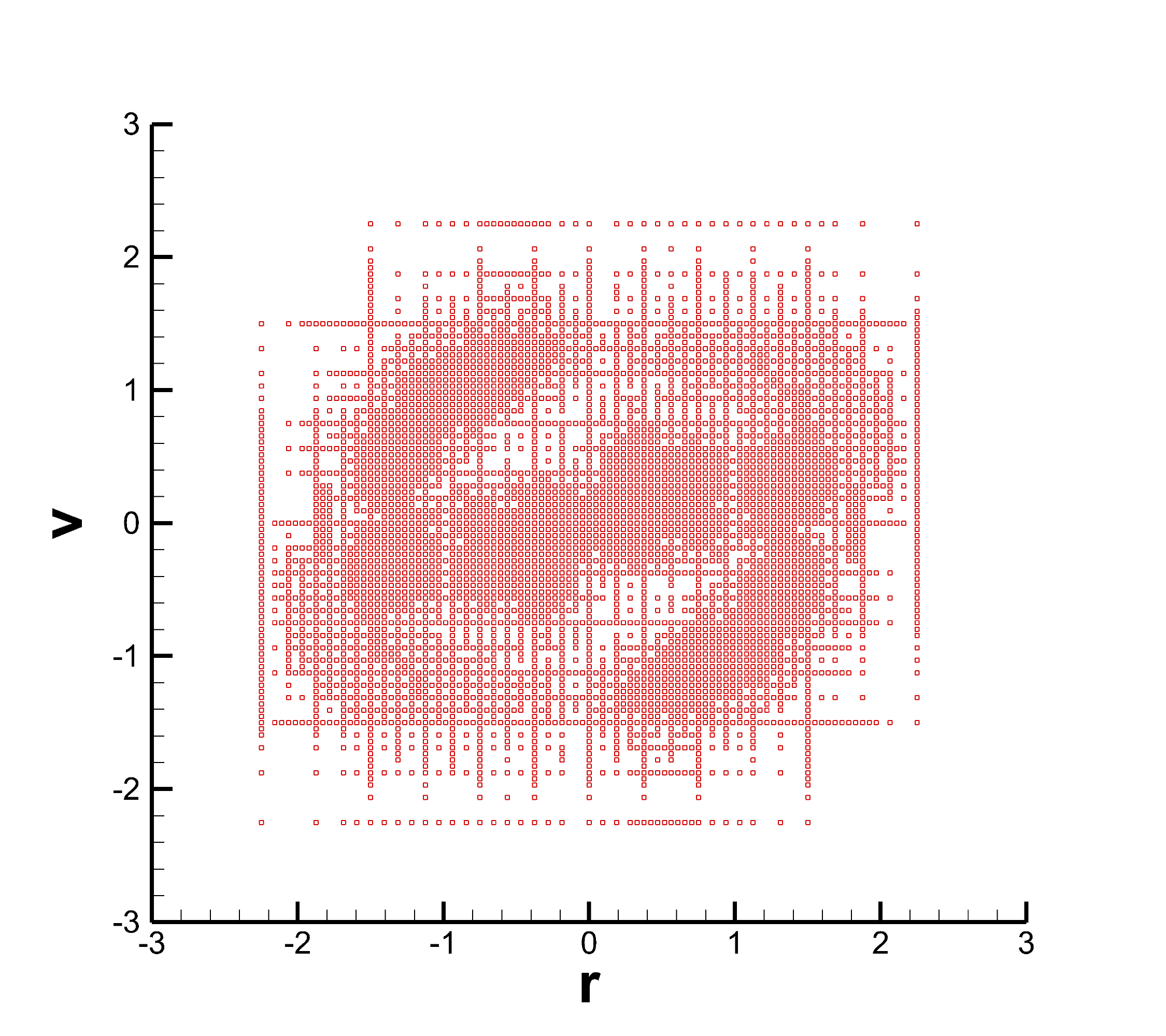}}	
	\end{center}
	\caption{Example \ref{exa:ovp}.  Phase space contour plots and the associated active elements at $t=5.1462$. $k=2$. $N=7$. $\varepsilon=10^{-4}$.  The $L^2$ norm based criteria are used as the adaptive indicator.}
	\label{fig:con_osc_l2}
\end{figure}

The phase space contours from our scheme agrees well with those in the literature \cite{crouseilles2016}.
In Figure \ref{fig:con_osc_l1}, we present the   contour plots and the associated active elements at three instances of time, for which the $L^1$ norm based criteria are used as the adaptive indicator. In Figure \ref{fig:con_osc_l2}, we also report the contour plot and associated adaptive mesh at final time with the $L^2$ norm based criteria to compare the performance of the two criteria as adaptive indicators. It is observed that the numerical results are qualitatively the same, but more elements are used by the scheme with the $L^2$ norm based criteria. In Figure \ref{fig:evo_vp_conser}, the time evolution of the number of active degrees of freedom are plotted. Again, when the solution develops   filaments, more degrees of freedom are added thanks to the adaptive mechanism.   In summary, for this example, the $L^1$ norm based criteria is preferred for the sake of efficiency.% Also note that the scheme with the $L^1$ norm based criteria employs less degrees of freedom than that with the $L^2$ norm based criteria, yet the quality of numerical results are almost the same. Hence, for this example with discontinuous solution structures, the $L^1$ norm based criteria is preferred for the sake of efficiency. %Lastly, we remark that the numerical results agree well with those reported in \cite{crouseilles2016}. 

\begin{figure}[htp]
	\begin{center}
		\subfigure[]{\includegraphics[width=.42\textwidth]{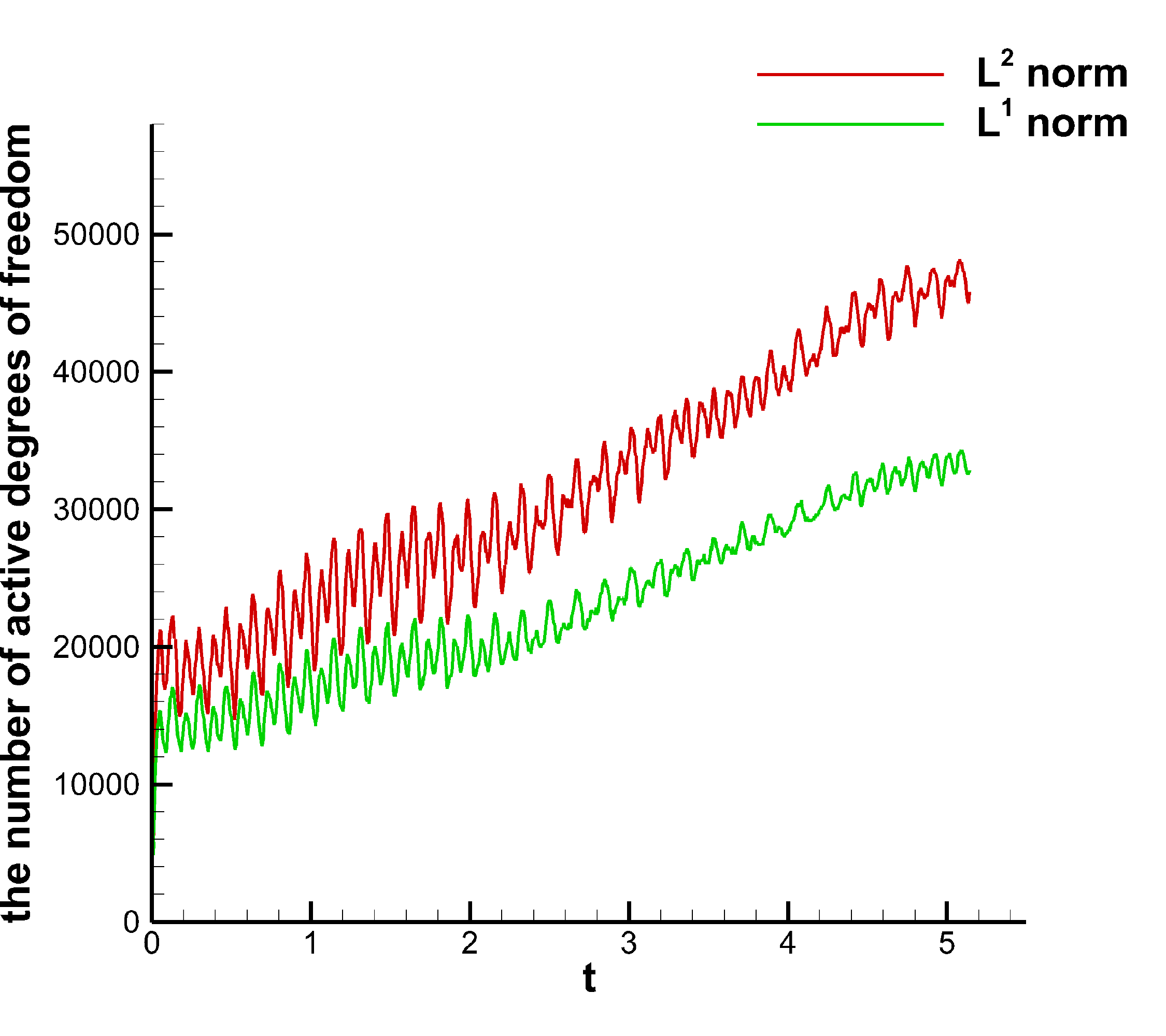}}
	\end{center}
	\caption{Example \ref{exa:ovp}.  Time histories of the number of active degrees of freedom. $N=7$. $k=2$.  $\varepsilon=10^{-4}$.}
	\label{fig:vp_osc_elements}
\end{figure}

\section{Conclusions and future work}
\label{sec:conclusion}
In this paper, we develop an adaptive multiresolution DG scheme for computing time-dependent transport equations. The key ingredients of the scheme are the weak formulation of the DG method and adaptive error thresholding based on hierarchical surplus. Extensive numerical tests show that our scheme performs similarly to a sparse grid DG method when the solution is smooth, and can automatically capture fine local structures when the solution is no longer smooth. Detailed comparison between several refinement/coarsening error indicators are performed. The method is demonstrated to work well for kinetic simulations. Future work consists of the study of limiters and further improvement of the scheme including local time stepping and adaptivity with both the mesh and polynomial degrees.

\bibliographystyle{abbrv}
\bibliography{ref_cheng,ref_cheng_2,adaptive}

\end{document}